\input amstex
\input xy
\xyoption{all}
\SelectTips {cm}{}

\baselineskip=24pt

\NoBlackBoxes

\documentstyle{amsppt}
\pagewidth{35pc}

\topmatter
\title
Concrete Duality for Strict $\infty $-Categories
\endtitle
\author
G.V. Kondratiev
\endauthor
\address
\endaddress
\email
gennadii\@\, hotmail.com
\endemail
\subjclass
CT Category Theory
\endsubjclass
\abstract
An elementary theory of strict $\infty $-categories with application to concrete duality is given. New examples of first and second order
concrete duality are presented.
\endabstract

\endtopmatter

\document

\vskip 0.35cm

\item{} \hskip 5.5cm {\bf \hskip 0.8cm Content} 
\vskip 0.2cm
\item{} {\bf 0. Introduction} \hfill 2 \ \ \ \ 
\item{} {\bf 1. Categories, functors, natural transformations, modifications} \hfill 5 \ \ \ \ 
\item{} {\bf 1a. Weak categories, functors, natural transformations, modifications} \hfill 11 \ \ \ \ 
\item{} {\bf 1.1. Fractal organization of the new universe}  \hfill 16 \ \ \ \ 
\item{} {\bf 1.2. Notes on Coherence Principle}               \hfill 16 \ \ \ \ 
\item{} {\bf 2. $(m,n)$-invariants}                          \hfill 16 \ \ \ \ 
\item{} {\bf 2.1. Homotopy groups associated to $\infty $-categories} \hfill 17 \ \ \ \ 
\item{} {\bf 2.2. Duality and Invariant Theory}           \hfill 19 \ \ \ \ 
\item{} {\bf 2.2.1. Classification of covariant geometric objects} \hfill 19 \ \ \ \ 
\item{} {\bf 2.2.2. Classification of smooth embeddings into Lie group}    \hfill 21 \ \ \ \ 
\item{} {\bf 2.3. Tangent functor for smooth algebras}   \hfill 21 \ \ \ \ 
\item{} {\bf 3. Representable $\infty $-functors}   \hfill 22 \ \ \ \ 
\item{} {\bf 4. (Co)limits}    \hfill 26 \ \ \ \ 
\item{} {\bf 5. Adjunction}    \hfill 28 \ \ \ \ 
\item{} {\bf 6. Concrete duality for $\infty$-categories} \hfill 31 \ \ \ \ 
\item{} {\bf 6.1. Natural and non natural duality} \hfill 32 \ \ \ \ 
\item{} {\bf 7. Vinogradov duality}  \hfill 36 \ \ \ \ 
\item{} {\bf 8. Duality for differential equations}  \hfill 37 \ \ \ \ 
\item{} {\bf 8.1. Cartan involution}     \hfill 39 \ \ \ \ 
\item{} {\bf 9. Gelfand-Naimark 2-duality}  \hfill 40 \ \ \ \ 
\item{} {\bf 10. Lukacs' extension of Pontryagin duality}   \hfill 45 \ \ \ \ 
\item{} {\bf 11. Differential algebras as a dual to Lie calculus} \hfill 46 \ \ \ \ 
\item{} {\bf Bibliography}   \hfill 48 \ \ \ \

\vskip 0.5cm

\head {} {\bf 0. Introduction}\endhead 

\vskip 0.5cm
\hskip 9.1cm\hbox{\it "One should jump over not just}

\hskip 8.6cm\hbox{\it theorems but whole theories as well"}

\hskip 12cm\hbox{V.\hskip 0.02cmI. Arnold}

\vskip 1cm

This work is an analysis of the concept duality being used in modern mathematics. 
The main theorem 6.1.1 is a generalization for strict $\infty $-categories of Porst-Tholen 
criterion of existence of a concrete natural dual adjunction for usual 1-categories.

\head {\bf 0.1. Why duality is effective} \endhead

Duality is one of the fundamental recurring ideas in all of mathematics, and category theory provides the appropriate framework for defining and analyzing the idea that ``opposite structures can reflect one another". Some of the most famous theorems in mathematics are duality theorems. We are thinking in particular of Pontryagin duality, Gelfand -Naimark duality and Stone duality.

The computational power of such dualities increases since one can choose that side which works in the simplest
and most effective manner in a given situation.  Most mathematicians, when they develop practical techniques, use
some kind of duality  (such as  distributions in functional analysis or flows in differential topology) which
simplifies the main ideas and formulations. It would probably not be a sufficient reason for using duality
if everything was exactly mirrored. Some properties are not preserved under categorical equivalence. This
gives rise to an additional dimension for those new constructions which are not reducible to either of the
opposites. Notions such as schemes arise in this way. Historically, the abstract concept of duality was
introduced much later than numerous (famous) concrete examples. A deep categorical analysis of first order
duality was given in \cite{P-Th}.

\head {\bf 0.2. Not everything is preserved under duality. Bifurcation theory} \endhead

As an example, this phenomenon is well-known in the qualitative theory of differential equations when small changes of
parameter cause catastrophes in the solution space (under the general duality of differential equations to
their solution spaces). This is the subject of bifurcation theory. The same phenomenon holds for algebraic
equations and for any type of equations and their deformations in the previous sense. The reason of 
this phenomenon is that duality between 'big' categories of equations and their solution spaces is 
always 1-duality and 'never' 2-duality. The exact criterion for the order of dulity can be 
translated to a criterion for bifurcations. 

\head {\bf 0.3. Development of Modern Geometry}\endhead

\vskip 0.2cm
{\bf Duality} plays a central role in the principal steps of the development of Modern Geometry. 
It is now a standard tool to talk about spaces which are unknown but which are 
well representable by their dual objects. All the development of modern  algebraic  geometry can be 
regarded as a sequence of extensions of algebraic duals, which can be seen from the following diagram:

\vskip 0.3cm
\hskip -1cm$
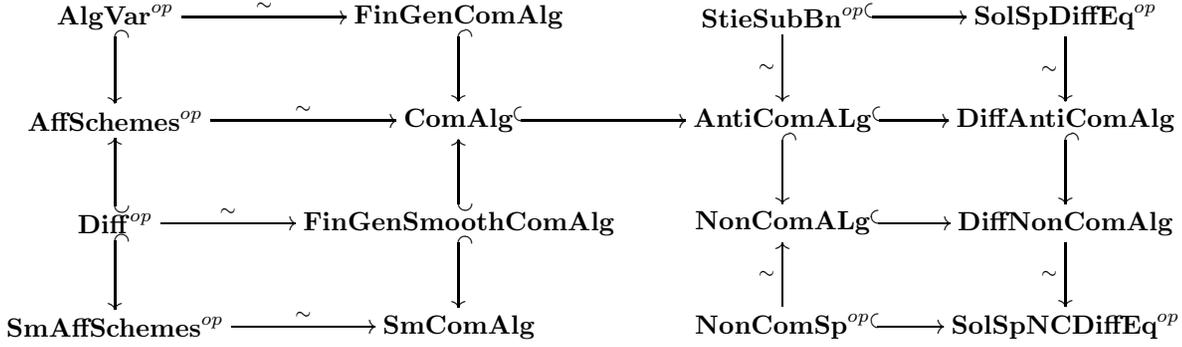
\xymatrix{\bold{AlgVar}^{op} \ar[r]^-{\sim }  \ar@{^{(}->}[d] &    \bold{FinGenComAlg} \ar@{^{(}->}[d] &   \bold{StieSubBn}^{op}  \ar[d]_-{\sim } \ar@{^{(}->}[r] & \bold{SolSpDiffEq}^{op}  \ar[d]_-{\sim } \\
\bold{AffSchemes}^{op}  \ar[r]^-{\sim }        &     \bold{ComAlg}    \ar@{^{(}->}[r]    &  \bold{AntiComALg} \ar@{^{(}->}[d] \ar@{^{(}->}[r] & \bold{DiffAntiComAlg}   \ar@{^{(}->}[d]  \\
\bold{Diff}^{op}   \ar[r]^-{\sim } \ar@{_{(}->}[u]  \ar@{^{(}->}[d]   &    \bold{FinGenSmoothComAlg}	\ar@{_{(}->}[u] \ar@{^{(}->}[d] & \bold{NonComALg}  \ar@{^{(}->}[r] &  \bold{DiffNonComAlg} \ar[d]_-{\sim } \\
\bold{SmAffSchemes}^{op}  \ar[r]^-{\sim } & \bold{SmComAlg} &  \bold{NonComSp}^{op} \ar[u]^-{\sim } \ar@{^{(}->}[r] &  \bold{SolSpNCDiffEq}^{op}   }$

\vskip 0.3cm
It is still a compact diagram, some steps are omitted, some extensions are not unique (e.g., for algebraic geometry, it is better 
to regard commutative algebras as anticommutative ones concentrated in degree 0. But for algebraic topology it is natural to regard 
graded commutative algebras as graded anticommutative ones with degrees of all elements doubled). 
One of the key ideas of this thesis is that the $\infty $-category setting allows us to expand this diagram in a new (homotopical) 
dimension. So that (monoidal) $\infty $-categories give an appropriate framework for Homotopical Algebraic Geometry.

\vskip 0.35cm
\head {\bf 0.4. Low-dimensional and $\infty $-dimensional Approaches to Homotopy Theory}\endhead

\vskip 0.2cm
Higher dimensional functors preserve ``homotopy" invariants but not in a canonical way, i.e. they usually do not preserve 
$\pi _*, \ H_*, \ H^*$, etc. This is because these ``homotopy" invariants are not formulated {\bf   internally}
in a category. 
For example, regard the classifying space functor $B:\bold{wTopGrp} \to \bold{Top}$. In $\bold{wTopGrp}$ (category of weak topological groups) there are two (noncomparable in general) 
2-categorical structures: when 2-cells are conjugations and when they are homotopy classes of homotopies; [the last structure is weaker
if we restrict ourselves to a subcategory of path connected weak topological groups].
We note:

\proclaim{\bf Proposition 0.4.1} The classifying space functor $B:\bold{wTopGrp}\to \bold{Top}$ is 
\item{$\bullet $} a 2-functor with respect to conjugations in $\bold{wTopGrp}$, 
\item{$\bullet $} a 2-functor and 2-equivalence (not 1-equivalence) with respect to homotopy classes of homotopies in $\bold{wTopGrp}$. Its quasiinverse is the
loop space functor $\Omega :\bold{Top}\to \bold{wTopGrp}$.      \hfill  $\square $
\endproclaim 

One  would expect that there are many relations between conjugation invariants and homotopy invariants for $\bold{wTopGrp}$ 
and $\bold{Top}$. Indeed, there are some such, but they are not straightforward.  For example, $H^*(BG)$ is a conjugation-invariant 
commutative anticommutative algebra. We would expect that it is a subalgebra (or, maybe a quotient algebra) of 
$\bold{AdInvPol}({\goth g})$ (the algebra of polynomials on the Lie algebra $\goth g$ invariant under conjugations) but this is not true in general, although  $\bold{AdInvPol}({\goth g})$ is isomorphic to $H^*(BG)$ for compact Lie groups
$G$ (Chern-Weil homomorphism). The relation between $H^*(G)$ and $H^*(BG)$ is rather complicated: it is given by the Eilenberg-Moore spectral sequence
$H^*(G)\otimes H^*(BG) \to 0$. Why does such a  nice equivalence $B:\bold{wTopGrp}@>\sim >>\bold{Top}:\Omega $ give such complicated relations 
between homotopy invariants? Because these homotopy invariants are not defined internally.

\vskip 0.2cm
The typical definitions of homotopy invariants, such as the functor $\pi _n(X)=[S^n,X]$,  are not invariant under 
2(and higher order)-functors, because the $n$-spheres $S^n$ are  not traditionally
determined categorically. One of the goals of this thesis is to introduce a new notion of 
homotopy group, which are invariant under higher-order categorical equivalences
Our reinterpretation of homotopy is as follows:

\vskip 0.2cm
{\bf Definition 0.4.1.} 
\item{$\bullet $} For an $\infty $-category $\bold{C}$, $I\in Ob\, \bold{C}$, and a point $x:I\to X$ ({\bf formal}) {\bf homotopy groups} of $X$ are defined 
as follows \, $\tilde \pi ^I_n(X,x):=\bold{Aut}(e^{n-1}x)/\sim $ \, (where $e^{n-1}$ is $n-1$ times application of the identity operation $e$). 
\item{$\bullet $} When functors $\bold{Aut}(e^{n-1}(-))/\sim $ are {\bf representable} the representing objects $\tilde S^n$ are called
({\bf formal}) {\bf spheres} (in this case we have $\tilde \pi ^I_n(X,x):=(\bold{Aut}(e^{n-1}(x))/\sim )@>\sim>> [\tilde S^n,X]$).      \hfill    $\square $ 

\vskip 0.2cm
For $\infty \text{-}\bold{Top}$ when $I=\bold{1}$ there is a homomorphism of the usual homotopy groups into our formal ones 
$\pi _n(X,x)\to \tilde \pi ^I_n(X,x)$ (induced by the quotient map $I^n/(I^{n-1}\times 0)\cup (I^{n-1}\times 1)\twoheadrightarrow S^n$). 
For a category $\infty \text{-}\bold{TopMan}_b$ of topological manifolds with boundary and homotopies relative to the  boundary, formal homotopy groups 
coincide with the usual ones $\tilde \pi ^I_n=\pi _n$ when $I=\bold{1}$. 
For $\infty \text{-}\bold{Top}_*$
(pointed spaces and maps) $\tilde \pi ^{\bold{1}}_n(X,x)=[\bold{1},X]$ are trivial for all $n$ although $[S^n,X]$ gives the usual homotopy groups.

\vskip 0.2cm
Both functors $B$ and $\Omega $ preserve the homotopy type of $\bold{1}$. So, $\tilde \pi ^{\bold{1}}_n(G)=\tilde \pi ^{\bold{1}}_n(BG)$. But these groups 
are trivial and they  give no information (if we change $\bold{1}$ to a more complicated object $I\in Ob\, \bold{Top}$ then the 
information can be very nontrivial).

\vskip 0.2cm
\proclaim{\bf Proposition 0.4.2} If $F:\Cal A \to \Cal B$ is an $\infty $-equivalence between {\bf full} topological subcategories 
of $\infty \text{-}\bold{TopMan}_b$ such that $\bold{1}\sim F(\bold{1})$ then $F$ preserves the usual homotopy groups.              \hfill     $\square $ 
\endproclaim 

\vskip 0.35cm
\head {\bf 0.5. Concrete Duality}\endhead

\vskip 0.2cm
The underlying philosophy of our theory of concrete duality is that the world is nonlinear and opposites converge rather  than diverge.

\vskip 0.2cm
{\bf Definition 0.5.1.} Two $n$-categories $\Cal A$ and $\Cal B$ are called {\bf concretely dual} if 
there exists a ``schizophrenic" object $D$ living in
both of theses categories such that hom-functors $\Cal A(-,D):\Cal A\to n\text{-}\bold{Cat}$ \, and \, 
$\Cal B(-,D):\Cal B\to n\text{-}\bold{Cat}$ factor through the other category, i.e. \hskip 0.8cm
$\vcenter{\xymatrix{\Cal A \ar[r]^-{F} \ar[dr]_-{\Cal A(-,D)} & \Cal B \ar[d]^-{V} \\
    &    n\text{-}\bold{Cat}}}$ \hskip 0.8cm and \hskip 0.8cm
$\vcenter{\xymatrix{\Cal B \ar[r]^-{G} \ar[dr]_-{\Cal B(-,D)} & \Cal A \ar[d]^-{U} \\
    &    n\text{-}\bold{Cat}}}$  \hskip 0.5cm where $F$, $G$ are equivalences quasiinverse to each other.    \hfill    $\square $

\vskip 0.2cm
The higher order the duality is,  the more (homotopy) invariants are preserved.

\vskip 0.2cm
If $\Cal A\hookrightarrow n\text{-}\bold{Top}$ then the forgetful functor $U:\Cal A\to n\text{-}\bold{Cat}$ is usually the composite of
inclusion and the $n$-groupoid functor 
$\xymatrix{\Cal A \ar@{^{(}->}[r] & n\text{-}\bold{Top} \ar[rr]^-{n\text{-}\bold{Top}(\bold{1}, - )} & & n\text{-}\bold{Cat} }$ 
[by Grothendieck's hypothesis $\infty \text{-}\bold{Top}(\bold{1}, - ):\infty \text{-}\bold{Top}\to \infty \text{-}\bold{Cat}$ is 
an equivalence with its image].

\vskip 0.2cm
The above factorization (lifting) of hom-functors is frequently {\bf initial}. For first 
order categories it was proven by Porst and Tholen \cite{P-Th} that initial means maximal and any other concrete duality factors through the 
initial (natural) one;  for higher order categories, the analogous statement has not been proven yet. We hope that the structures introduced here will be useful in extending this result to the higher order case. 

\vskip 0.35cm
\proclaim{\bf Proposition 0.5.1} 
\item{$\bullet $} Every (weak) duality (adjunction) $\xymatrix{\Cal A \ar@/^/[r]^-{F} & \Cal B \ar@/^/[l]^-{G} }$ is concrete 
(over $\Cal C$) if there are {\bf representable} forgetful functors $U:\Cal A\to \Cal C$ and $V:\Cal B\to \Cal C$.
The dualizing object $D$ is both $FI$ and $GJ$ in a sense to be made precise, where $I,J$ are representing objects for $U,V$.
\item{$\bullet $} If $\Cal A$ and $\Cal B$ have representable forgetful functors over $\Cal C$ and a dualizing object 
$D$ such that the corresponding hom-functors $\Cal A( - ,D)$, $\Cal B( - ,D)$ satisfy the {\bf initial lifting condition} 
(essentially, the arrow $f^n:VX\to \Cal A(Y,D)$ is a $\Cal B$-arrow iff the composite $VX@>f^n>>\Cal A(Y,D)@>ev_{x^n}>> \Cal A(I,D)$
is a $\Cal B$-arrow $\forall \, x^n:I\to Y$, and similarly for $\Cal B( - ,D)$)
then there exists a concrete dual adjunction between $\Cal A$ and $\Cal B$ which is natural and strict.
\item{$\bullet $} Concrete natural duality is a strict adjunction. [\, Higher order duality need not be an adjunction at all]  \hfill  $\square $
\endproclaim 

\vskip 0.2cm
Point 2 of the above proposition is a generalization (for $n$-categories) of the Porst-Tholen theorem about concrete duality for 
first order categories.

\vskip 0.2cm
The main and  most interesting interplay for duality is between algebra and geometry. Certain complicated colimits in 
algebraic categories are often easily viewed via duality as geometric limits (e.g. the notion of tensor product of algebras is more
understandable via the notion of product of manifolds).

\vskip 0.2cm
Examples of {\bf well-known dualities} are those between algebraic varieties and finitely-generated commutative algebras, 
between affine schemes and commutative rings (Grothendieck), compact abelian groups and abelian groups (Pontryagin), Boolean algebras
and Boolean spaces (Stone), commutative $C^*$-algebras and compact Hausdorff spaces (Gelfand-Naimark), and others.   

\vskip 0.2cm
In this paper several new examples of concrete duality are introduced. These include duality for differential equations (introducing 
anticommutative geometry of solution spaces), Vinogradov duality (formalizing the well-known duality between modules of linear differential 
operators and jet modules of sections), Gelfand-Naimark 2-duality (extending the usual one to homotopy classes of homotopies), 
Pontryagin-Lukacs duality (Lukacs' extension of Pontryagin duality to locally precompact abelian groups).

\vskip 0.5cm

\head {} {\bf 1. Categories, functors, natural transformations, modifications}\endhead

\vskip 0.2cm
{\bf Definition 1.1.} 
\item{$\bullet $} {\bf $\infty$-precategory} is a (big) set $L$ endowed with
\roster
\item a grading $L=\coprod\limits _{n\ge 0} L^n$
\item unary operations $d,c:\coprod\limits _{n\ge 1} L^n\to \coprod\limits _{n\ge 0} L^n$, $deg(d)=deg(c)=-1$, $dc=d^2$, $cd=c^2$
\item a unary operation $e:\coprod\limits _{n\ge 0} L^n\to \coprod\limits _{n\ge 0} L^n$, $deg(e)=1$, $de=1$, $ce=1$
\item partial binary operations $\circ _k$, $k=1,2,...$, of degree $0$. $f\circ _k g$ is determined iff $d^kf=c^kg$
\endroster
\item{}such that each {\bf hom-set} $L(a,a'):=\{f\in L\ |\ \exists \, k\in {\Bbb N} \ d^kf=a, \ c^kf=a'\}$, $deg(a)=deg(a')$, inherits all properties (1)-(4).
\item{$\bullet $} $\forall \, a,a',a'' \in L^m$ there are maps $\mu _{a,a',a''} :\coprod\limits _{n\ge 0}L^n(a',a'')\times L^n(a,a')\to L(a,a'')$ 
such that if the bottom composite is determined then
$$\xymatrix{
\coprod\limits _{n\ge 0}L^n(a',a'')\times L^n(a,a') \ar[rr]^-{\mu _{a,a',a''}} && L(a,a'')\\
L^n(a',a'')\times L^n(a,a') \ar[u]^{i\times i} \ar[rr]_-{\circ _{n+1}} && L^n(a,a'') \ar[u]_{i}\\
}$$
$\mu _{a,a',a''}$ are called {\bf horizontal composites} on level $deg(a)$; all composites inside of $L(a,a')$ are {\bf vertical}. \hfill $\square $
\vskip 0.2cm
{\bf Remarks.}
\item{$\bullet $} Our definition of {\bf $\infty$-precategory} coincides with what Penon calls 
a {\it magma}; essentially it is a reflexive globular set with all possible binary composites
\cite{Lei}.
 \item{$\bullet $} If $\alpha ^n, \beta ^n \in L^n, n>0$, such that $d\alpha ^n\ne d\beta ^n$ or $c\alpha ^n\ne c\beta ^n$, then 
$L(\alpha ^n,\beta ^n)=\emptyset$ (because of $d^2=dc, \, c^2=cd$). So, $\mu _{a,a',a''}$ can be the empty map $\emptyset :\emptyset \to \emptyset$. 
\item{$\bullet $} It is convenient to use a letter with appropriate superscript, like $x^m, \alpha ^k$, etc., as an element (or sometimes as a variable) 
with domain $L^m, L^k$, etc. respectively (or with domain $L^m(a,b), L^k(x,y)$, etc.) \, Also, the grading can be taken to range over  $\Bbb Z$ 
under the assumption that  $L^{-m}:=\emptyset , \ m>0$. 
\item{$\bullet $} Call elements $a\in L^0$ of degree $0$  {\bf objects} of $L$, elements $f^n\in L^n(a,a'), \, a,a'\in L^0$, {\bf arrows of degree $n+1$ from $a$ to $a'$}. 
\item{$\bullet $} Denote {\bf horizontal composites} by $*$\, , and extend it over arrows {\bf of different degrees} 
by the rule $*:L(b,c)\times L(a,b)\to L(a,c):(g^n,f^m)\mapsto \mu _{a,b,c}(e^{max(m,n)-n}g^n,e^{max(m,n)-m}f^m)=:g^n*f^m$
($f^m\in L^m(a,b), g^n\in L^n(b,c)$).  \hfill $\square $

The following definition of  equivalence is given  ``coinductively"  (see \cite{J-R})

{\bf Definition 1.2.} For $a,b\in L^n$ $a\sim b$ iff \ $\exists \
\xymatrix{a \ar@/^/[r]^f & b \ar@/^/[l]^g}$ such that $e(a)\sim g\circ _1f$ \  and  \ $f\circ _1g \sim e(b)$
(it means that there exists  an  $f\in L^0(a,b)$, $g\in L^0(b,a)$ and two  infinite sequences of arrows of higher order,
one in $L(a,a)$ and the other in $L(b,b)$; all this data we will call {\it arrows representing
the given equivalence}).        \hfill $\square $

\vskip 0.1cm
$\sim $ \, is reflexive and symmetric, but may be not transitive.

\proclaim{\bf {Lemma 1.1}} If $L$ is an $\infty $-precategory such that 
\item{}$\circ _1$ is weakly associative:   $f\circ _1(g\circ _1h)\sim (f\circ _1g)\circ _1h$ (for composable arrows),
\item{}$\circ _1$ satisfies the weak unit law:  
\ $\forall f\in \coprod\limits _{n\ge 1}L^n$ $\cases f\circ _1 edf\sim f & \\
ecf\circ _1 f\sim f & \endcases $\hskip -0.3cm,
\item{}$\sim $ \, is compatible with $\circ _1$ \, , i.e. $(f\sim g)\, \& \, (h\sim k)\Rightarrow (f\circ _1h)\sim (g\circ _1k)$ (for composable arrows), 
\item{}$\sim $ \, is transitive in higher orders: i.e.  there exists  $m>0$ such that if  $\sim $ is transitive for $\coprod\limits _{n\ge m}L^n$, 
\item{}then $\sim $ \, is transitive in all orders.
\endproclaim
\demo{Proof} Let $\vcenter{\xymatrix{a \ar@/^1ex/[r]^-{f} & b \ar@/^1ex/[l]^-{g} \ar@/^1ex/[r]^-{f'} & c \ar@/^1ex/[l]^-{g'}}}$ be the given 
equivalences, i.e. $ea\sim g\circ _1f$, \, $eb\sim f\circ _1g$, \, $eb\sim g'\circ _1f'$, \, $ec\sim f'\circ _1g'$. Then 
$\vcenter{\xymatrix{a \ar@/^1ex/[r]^-{f'\circ _1f} & c \ar@/^1ex/[l]^-{g\circ _1g'}}}$ is the required equivalence since 
$ea\sim g\circ _1f\sim g\circ _1(eb\circ _1f)\sim g\circ _1((g'\circ _1f')\circ _1 f)\sim (g\circ _1g')\circ _1(f'\circ _1f)$ and similarly
$ec\sim (f'\circ _1f)\circ _1(g\circ _1g')$. \vskip 0.0cm \hfill $\square $
\enddemo
{\bf Remarks.} 
\item{$\bullet $} Transitivity in higher orders trivially holds for $n$-categories (starting from level $n$), taking $\sim
$ as the identity. For proper 
$\infty $-categories it is better to make the assumption  ``\, $\sim $ \, is transitive in all orders" from the
beginning.
\item{$\bullet $} This lemma shows that although transitivity of \, $\sim $ \, is not automatic for $\infty $-precategories, it is indeed
consistent with (weak) associativity,  the unit law, and compatibility of \, $\sim $ \, with composites.   \hfill    $\square $
\vskip 0.2cm
{\bf Definition 1.3.} An $\infty $-precategory $L$ with relation \, $\sim $ \, as above is called a (weak) {\bf $\infty $-category} iff
\item{$\bullet $} \, $\sim $ \, is transitive:  \, $\alpha \sim \beta \sim \gamma \Rightarrow \alpha \sim \gamma $,
\item{$\bullet $} \, $\sim $ \, is compatible with all composites:  \, $(f\sim g)\, \& \, (h\sim k)\Rightarrow (f\circ _nh)\sim (g\circ _nk)$ 
(when they are defined),
\item{$\bullet $} horizontal composites preserve properties (1)-(2) and weakly preserve properties (3)-(4) of $\infty $-precategories in the following sense:
\roster
\item grading $deg_{L(a,a'')}(\mu _{a,a',a''}(f,g))=deg_{L(a',a'')}(f)=deg_{L(a,a')}(g)$
\item $\mu _{a,a',a''}(df,dg)=d\mu _{a,a',a''}(f,g)$, \ $\mu _{a,a',a''}(cf,cg)=c\mu _{a,a',a''}(f,g)$
\item $\mu _{a,a',a''}(ef,eg)\sim \, e\mu _{a,a',a''}(f,g)$
\item $\mu _{a,a',a''}(f\circ _k f',g\circ _k g')\sim \mu _{a,a',a''}(f,g)\circ _k \mu _{a,a',a''}(f',g')$ \ ("interchange law")
\endroster
\item{$\bullet $} each $\circ _k, k\in \Bbb N$, is {\bf weakly associative}:  $(f\circ _k g)\circ _k h\sim f\circ _k(g\circ _k h)$ 
(for composable elements),
\item{$\bullet $} The {\bf weak unit law} holds:   $e^{k}c^kf\circ _k f\sim f$, \ $f\circ _k e^kd^kf\sim f$ (when all operations are defined).  \vskip 0.0cm \hfill $\square $

\vskip 0.2cm
{\bf Remarks.} 
\item{$\bullet $} It is instructive to see what goes wrong if we attempt to consider a bicategory as
an instance of this definition. One would think that we could obtain an example by defining
$\sim$ on 1-cells as isomorphism of 1-cells and as equality for 2-cells. However, the problem lies in the horizontal composition of 2-cells  which would be required to be strictly associative, whereas in general
the horizontal composite of 2-cells is not. 
\item{$\bullet $} By lemma 1.1, \, for $n$-categories, the transitivity condition on \, $\sim $ \, follows from
the others.
\item{$\bullet $} Hom-sets in an $\infty $-category $L$ are $\infty $-categories themselves, and horizontal
composites \,  
$*:L(b,c)\times L(a,b)\to L(a,c)$, are $\infty $-functors.    
\item{$\bullet $} Since strict functors preserve the equivalences $\sim $ for categories in which horizontal composites preserve identity 
and composites strictly, the compatibility condition on $\sim $ with composites holds automatically.   \hfill  
$\square $

\vskip 0.2cm
A category is called {\bf strict} if the associativity and unit laws hold for elements (not just for $\sim $\,
-equivalence classes) and  horizontal composites preserve identities and composites strictly. Note that $\sim $ \, still
makes sense for strict categories.

\vskip 0.2cm
\proclaim{\bf {Proposition 1.1}} In a strict $\infty $-category $L$, arrows of degree $n$ (i.e., $L^n$) form 
a $1$-category with objects 
$L^0$, arrows $L^n$, domain function $d^n$, codomain function $c^n$. Observe that $d,c:L^n\to L^{n-1}$ are $1$-functors. \hfill $\square $
\endproclaim

\proclaim{\bf {Lemma 1.2}} 
\item{$\bullet $} In the strict $\infty $-category $L$ \ $e^k(f\circ _ng)=e^kf\circ _{n+k}e^kg$ (when either side is defined).
\item{$\bullet $} $\sim $ is preserved under $\sim $\, , i.e., if $\xymatrix{a \ar[r]^-{\sim }_-{f} & a'}$ is an equivalence with
$\xymatrix{a' \ar[r]^-{\sim }_-{g} & a}$, its quasiinverse (i.e. $ea\sim g\circ f$, $ea'\sim f\circ g$), and if $f'\sim f$ then $g$ is 
quasiinverse of $f'$ as well.
\item{$\bullet $} A quasiinverse is determined up to $\sim $\, , i.e. if 
$\xymatrix{a \ar@/^0.6pc/[r]^-{f}  &   b   \ar@/^0.5pc/[l]_-{g} \ar@/^1pc/[l]^-{g'} }$ and $g'\circ _1f\sim ea\sim g\circ _1f$ and
$f\circ _1g'\sim eb\sim f\circ _1 g$ then $g'\sim g$.
\vskip 0.1cm\item{$\bullet $} All $n+1$ composites in $\bold{End}(e^na):=L^{0}(e^na,e^na), \, n\ge 0$ coincide up to equivalence $\sim $\, .
\endproclaim
\demo{Proof}
\item{$\bullet $} Assume $f,g\in L^m, \ m\ge n$. Then $f\circ _ng=\mu _{d^ng,c^ng,c^nf}(f,g)$, which preserves $e$.
\item{$\bullet $} $ea\sim g\circ f\sim g\circ f'$, \ $ea'\sim f\circ g\sim f'\circ g$.
\item{$\bullet $} $g'=g'\circ _1eb\sim g'\circ _1f\circ _1g\sim g\circ _1f\circ _1g\sim g\circ _1eb=g$.
\item{$\bullet $} $f\circ _{n+1}g=\mu _{a,a,a}(f,g)\sim \mu _{a,a,a}(f\circ _k e^{n+1}a,e^{n+1}a\circ _k g)\sim \mu _{a,a,a}(f,e^{n+1}a)\circ _k\mu _{a,a,a}(e^{n+1}a,g)\sim f\circ _kg$, \, $1\le k\le n+1$.
\hfill $\square $
\enddemo

\vskip 0.1cm
{\bf Definition 1.4.} An arrow $(f:a\to a')\in L^0(a,a')$, \, $deg(a)=deg(a')=m\ge 0$, \, is called
\item{$\bullet $} {\bf monic} if $\forall g,h:z\to a$ if $f\circ _1g\sim f\circ _1h$ then $g\sim h$
\item{$\bullet $} {\bf epic} if $\forall g',h':a'\to w$ if $g'\circ _1f\sim h'\circ _1f$ then $g'\sim h'$
\item{$\bullet $} an {\bf equivalence} if  there exists $f':a'\to a$ such that $edf\sim f'\circ _1f$ and $edf'\sim f\circ _1f'$ \hfill $\square $

\proclaim{\bf Proposition 1.2} For composable arrows 
\item{$\bullet $} If $f,g$ are monics then $f\circ _1g$ is monic. If $f\circ _1g$ is monic then $g$ is monic
\item{$\bullet $} If $f,g$ are epics then $f\circ _1g$ is epic. If $f\circ _1g$ is epic then $f$ is epic
\item{$\bullet $} If $f,g$ are equivalences then $f\circ _1g$ is an equivalence \hfill $\square $
\endproclaim

\vskip 0.1cm
\proclaim{\bf Proposition 1.3} All arrows representing equivalence $a\sim b$ are equivalences.     \hfill    $\square $
\endproclaim 

\vskip 0.1cm
{\bf Definition 1.5.} An {\bf $\infty $-functor} $F:L\to L'$ is a function which 
strictly preserves the following properties (1)-(2) of precategories:
\roster
\item if $a\in L^n$ then $F(a)\in L^{'n}$
\item $F(da)=dF(a), F(ca)=cF(a)$
\endroster
and weakly preserves the following properties (3)-(4):
\roster
\item[3] $F(ea)\sim eF(a)$
\item $F(a\circ _k b)\sim F(a)\circ _k F(b)$ \hfill $\square $
\endroster

\vskip 0.1cm
{\bf Remark.} 
\item{$\bullet $} We do not require the functor $F$ to preserve equivalences $\sim $ because it is not automatic and
can be too  restrictive. However,   the functors preserving $\sim $ are very important (e.g., see point 1.2).  
\item{$\bullet $} The inverse map $F'$ for a bijective weak functor $F$ is not a functor, in general. If $F$
preserves
$\sim $ then  to say  the inverse map $F'$ is a (weak) functor is equivalent to saying $F'$ preserves $\sim $.
The inverse of a strict functor is  always a strict functor.    \hfill   $\square $  

\vskip 0.2cm
\proclaim{\bf Lemma 1.3}
\item{$\bullet $} Strict functors preserve equivalences $\sim $\, .
\item{$\bullet $} If functor $F\hskip -0.035cm:\hskip -0.035cmL\hskip -0.035cm\to \hskip -0.035cmL'$ is such that each  
restriction on hom-sets $F_{a,b}\hskip -0.035cm:\hskip -0.035cmL(a,b)\hskip -0.035cm\to \hskip -0.035cmL'(F(a),F(b))$, 
$a,b\in L^0$, preserves equivalences $\sim $\, , then $F$ preserves equivalences $\sim $\, . 
\item{$\bullet $} If $F\hskip -0.035cm:\hskip -0.035cmL\hskip -0.035cm\to \hskip -0.035cmL'$ is an embedding (injective map) 
such that $\forall a,b\hskip -0.0125cm\in \hskip -0.025cmL^0$ $F_{a,b}\hskip -0.035cm:\hskip -0.035cmL(a,b)\hskip -0.035cm\to \hskip -0.035cmL'(F(a),F(b))$ 
is a strict isomorphism and inverse $F'$ to codomain restriction of 
$\xymatrix{F:L \ar[r]_-{F\hskip -0.3cm\underset {Im(F)}\to |} & Im(F) \ar@{-->}@/_/[l]_-{F'} \ar@{^{(}->}[r] & L'}$
is a functor, then $F$ reflects $\sim $\, .
\endproclaim 
\demo{Proof}
\item{$\bullet $} Each arrow presenting a given equivalence $x\sim y$ is between a domain and a codomain which are constructed in a certain way
only by composites and identity operations from arrows of smaller degree presenting the given equivalence and from elements $x$ and $y$. 
A strict functor keeps the structure of the domains and codomains of arrows presenting the equivalence $x\sim y$. So, the image of 
arrows presenting an equivalence $x\sim y$ will be a family of arrows presenting an equivalence $F(x)\sim F(y)$.
\item{$\bullet $} For arrows of degree $>0$ equivalences are preserved by assumption.
Let $\xymatrix{a \ar@/^0.58pc/[r]^-{f}_-{\sim } & b \ar@/^0.58pc/[l]^-{g}_-{\sim } }$, $a,b\in L^0$, be an equivalence for objects in $L$, i.e. 
$ea\sim g\circ _1f$, $eb\sim f\circ _1g$. Then there are two opposite arrows $\xymatrix{F(a) \ar@/^0.58pc/[r]^-{F(f)} & F(b) \ar@/^0.58pc/[l]^-{F(g)} }$.
By assumption, $F(ea)\sim F(g\circ _1f)$, $F(eb)\sim F(f\circ _1g)$. So, $eF(a)\sim F(ea)\sim F(g\circ _1f)\sim F(g)\circ _1F(f)$ 
and $eF(b)\sim F(eb)\sim F(f\circ _1g)\sim F(f)\circ _1F(g)$. Therefore, $\xymatrix{F(a) \ar@/^0.58pc/[r]^-{F(f)}_-{\sim } & F(b) \ar@/^0.58pc/[l]^-{F(g)}_-{\sim } }$ 
is an equivalence.
\item{$\bullet $} The inverse of a strict isomorphism is a strict isomorphism, i.e. preserves equivalences. So, $F'$
is a functor  which preserves equivalences in all hom-sets and, consequently, preserves all equivalences.
Preservation of equivalences for $F'$  is exactly reflection of equivalences for $F$.         \hfill    $\square $
\enddemo 

\proclaim{\bf Lemma 1.4}
\item{$\bullet $} $x=y$ \, iff \, $ex\sim ey$ \, [in particular, \, $=$ \, is definable via \, $\sim $\, ].
\item{$\bullet $} Functors  preserving \, $\sim $\,   strictly preserve all composites \, $\circ _k$, $k\ge 1$.
\item{$\bullet $} Functors weakly preserving $e^2$ strictly preserve $e$, i.e. $e^2F(a)\sim F(e^2a)$ $\Rightarrow $ $eF(a)=F(ea)$.
\item{$\bullet $} {\bf Quasiequal} functors (i.e. $F(f^n)\sim G(f^n)$ for all $f^n\in L^n$, $n\ge 0$) are equal.
\endproclaim 
\demo{Proof}
\item{$\bullet $} $x=y$ $\Rightarrow $ $ex=ey$ $\Rightarrow $ $ex\sim ey$. Conversely, $ex\sim ey$ $\Rightarrow $ $dex=dey$ 
$\Rightarrow $ $x=y$. 
\item{$\bullet $} Sufficient to prove $eF(f\circ _kg)\sim e(F(f)\circ _kF(g))$, but it holds 
$eF(f\circ _kg)\sim F(e(f\circ _kg))\sim \text{($F$ preserves $\sim $) }F((ef)\circ _{k+1}(eg))\sim F(ef)\circ _{k+1}F(eg)\sim eF(f)\circ _{k+1}eF(g)\sim e(F(f)\circ _kF(g))$.
\item{$\bullet $} $e^2F(a)\sim F(e^2a)$ $\Rightarrow $ $de^2F(a)=dF(e^2a)$ $\Rightarrow $ $eF(a)=F(ea)$.
\item{$\bullet $} Again, it is sufficient to prove $eF(f^n)\sim eG(f^n)$. 
\item{}$eF(f^n)\sim F(ef^n)\sim \text{(by assumption) }G(ef^n)\sim eG(f^n)$.   \hfill   $\square $
\enddemo 

\vskip 0.0cm
{\bf Corollary.} {\it $\infty $-categories in the sense of definition 1.1.3 are almost strict}, namely, with strict associativity, 
identity, and interchange laws.
\demo{Proof} Strict associativity and strict identity laws hold because,  by the axioms, the functors $L(x,y)\times L(y,z)\times L(z,t)\to L(x,t):(f^n,g^n,h^n)\mapsto (h^n*g^n)*f^n$ and
$L(x,y)\times L(y,z)\times L(z,t)\to L(x,t):(f^n,g^n,h^n)\mapsto h^n*(g^n*f^n)$, $deg(x)=deg(y)=deg(z)=deg(t)$, are quasiequal, and, respectively, functors 
$L(x,y)\to L(x,y):f\mapsto f$ and $L(x,y)\to L(x,y):f\mapsto ey*f$, $deg(x)=deg(y)$ (similarly for the right identity), are quasiequal.
The strict interchange law holds because the functor $L(x,y)\times L(y,z):(f,g)\mapsto g*f$ preserves $\sim $.    
\hfill 
$\square $
\enddemo 

\vskip 0.1cm
{\bf Definition 1.6.} For two given functors $F$, $G$ , \ an  {\bf $\infty $-natural transformation} $\alpha :F\to
G$
\ is a function
$\alpha :L^0\to L^{'1}:a\mapsto (\xymatrix{F(a)\ar[r]^-{\alpha (a)} & G(a)})$ such that
$$\mu _{F(a),F(b),G(b)}(e^k\alpha (b), F(f))\sim \mu _{F(a),G(a),G(b)}(G(f),e^k\alpha (a))$$ for all $f\in L^{k}(a,b)$,
$k=0,1,...$ \hfill $\square $

{\bf Definition 1.7.} For two given functors $F$, $G$ and two natural transformations
$\xymatrix{F \ar@<1ex>[r]^-{\alpha }
\ar@<-1ex>[r]_-{\beta } & G}$ an 
{\bf $\infty $-modification} $\lambda :\alpha \to \beta $ is a function $\lambda :L^0 \to L^{'2}:a\mapsto (\xymatrix{\alpha (a) \ar[r]^-{\lambda (a)} & \beta (a)})$ such that
$$\mu _{F(a),F(b),G(b)}(e^k\lambda (b), F(f))\sim \mu _{F(a),G(a),G(b)}(G(f),e^k\lambda (a))$$
for all $f\in L^{k+1}(a,b)$,
$k=0,1,...$ \hfill $\square $
\vskip 0.3cm
\hskip -0.39cmAnalogously, modifications of higher order are introduced. We call modifications  
$1$-modifications, natural transformations   $0$-modifications.

{\bf Definition 1.8.} Given two functors $F$, $G$, two $0$-modifications $\xymatrix{F \ar@<1ex>[r]^-{\alpha ^0_1} \ar@<-1ex>[r]_-{\alpha ^0_2} & G}$, \newline
two $1$-modifications
$\xymatrix{\alpha ^0_1 \ar@<1ex>[r]^-{\alpha ^1_1} \ar@<-1ex>[r]_-{\alpha ^1_2} & \alpha ^0_2}$,..., two $n-1$-modifications
$\xymatrix{\alpha ^{n-2}_1 \ar@<1ex>[r]^-{\alpha ^{n-1}_1} \ar@<-1ex>[r]_-{\alpha ^{n-1}_2} & \alpha ^{n-2}_2}$ \newline
{\bf $\infty $-$n$-modification} $\alpha ^n:\alpha ^{n-1}_1\to \alpha ^{n-1}_2$ is a function
$\alpha ^n:L^0\to L^{'n+1}:$\newline
$a\mapsto (\xymatrix{\alpha ^{n-1}_1(a) \ar[r]^-{\alpha ^n(a)} & \alpha ^{n-1}_2(a)})$ such that
$$\mu _{F(a),F(b),G(b)}(e^k\alpha ^n (b), F(f))\sim \mu _{F(a),G(a),G(b)}(G(f),e^k\alpha ^n (a))$$
for all $f\in L^{k+n}(a,b)$,
$k=0,1,...$ \hfill $\square $

\vskip 0.2cm
{\bf Corollary.} {\it All $n$-modifications in the sense of Definition 1.1.8 are strict, i.e. 
all naturality squares commute strictly.} 
\demo{Proof} By the conditions in Definition 1.1.8, two functors $\alpha ^n(b)*F(-):L^{\ge n}(a,b)\to L^{'\ge n}(F(a),G(b))$ and 
$G(-)*\alpha ^n(a):L^{\ge n}(a,b)\to L^{'\ge n}(F(a),G(b))$ are quasiequal and, so, equal.  \vskip 0.0cm \hfill   $\square $
\enddemo

\vskip 0.2cm
{\bf Definition 1.9.} {\bf $\infty $-CAT} is an $\infty $-category consisting of
\item{$\bullet $} A  graded set $C=\coprod\limits _{n\ge 0}C^n$, where $C^0$ are categories, $C^1$ functors, $C^n$ $(n-2)$-modifications
\item{$\bullet $} if $\alpha ^n:\alpha ^{n-1}_1\to \alpha ^{n-1}_2\in C^n$ then $d\,\alpha ^n=\alpha ^{n-1}_1$, $c\,\alpha ^n=\alpha ^{n-1}_2$
\item{$\bullet $} $e\,\alpha ^n\in C^{n+1}$ is the map $L^0\to L^{'(n+1)}:a\mapsto e(\alpha ^n(a))$
\item{$\bullet $} for given two n-modifications $\alpha ^n_1$, $\alpha ^n_2$ such that $d^k\alpha ^n_1=c^k\alpha ^n_2$
$$\hskip -0.05cm\alpha ^n_1\circ _k\alpha ^n_2:=\cases
a\mapsto (\alpha ^n_1(a)\circ _k \alpha ^n_2(a))  & \text{if} \ \ k<n+2 \\
a\mapsto (\alpha ^n_1(F'(a))\circ _{(n+1)}G(\alpha ^n_2(a)))  & \text{if}
\ \ k=n+2, \ F'=c^{(n+1)}\alpha ^n_2, \ G=d^{(n+1)}\alpha ^n_1
\endcases $$
\hskip -0.35cmThe first composite works when $\alpha ^n_1, \alpha ^n_2\in \infty \text{\bf -CAT}(L,L')$, the second when 
$\alpha^n_1\in \infty \text{\bf -CAT}(L',L'')$
\vskip 0.0cm \hskip -0.35cm$\alpha ^n_2\in \infty \text{\bf -CAT}(L,L')$, where $L,L',L''$ are categories. \hfill $\square $

\vskip 0.2cm
\proclaim{\bf Lemma 1.5} In $\infty \text{-}\bold{CAT}$ there are two ways of taking horizontal composites (and
they are equal): 
$\alpha ^n*\beta ^n:=\alpha ^nF'\circ _{n+1}G\beta ^n=G'\beta ^n\circ _{n+1}\alpha ^nF$ (where $F:=d^{n+1}\beta ^n$, 
$F':=c^{n+1}\beta ^n$, $G:=d^{n+1}\alpha ^n$, $G':=c^{n+1}\alpha ^n$).
\endproclaim 
\demo\nofrills{Proof\ \ } follows from the naturality square for $\alpha ^n$ \hskip 0.5cm
$\vcenter{\xymatrix{G'F(a) \ar[r]^-{G'(\beta ^n(a))} &  G'F'(a) \\
   GF(a)  \ar[u]^-{\alpha ^n(F(a))} \ar[r]_-{G(\beta ^n(a))} &   GF'(a) \ar[u]_-{\alpha ^n(F'(a))} }}$    \hfill   $\square $
\enddemo 

\vskip 0.35cm
\proclaim{\bf {Proposition 1.4}} Categories, functors, natural transformations, modifications, etc. 
form the $\infty $-category {\bf $\infty $-CAT} of $\infty $-categories. \hfill $\square $
\endproclaim
\demo\nofrills{Proof\ \ } is similar to that for $2\text{-}\bold{CAT}$.
\item{$\bullet $} Horizontal composites preserve grading (obvious). 
\item{$\bullet $} $d$, $c$, $e$ are preserved for a similar reason, e.g., take $d$: 
$(d\alpha ^n)*(d\beta ^n)(a):=(d\alpha ^n)(F'(a))\circ _{n}G(d\beta ^n(a))=d(\alpha ^n(F'(a)))\circ _nd(G(\beta ^n(a)))=d(\alpha ^n(F'(a))\circ _{n+1}G(\beta ^n(a)))=d((\alpha ^n*\beta ^n)(a))=(d(\alpha ^n*\beta ^n))(a)$ 
(where $F':=c^{n+1}\beta ^n$, \, $G:=d^{n+1}\alpha ^n$).
\item{$\bullet $} (interchange law) 
$\vcenter{\xymatrix{L  \ar@/^2.3pc/[r]^-{F} \ar@/^1pc/[r]_-{F'}^-{\Downarrow \gamma ^n} \ar@/_1pc/[r]^-{F''} \ar@/_2pc/[r]_-{F'''}^-{\Downarrow \delta ^n} &  L'  \ar@/^2.3pc/[r]^-{G} \ar@/^1pc/[r]_-{G'}^-{\Downarrow \alpha ^n} \ar@/_1pc/[r]^-{G''} \ar@/_2pc/[r]_-{G'''}^-{\Downarrow \beta ^n} &   L''}}$
(by condition $d^k\delta ^n=c^k\gamma ^n$, $d^k\beta ^n=c^k\alpha ^n$, $k<n+2$, all $F$'s and $G$'s are functors)
\item{} $(\beta ^n\circ _k\alpha ^n)*(\delta ^n\circ _k\gamma ^n)(a):=(\beta ^n\circ _k\alpha ^n)(F'''(a))\circ _{n+1}G((\delta ^n\circ _k\gamma ^n)(a))=
(\beta ^n(F'''(a))\circ _k\alpha ^n(F'''(a)))\circ _{n+1}(G(\delta ^n(a))\circ _kG(\gamma ^n(a)))=$\newline 
\vskip 0.0cm \hskip -0.6cm$\cases \hskip -0.2cm(\beta ^n(F'''(a))\circ _{n+1}G(\delta ^n(a)))\circ _k(\alpha ^n(F'''(a))\circ _{n+1}G(\gamma ^n(a)))=(*)   &  \\
	\hskip -0.2cm\beta ^n(F'''(a))\circ _{n+1}\bigl (\alpha ^n(F'''(a))\circ _{n+1}G(\delta ^n(a))\bigr )\circ _{n+1}G(\gamma ^n(a))=(**) &     \endcases $\newline 
\vskip 0.0cm \hskip -0.6cm$\cases  \hskip -0.2cm(*)=((\beta ^n*\delta ^n)\circ _k(\alpha ^n*\gamma ^n))(a)  \  \text{ if } k<n+1 \ (\text{in this case } G=G'', G'=G''', F=F'', F'=F''')   &     \\
\hskip -0.2cm(**)=\beta ^n(F'''(a))\circ _{n+1}\bigl (G'(\delta ^n(a))\circ _{n+1}\alpha ^n(F''(a))\bigr )\circ _{n+1}G(\gamma ^n(a))=(**)     &        \endcases $\newline 
\vskip -0.8cm \hskip -0.6cm$\cases \hskip -0.2cm      &        \\
\hskip -0.2cm(**)=((\beta ^n*\delta ^n)\circ _{n+1}(\alpha ^n*\gamma ^n))(a)  \   \text{ if } k=n+1 \ (\text{in this case } F'=F'', G'=G'')                   &       \endcases $\newline 
\vskip -0.05cm\item{$\bullet $} The associativity law for vertical composites and the identity law hold essentially because of the componentwise 
definition of vertical composites. The associativity law for horizontal composites is due to the  interchange law and lemma 1.5. \hfill $\square $
\enddemo

\vskip 0.1cm
{\bf Definition 1.10.} A category $L$ is called an  {\bf $\infty $-$n$-category} if $L^{j+1}=e(L^j)$ for $j\ge n$.
\hfill $\square $

A quotient $L/\hskip -0.1cm\sim $ is not a category in general since $\sim $ is not compatible with $e$. However, if we take the quotient only on a fixed 
level $n$ and make all higher arrows identities we get $\infty $-$n$-category $L^{(n)}$, $n$-th approximation of $L$. 
Generally there are no functors $\xymatrix{L^{(n)} \ar@{^{(}->}[r] & L}$, $\xymatrix{L \ar@{->>}[r] & L^{(n)}}$ (except for the last 
surjection if $L$ is a weak $\infty $-$(n+1)$-category and all $(n+1)$-arrows are isomorphisms's).

\subhead {} 1.a. Weak categories, functors, natural transformations, modifications \endsubhead 

As we saw above, using a weak language (substituting $\sim $ for $=$) does not give a weak category theory.
The only  advantage was that we could deal with $\sim $ instead of $=$, which is important for the classification
problem (that still makes sense  for strict $\infty $-categories). All known definitions of weak categories
\cite{C-L, Lei, Koc} are nonelementary (at least, they use functors, natural transformations, operads, monads just 
for the very definition). Probably, this is a fundamental feature of weak categories. To introduce them we also
need the whole universe 
$\infty \text{-}\bold{PreCat}$ of $\infty $-precategories.

\vskip 0.2cm
{\bf Definition 1.a.1.} $\infty \text{-}\bold{PreCat}$ consists of 
\item{$\bullet $} {\bf $\infty $-precategories} (definition 1.1) together with $\sim $-relation in each [$\sim $ may be not transitive],
\item{$\bullet $} {\bf $\infty $-functors} (definition is like 1.5 for $\infty $-categories), i.e. functions $F:L\to L'$ of degree $0$
preserving $d$ and $c$ strictly, and $e$ and $\circ _k$, $k\ge 1$, weakly, 
\item{$\bullet $} {\bf lax $\infty \text{-}n$-modifications}, $n\ge 0$, i.e. {\bf total} maps $\alpha ^n:L\to L'$ (with variable degree on 
different elements, but $\le n+1$, more precisely, the induced map $\Bbb N\to \Bbb N:deg(x)\mapsto (deg(\alpha ^n(x))-deg(x))$ 
is an antimonotone map, decreasing by 1 at each step from $n+1$ at $deg(x)=0$ to 1 at $deg(x)=n$ and remaining constant 1 after) 
being defined for a given sequence of two functors $F,G:L\to L'$, two $0$-modifications (natural transformations) $\alpha ^0_1,\alpha ^0_2:F\to G$, ...,
two $(n-1)$-modifications $\alpha ^{n-1}_1,\alpha ^{n-1}_2:\alpha ^{n-2}_1\to \alpha ^{n-2}_2$ \ \ \ as \hskip 0.5cm $\alpha ^n:=$
\vskip 0.25cm
\item{}\hskip -0.25cm$\cases  \hskip -0.05cm(\alpha ^n(x):\alpha ^{n-1}_1(x)\to \alpha ^{n-1}_2(x))\in L^{'n}(F(x),G(x))  &  \hskip 0.25cmx\in L^0 \hskip -0.3cm \\
\hskip -0.05cm\alpha ^n(x):=\alpha ^n(e^{n+1-k}x)\in L^{'n+1}(F(d^kx),G(c^kx))  &   \hskip 0.15cm\underset {0<k<n+1}\to {x\in L^k}     \hskip -0.3cm  \\
\hskip -0.05cm(\alpha ^n(x):\alpha ^n(c^{n+1}x)\circ _{n+1}F(x)\to G(x)\circ _{n+1}\alpha ^n(d^{n+1}x))\in & \hskip 0.0cmx\in L^{n+1} \hskip -0.3cm \\        
\hskip 4.4cm\in L^{'n+1}(F(d^{n+1}x),G(c^{n+1}x))	&     \hskip -0.3cm  \\
\hskip -0.05cm(\alpha ^n(x):\alpha ^n(cx)\circ _1(e\alpha ^n(c^{n+2}x)\circ _{n+2}F(x))\to     & \hskip 0.0cmx\in L^{n+2} \hskip -0.3cm \\
\hskip -0.05cm(G(x)\circ _{n+2}e\alpha ^n(d^{n+2}x))\circ _1\alpha ^n(dx))\in L^{'n+2}(F(d^{n+2}x),G(c^{n+2}x))   &          \hskip -0.3cm  \\
\hskip -0.05cm\alpha ^n(x):\alpha ^n(cx)\circ _1(e\alpha ^n(c^2x)\circ _2(e^2\alpha ^n(c^{n+3}x)\circ _{n+3}F(x)))\to     &  \hskip 0.0cmx\in L^{n+3}   \hskip -0.3cm  \\
\hskip -0.05cm((G(x)\circ _{n+3}e^2\alpha ^n(d^{n+3}x))\circ _2e\alpha ^n(d^2x))\circ _1\alpha ^n(dx)\in L^{'n+3}(F(d^{n+3}x),G(c^{n+3}x)) &             \hskip -0.3cm  \\
  \dots                      &                            \\
\hskip -0.05cm\alpha ^n(x):    &   \hskip 0.0cmx\in L^{n+m}   \hskip -0.3cm  \\
\hskip -0.05cm\alpha ^n(cx)\circ _1\dots \circ _{m-2}(e^{m-2}\alpha ^n(c^{m-1}x)\circ _{m-1}(e^{m-1}\alpha ^n(c^{n+m}x)\circ _{n+m}F(x)\underbrace{)\dots )}_{m-1}\to &    \hskip -0.3cm  \\
\hskip -0.05cm\underbrace{(\dots (}_{m-1}G(x)\circ _{n+m}e^{m-1}\alpha ^n(d^{n+m}x))\circ _{m-1}e^{m-2}\alpha ^n(d^{m-1}x))\circ _{m-2}\dots \circ _1\alpha ^n(dx)\in  &       \hskip -0.3cm  \\
\hskip 4.4cm\in L^{'n+m}(F(d^{n+m}x),G(c^{n+m}x))  &       \hskip -0.3cm \\
\dots               &                \hskip -0.3cm \endcases $   
\item{}$d\alpha ^n:=\alpha ^{n-1}_1, \ c\alpha ^n:=\alpha ^{n-1}_2$ \ [$(d\alpha ^n)(x)\ne d(\alpha ^n(x)), \ (c\alpha ^n)(x)\ne c(\alpha ^n(x))$ if $deg(x)>0$].   \hfill    $\square $

\vskip0.2cm
{\bf Remarks.} 
\item{$\bullet $} $\infty \text{-}n$-modifications look terrible but they are the weakest form of naturality
(infinite sequences  of naturality squares arising by considering naturality squares given by
equations $e_1(x)\sim e_2(x)$ which express $\sim$-naturality in $x$.  This leads to an infinite
sequence of naturality squares). To deal with
such entities a kind of operad  is needed.
\item{$\bullet $} To give an $n$-modification $\alpha ^n$ is the same as to give a map $\alpha ^n\underset {L^0}\to
|:L^0\to L'$ of degree
$n+1$ and $\forall a,b\in L^0$ a natural transformation $\nu ^{\alpha ^n}_{a,b}:\alpha ^n(b)*F(-)\to G(-)*\alpha ^n(a):L^{\ge n}(a,b)\to L^{'\ge n}(F(a),G(b))$, 
where $F=d^{n+1}\alpha ^n$, $G=c^{n+1}\alpha ^n$.
\item{$\bullet $} When $\alpha ^n(x)$, $deg(x)>0$, are all identities (of the required types) $\infty \text{-}n$-modifications are 
called {\bf strict}. They are the usual modifications and composable as in definition 1.1.9 when 
the universe $\infty \text{-}\bold{CAT}$ 
is strict (in that case strict modifications are weak as well). In a weak universe $\infty \text{-}\bold{CAT}$ strict modifications
need not to be weak (i.e. to be modifications at all).
\item{$\bullet $} $\infty \text{-}\bold{PreCAT}$ is not an $\infty $-precategory itself because there are no identities and composites 
for weak $n$-modifications. 
The problem here with identities and composites is not clear, for example if they exist at all without making
either naturality condition or 
$\infty $-categories stricter. 
\item{$\bullet $} In general, these two sides ``categories and functors" and  ``$n$-modifications"
form a strange pair. If we weaken
one of these sides, the other one becomes stricter (under condition that $\infty \text{-}\bold{CAT}$ is a (let it be very weak) {\bf category}). 
So, the following {\bf hypothesis} holds:\vskip 0.1cm
\item{}{\it There is no $\infty \text{-}\bold{CAT}$ with {\bf simultaneously weak} categories, functors, and $n$-modification}. 
\vskip 0.1cm\item{} For example, if we want weak modifications and want them to be composable we need to introduce several axioms on categories, 
one of which is like '\, {\it $\forall a,b\in L^0$ and  $\forall $ functors $F,G:L\to L'$ if $\exists $ natural transformations 
$\alpha :f_1*F(-)\to G(-)*g_1:L^{\ge n}(a,b)\to L^{'\ge n}(F(a),G(b))$ and 
$\beta :f_2*F(-)\to G(-)*g_2:L^{\ge n}(a,b)\to L^{'\ge n}(F(a),G(b))$ and $n+1$-cells 
$f_1,f_2$ and $g_1,g_2$ are $\circ _k$-composable then $\exists $ a natural transformation ($k$-composite)
$\gamma :(f_1\circ _kf_2)*F(-)\to G(-)*(g_1\circ _kg_2):L^{\ge n}(a,b)\to L^{'\ge n}(F(a),G(b))$}\, '. But such axioms make very special 
categories. From the other side, if we want categories to be weak we need to make stricter (maybe, strict) $n$-modifications
in order that they would be composable. The problem is in existence of composites (and units) for weak $n$-modifications. 
\item{$\bullet $} Instead of lax $n$-modifications we could use modifications with $\alpha ^n(x)$ being $\sim $ for $deg(x)>0$ in $L'$. In both cases 
in order to make horizontal composites (at least, $F*\alpha ^n:=F\circ _{\bold{SET}}\alpha ^n$) we need functors preserving composites
(or composites and $\sim $), i.e. 'weak modifications' $\Rightarrow $ 'strict functors'.
\item{$\bullet $} If the above hypothesis was true it would be nice, e.g. a universe where $\infty \text{-}\bold{Top}$ 
lives would contain only strict $n$-modifications.  \hfill   $\square $

\vskip 0.2cm
{\bf Definition 1.a.2.} A {\bf weak $\infty $-category $L$} is an $\infty $-precategory (see definition 1.1)
such that 
\item{$\bullet $} $\sim $ is transitive \, $x\sim y\sim z$ $\Rightarrow $ $x\sim z$,
\item{$\bullet $} horizontal composites \, $*$ \, strictly preserve properties (1)-(2) of precategories
\itemitem{(1)} $deg(x*y)=deg(x)=deg(y)$ \, if \, $deg(x)=deg(y)$ \, (interchange law for degree)
\itemitem{(2)} $d(x*y)=(dx)*(dy)$, $c(x*y)=(cx)*(cy)$ \, if \, $deg(x)=deg(y)$ \, (interchange law for domain and codomain)
\item{}and weakly preserve properties (3)-(4) of precategories
\itemitem{(3)} $e(x*y)\sim (ex)*(ey)$ \, if \, $deg(x)=deg(y)$ \, (interchange law for identity)
\itemitem{(4)} $(x\circ _ky)*(z\circ _kt)\sim (x*z)\circ _k(y*t)$ \, if \, $deg(x)=deg(y)=deg(z)=deg(t)$ \, 
(interchange law for composites) \, [$\circ _k$ has smaller 'deepness' $k$ than the given $*=\circ _n$, $n>k$],
\item{$\bullet $} ({\bf weak associativity}) 
\item{}$\forall x,y,z,t\in L^n$ for two functors 
$l_{x,y,z,t}:L(x,y)\times L(y,z)\times L(z,t)\to L(x,t):(f,g,h)\mapsto (h*g)*f$ and 
$r_{x,y,z,t}:L(x,y)\times L(y,z)\times L(z,t)\to L(x,t):(f,g,h)\mapsto h*(g*f)$ 
$\exists $ natural transformation $\alpha _{x,y,z,t}:l_{x,y,z,t}\to r_{x,y,z,t}$, 
\item{$\bullet $} ({\bf weak unit}) 
\item{}$\forall x,y\in L^n$ and functors $u^{l}_{x,y}:L(x,y)\to L(x,y):f\mapsto ey*f$ and 
$u^{r}_{x,y}:L(x,y)\to L(x,y):f\mapsto f*ex$ $\exists $ natural transformations
$\epsilon ^l_{x,y}:u^l_{x,y}\to Id$ and $\epsilon ^r_{x,y}:Id\to u^r_{x,y}$.      \hfill   $\square $

\vskip 0.2cm
{\bf Remarks.} 
\item{$\bullet $} We do not introduce a universe $\infty \text{-}\bold{CAT}$ with weak categories, functors and $n$-modifications because
there are no (at least, obvious) units and composites for $n$-modifications (however, identity natural transformations exist if only the 
vertical composites of natural transformations are defined, for if $F:L\to L'$ is a functor take $(eF)(a):=e(F(a)), \, a\in L^0$ and by the weak unit law $\forall a,b\in L^0$ 
$\exists $ a natural transformation $\nu _{a,b}:e(F(b))*F(-)\to F(-)*e(F(a)):L^{\ge 0}(a,b)\to L^{'\ge 0}(F(a),F(b))$, take 
$\nu _{a,b}:=(\epsilon ^{u^r}_{F(a),F(b)}\circ _1\epsilon ^{u^l}_{F(a),F(b)})*F:=(\epsilon ^{u^r}_{F(a),F(b)}\circ _1\epsilon ^{u^l}_{F(a),F(b)})\circ _{\bold{SET}}F$). 
{\bf The problem} is what are the weakest conditions on categories,
functors and $n$-modifications in order that they form a category. Maybe  there are several independent such
conditions and, so, several categories $\infty \text{-}\bold{CAT}$ with weakest entities.
\item{$\bullet $} To keep a usual form of (weak) associativity and (weak) unit we could introduce relations $\underset k\to {\sim }$
for elements of images of two functors $F,G:L\to L'$ connected by a natural transformation $\alpha :F\to G$, namely, 
$x\underset k\to {\sim }y$ if $\exists $ $z\in L^k$ such that $x=F(z), \, y=G(z)$. These relations are not
reflexive, symmetric or  transitive. Then we could write associativity and unit laws as $(x\circ _ky)\circ
_kz\underset {k-1}\to {\sim }x\circ _k(y\circ _kz)$ and $e^kc^kx\circ _kx\underset {k-1}\to {\sim }x$, \,
$x\underset {k-1}\to {\sim }x\circ _ke^kd^kx$. Under assumption that composites and  units exist in an $\infty
\text{-}\bold{CAT}$ we could choose more sensible piece of $\infty \text{-}\bold{CAT}$ with categories in  which
$\underset 0\to {\sim }\, \equiv \, \sim $ and all $\underset k\to {\sim }$ are symmetric and transitive by the
requirement that
$\alpha _{x,y,z,t}, \, \epsilon ^l_{x,y}, \, \epsilon ^r_{x,y}$ are equivalences.          \hfill    $\square $

\vskip 0.35cm
\centerline{\bf Examples}
\vskip 0.2cm
\item{1.} {\bf $2$-Top} is a strict $\infty $-$2$-category with 2-cells, as homotopy classes of homotopies, and just identities in 
higher order (\, $\sim $ \, on the level of objects means homotopy equivalence of spaces, on the level of $1$-arrows homotopies of maps, and on the
level $\ge 2$ coincidence). $2\text{-}\bold{Cat}$ is similar.
\item{2.} It is widely believed that {\bf $\infty $-Top} is a (weak) $\infty $-category with homotopies between homotopies as higher order cells. It is hoped that this notion of (weak)  $\infty $-category (as above) will be useful in clarifying this issue. Assuming this, we can give two further
examples, as follows.
\item{3.}   {\bf $\infty $-Diff} is an $\infty $-category of differentiable manifolds in the same way as {\bf $\infty $-Top}.
\item{4.} {\bf $\infty $-TopALg} is an $\infty $-category of topological algebras in the same way as {\bf $\infty $-Top} where each instance of
homotopy is a homomorphism of topological algebras.

\item{5.} {\bf $\infty $-Compl} is an $\infty $-category of (co)chain complexes with (algebraic) homotopies for homotopies as higher order cells (see \cite{Lei} ).
\item{6.} For a $1$-category $A$, \ $A_{equiv}$ is a strict $\infty $-$2$-category such that $A^0_{equiv}=A^0$, $A^1_{equiv}=\Bigl\{f\in A \ \Bigl|\Bigr.\vcenter{\xymatrix{{\bullet } \ar[r]^-{\exists H}_-{\sim } & {\bullet }\\
{\bullet } \ar[u]^{f} \ar[r]_-{\forall h}^-{\sim } & {\bullet } \ar[u]_{f}}}
\Bigr\}$, $A^2_{equiv}=\Bigl\{\text{isomorphisms's} \ \Bigl|\Bigr. \ \forall f,g\in A^1_{equiv} \ \exists! \ \xymatrix{f \ar[r]^-{\gamma }_-{\sim } & g} \ \text{iff} \ \vcenter{\xymatrix{{\bullet } \ar[r]^-{\exists H}_-{\sim } & {\bullet }\\
{\bullet } \ar[u]^{f} \ar[r]_-{\forall h}^-{\sim } & {\bullet } \ar[u]_{g}}} \Bigr\}$. \newline
$A_{equiv}$ contains all equivariant maps $f:X\to Y$ with respect to a group homomorphism $\rho :\bold{Aut}(X)\to \bold{Aut}(Y)$.
\item{7.} The (weak) {\bf covariant} {\bf $\infty $-Hom}-functor $L(a,-):L\to \infty \text{\bf -CAT}:$
\vskip -0.3cm
$$\cases
b\mapsto L(a,b) & \ b\in L^0 \\
(f:b\to b')\mapsto (L(a,f):g\mapsto \mu (e^kf,g))& \ f\in L^0(b,b'), \ g\in L^{k}(a,b) \\
(\alpha :f\to f') \mapsto (L(a,\alpha ):x\mapsto \mu (\alpha ,ex ) & \alpha \in L^1(b,b'), \ x\in L^0(a,b) \\
(\delta :\alpha \to \alpha ') \mapsto (L(a,\delta ):x\mapsto \mu (\delta ,e^2x)) & \delta \in L^2(b,b'), \ x\in L^0(a,b) \\
\dots & \\
(\alpha ^n:\alpha ^{(n-1)}_1\to \alpha ^{(n-1)}_2) \mapsto (L(a,\alpha ^n):x\mapsto \mu (\alpha ^n,e^nx)) & \alpha ^n \in L^n(b,b'), \ x\in L^0(a,b)\\
\dots &
\endcases$$
\vskip -0.05cm
\item{8.} The {\bf opposite category} $L^{op}$ is an $\infty $-category such that
\itemitem{$\bullet $} $(L^{op})^n=L^n$, $n\ge 0$
\itemitem{$\bullet $} \hbox{$d^{op}(\alpha ^n)=\cases d(\alpha ^n) & \text{if} \ n\ge 2\\
c(\alpha ^n) & \text{if} \ n=1 \endcases $} \hskip 1cm
\hbox{$c^{op}(\alpha ^n)=\cases c(\alpha ^n) & \text{if} \ n\ge 2\\
d(\alpha ^n) & \text{if} \ n=1 \endcases $}
\itemitem{$\bullet $} $e^{op}=e$
\itemitem{$\bullet $} \hbox{$\beta ^n\circ ^{op}_k\alpha ^n=\cases
\beta ^n\circ _k\alpha ^n & \text{if } \ \alpha ^n,\beta ^n\in L^n, \ k<n \\
\alpha ^n\circ _k\beta ^n & \text{if } \ \alpha ^n,\beta ^n\in L^n, \ k=n
\endcases $}
\hskip 0.3cm
\hbox{(for composable elements)}
\vskip 0.2cm
\item{9.} The (weak) {\bf contravariant} {\bf $\infty $-Hom}-functor $L(-,b):L^{op}\to \infty \text{\bf -CAT}:$
$$\cases
a\mapsto L(a,b) & \ a\in L^0 \\
(f:a\to a')\mapsto (L(f,b):g\mapsto \mu (g,e^kf))& \ f\in L^0(a,a'), \ g\in L^{k}(a',b) \\
(\alpha :f\to f') \mapsto (L(\alpha ,b):x\mapsto \mu (ex,\alpha ) & \alpha \in L^1(a,a'), \ x\in L^0(a',b) \\
(\delta :\alpha \to \alpha ') \mapsto (L(\delta ,b):x\mapsto \mu (e^2x,\delta )) & \delta \in L^2(a,a'), \ x\in L^0(a',b) \\
\dots & \\
(\alpha ^n:\alpha ^{(n-1)}_1\to \alpha ^{(n-1)}_2) \mapsto (L(\alpha ^n ,b):x\mapsto \mu (e^nx,\alpha ^n )) & \alpha ^n \in L^n(a,a'), \ x\in L^0(a',b)\\
\dots &
\endcases$$
\item{10.} The {\bf Yoneda embedding} \ $\bold Y:L\to \infty\text{\bf -CAT}^{L^{op}}:\alpha \mapsto L(-,\alpha )$, \ $\alpha \in L$, 
where $L$ is an $\infty $-category.
\item{11.} $\bold{Set}$ is simultaneously an object and a full subcategory of $\infty \text{-}\bold{CAT}$.
\item{12.} A (big) set $L_{\sim }:=\coprod\limits _{n\ge 0}L^n_{\sim }$, where $L^n_{\sim }$ are defined recursively as 
$L^0_{\sim }:=L^0$ and $L^n_{\sim }$ are all equivalences from $L^n$ with domain and codomain in $L^{n-1}_{\sim }$, is a subcategory
of $L$. Similarly, $L_{k\sim }:=\coprod\limits _{n\ge 0}L^n_{k\sim }, \, k\ge 0$, where 
$L^n_{k\sim }:=\cases \hskip -0.05cmL^n & n\le k \hskip -0.05cm\\
\hskip -0.05cm\text{equivalences from }L^n\text{ with dom and codom in }L^{n-1}_{k\sim } & n>k \hskip -0.05cm\endcases $, is a subcategory of $L$. From this point
$L_{\sim }=L_{0\sim }$. Such categories are most important for the classification problem (up to $\sim $).
Sometimes, 'invariants' can be constructed only for $L_{\sim }$ (see point 1.2.1).
\item{13.} {\bf Higher order concepts} can {\bf simplify} proof of first order facts. E.g., each strict $2$-functor 
$\Phi :2\text{-}\bold{CAT}\to 2\text{-}\bold{CAT}$, where $2\text{-}\bold{CAT}$ is the usual strict category of categories, functors, and natural transformations, preserves adjunction $\biggl ($indeed, triangle identities $\cases \varepsilon G\circ G\eta =1_G   & \hskip -0.35cm\\
F\varepsilon \circ \eta F=1_F  & \hskip -0.35cm\endcases $ are respected by $\Phi $ $\cases \Phi (\varepsilon)\Phi (G)\circ \Phi (G)\Phi (\eta )=1_{\Phi (G)}   & \hskip -0.35cm\\
\Phi (F)\Phi (\varepsilon )\circ \Phi (\eta )\Phi (F)=1_{\Phi (F)}  & \hskip -0.35cm\endcases $ $\biggr )$. It gives short proofs of the following results.
\item{a)} {\it Right adjoints preserve limits (left adjoints preserve colimits).}
\item{}{\smc Proof}. \hskip 1cm$\vcenter{\xymatrix{\bold{A}^{\bold{I}} \ar[rrr]_-{F^{\bold{I}}} \ar@/^1.25pc/[dd]^-{lim} \ar@/_1.6pc/[dd]_-{colim}^-{\dashv } & & & \bold{B}^{\bold{I}} \ar@/^1.25pc/[dd]^-{lim} \ar@/_1.6pc/[dd]_-{colim}^-{\dashv } \ar@/_1.25pc/[lll]^-{\perp }_-{G^{\bold{I}}} \\
    & & & \\
\bold{A} \ar[uu]^-{\Delta }_-{\, \dashv } \ar[rrr]_-{F} & & & \bold{B} \ar[uu]^-{\Delta }_-{\, \dashv } \ar@/_1.25pc/[lll]_-{G}^-{\perp } }}$
\item{}where \, $(-)^{\bold{I}}\equiv 2\text{-}\bold{CAT}(\bold{I},-):2\text{-}\bold{CAT}\to 2\text{-}\bold{CAT}$ is a hom-$2$-functor.
\item{}Now, $G^{\bold{I}}\circ \Delta =\Delta \circ G$ (obvious). Taking right adjoints of both sides completes the proof 
$lim\circ F^{\bold{I}}\simeq F\circ lim$ (for colimits the same argument works $F^{\bold{I}}\circ \Delta =\Delta \circ F\, \Rightarrow \, colim \circ G^{\bold{I}}\simeq G\circ colim$).   \hfill $\square $
\item{b)} {\it Each $1\text{-}\bold{Cat}$-valued presheaf admits a sheafification ($1\text{-}\bold{Cat}$ is a category of small categories and functors between them)}.
\item{}{\smc Proof}. $1\text{-}\bold{Cat}$-valued presheaf on $\bold{C}$ is the same as an internal category object in $\bold{Set}^{\bold{C}^{op}}$.
There is an adjoint situation $\xymatrix{\bold{Sh}(\bold{C}) \ar@{^{(}->}[r] & \bold{Set}^{\bold{C}^{op}} \ar@/_1.1pc/[l]^-{\perp } }$
in $\bold{LEX}$, where $\bold{LEX}\hookrightarrow 2\text{-}\bold{CAT}$ is a $2$-category of finitely complete categories, functors 
preserving finite limits, and (arbitrary) natural transformations. There is a $2$-functor $\bold{CAT}(-):\bold{LEX}\to 2\text{-}\bold{CAT}$
assigning to each category in $\bold{LEX}$ the category of its internal category objects and to each functor and natural transformation 
the induced ones. Then $\exists $ an adjunction $\xymatrix{\bold{CAT}(\bold{Sh}(\bold{C})) \ar[r] & \bold{CAT}(\bold{Set}^{\bold{C}^{op}}) \ar@/_1.1pc/[l]^-{\perp } }$
which means that each $1\text{-}\bold{Cat}$-valued presheaf can be sheafified by the top curved arrow.    \hfill  $\square $

\subhead {} 1.1. Fractal organization of the new universe\endsubhead

\vskip 0.2cm
{\bf Fractal Principle.} Object $A$ with properties $\{P_i\}_I$ has fractal structure if there are subobjects $\{A_j\}_J$ which relate
to each other in a certain way (express it by additional property $P=$'to have $|J|$ subobjects which relate in the certain way') and
each $A_j$ inherits all properties $\{P_i\}_I\&P$.

\hfill $\square $

In spite of its complicated structure, each $\infty $-category and even {\bf $\infty $-CAT} itself, has a regular structure
which is repeated for certain arbitrary small pieces. Such pieces are, of course, the hom-sets $L(a,b)$ which inherit all properties (1)-(4),
associativity and identity laws, and each piece of which still has the same structure. In particular, $L(a,b)(c,d)=L(c,d)$. An $\infty $-functor restricted to such a piece is again an $\infty $-functor.
Moreover, each $\infty $-category can be regarded as a hom-set of a little bit bigger category if we formally attach two distinct elements 
$\alpha , \beta \in L^{-1}$ with their identities of higher order $e^n(\alpha ),\ e^n(\beta ), \ n\ge 1$ 
(such that $d(L^0)=\alpha ,\ c(L^0)=\beta $ and composites with these identities of other elements hold strictly).
Other natural pieces of $L$ which inherit all properties and are $\infty $-categories are $L^{\ge n}$, $L^{\ge n}(a,b)$ (elements of degree not lower than $n$).

\subhead {} 1.2. Notes on Coherence Principle\endsubhead

This principle is an axiom to deal with the equivalence relation $\sim $. It is not logically necessary for higher order category theory itself. There
can be categories in which it does not hold.

{\bf Coherence Principle.} For a given set of cells $\{a_i\}_I$ and a given set of base equivalences
$\{t_j(\{a_i\}_I)\sim s_j(\{a_i\}_I)\}_J$ for any two constructions 
$F_1(\{a_i\}_I)$ and $F_2(\{a_i\}_I)$ and any two derived
equivalences $\varepsilon ^0_i:F_1(\{a_i\}_I)\sim F_2(\{a_i\}_I)$, $i=1,2$ there are derived equivalences $\varepsilon ^1_{m}:
\varepsilon ^0_1\sim \varepsilon ^0_2$, $m\in M^1$, such that for any two of them $\varepsilon ^1_{m_1}$, $\varepsilon ^1_{m_2}$ there are derived
equivalences $\varepsilon ^2_{m}:\varepsilon ^1_{m_1}\sim \varepsilon ^1_{m_2}$, $m\in M^2$ again such that for any pair of them $\varepsilon ^2_{m_1}$, $\varepsilon ^2_{m_2}$
there are derived equivalences of higher order, etc. \hfill $\square $

Here constructions mean application of composites, functors, natural transformations,.. to $\{a_i\}_I$. Derived equivalences mean
equivalences obtained from base ones by virtue of the categorical axioms.

\head {\bf 2. $(m,n)$-invariants} \endhead 

\vskip 0.2cm
{\bf Definition 2.1.} 
\item{$\bullet $} {\bf Equivalence} \, $x^k\sim y^k$, \, $x^k,y^k\in L^k$, \, $k\ge 0$, is called {\bf of degree $l$}, \, $deg(\sim ):=l$, $l\ge 0$, if all arrows 
representing it (starting from order $k+l+1$ and higher) are identities and for $l>0$ there is at least one nonidentity arrow on level 
$k+l$. If there is no such \, $l\in \Bbb N$, \, $deg(\sim ):=\infty $. Denote $\sim $ of degree $l$ by $\sim _l$.
\item{$\bullet $} A pair of equivalent elements \, $x^k\sim y^k$, \, $k\ge 0$, is called {\bf of degree $l$}, \, $deg(x^k\sim y^k):=l$, $l\ge 0$,  
if the lowest degree of equivalences existing between $x^k$ and $y^k$ is \, $l$.
\item{$\bullet $} An {\bf $\infty $-category} $L$ is called {\bf of degree l}, \, $deg(L)=l$, \, $l\ge 0$, if for any pair of 
equivalent objects $a\sim a'$, $a,a'\in L^0$, there exists an equivalence $a\sim _ka'$ of degree $k\le l$ and there exists at least
one pair of equivalent objects from $L$ of degree $l$.
\item{$\bullet $} A {\bf Functor} $F:L\to L'$ is called {\bf $(m,n)$-invariant} if $F$ preserves equivalences $\sim $\, ,
$m=deg(L)$, \ $0\le n\le deg(L')$ and $F$ maps every
pair of equivalent objects of degree $\le m$ to a pair of equivalent objects of degree $\le n$, i.e. $deg(a\sim a')\le m\Rightarrow deg(F(a)\sim F(a'))\le n$,
and boundary $n$ is actually achieved on a pair of equivalent objects of $L$.       \hfill        $\square $

\vskip 0.2cm
{\bf Remarks.} 
\item{$\bullet $} $(m,n)$-invariants are important for the classification problem (up to $\sim $). If $n<m$ an $(m,n)$-invariant decreases
complexity of the equivalence relation, i.e. partially resolves it.   
\item{$\bullet $} There can be trivial invariants which do not distinguish anything and do not carry any information 
such as constant functors $c:L\to L'$ (although they are $(deg(L),0)$-invariants).  \vskip -0.35cm   \hfill   $\square $

\vskip 0.3cm
\centerline{\bf Examples}
\vskip 0.25cm
\item{1.} $deg(ea)=0$; \, $deg(f:a@>\sim >isomorphisms>a')=1$; \, $deg(\bold{Set})=1$; \, $deg(\infty \text{-}\bold{Top})=2$; \, $deg(\infty \text{-}\bold{CAT})=\infty \, (?)$. 
\vskip 0.1cm
\item{2.} Homology and cohomology functors $H_*,H^*:\infty \text{-}\bold{Top}\to \bold{Ab}$ (trivially extended over higher order cells) are $(2,1)$-invariants.
\vskip 0.1cm
\item{3.} $\tilde \pi ^I_n/\hskip -0.1cm\sim \hskip 0.05cm:L^*_{1\sim }\to \bold{Grp}$ is an $(\infty ,1)$-invariant (see proposition 2.1.2).
\vskip 0.1cm
\item{4.} Let $X$ be a smooth manifold with Lie group action $\rho :G\times X\to X$, \, $L$ be a category with $L^0$, the set of 
submanifolds of $X$, $L^1(a,b):=\{(a,g,b)\in L^0\times G\times L^0\, |\, \rho (g,a)=b\, \}$, \, $L^n:=eL^{n-1}$ for  $n\ge 2$, \, 
$L'$ be a category with $L^{'0}:=C^{\infty }(X,\Bbb R)$ (smooth functions), 
$L^{'1}(f,h):=\{(f,g,h)\in L^{'0}\times G\times L^{'0}\, |\, f\circ \rho (g^{-1},-)=h\, \}$, \, $L^{'n}:=eL^{'n-1}$ for  $n\ge 2$.
If $F:L\to L'$ is a construction (functor) assigning invariant functions to objects from $L$ then $F$ is a $(1,0)$-invariant.  
\item{5.} Each equivalence $L@>\sim >>L'$ is \, $(deg(L),deg(L'))$-invariant with \, $deg(L)=deg(L')$. 

\vskip 0.3cm
\subhead 2.1. Homotopy groups associated to $\infty $-categories \endsubhead

\vskip 0.2cm
Let $L$ be an $\infty $-category in which $*$ strictly preserves $e$ and $\sim $ (i.e. $*$ is a strict functor).
Denote by $eqL:=\{f\in L\ |\ \exists \, g. \ edf\sim g\circ _1f, edg\sim f\circ _1g\}$ the subset of equivalences of the $\infty $-category $L$. 
It may not be   a category (because it is not closed under $d,c$, in general).
\vskip 0.2cm
{\bf Definition 2.1.1.} Assume, $L(I,-):L\to \infty \text{\bf -CAT}$, \, $x\in L^0(I,a)$. Then $\tilde \pi ^I_{n}(a,x):=$\newline
$$\cases (L^0(I,a),x) & \ \ \text{if} \ $n=0$ \\
\bold{Aut}_{L(I,a)}(e^{n-1}x):=eqL(I,a)(e^{n-1}x,e^{n-1}x)\cap (L(I,a))^0(e^{n-1}x,e^{n-1}x)= & \\
=eqL(e^{n-1}x,e^{n-1}x)\cap L^{n+1} & \ \ \text{if} \ n>0
\endcases$$
are (weak) {\bf homotopy groups} of object $a$ at point $x$ with representing object $I\in L^0$.  \hfill $\square $
\vskip 0.2cm
$\tilde \pi ^I_{0}(a,x)$ or $\tilde \pi ^I_{0}(a,x)/\sim $ are just pointed sets, \ $\tilde \pi ^I_{n}(a,x)/\sim , n>0$ are strict groups.
\vskip 0.2cm
{\bf Remarks.} 
\item{$\bullet $} If $L=\infty \text{-}\bold{Top}$, \ $I=\bold 1$ then $\tilde \pi ^I_n(X,x)=[I^n/(I^{n-1}\times 0)\cup (I^{n-1}\times 1),X]$.  
The quotient map $(I^{n-1}\times 0)\cup (I^{n-1}\times 1)\twoheadrightarrow S^n$ induces a group homomorphism $\pi _n(X,x)\to \tilde \pi ^I_n(X,x)$.
\item{}If $L=\infty \text{-}\bold{TopMan}_b$ (the infinity category of topological manifolds with boundary as objects, and homotopies 
relative to the boundary as higher order cells), $I=\bold{1}$ then $\tilde \pi ^I_n(X,x)=[I^n/(I^{n-1}\times 0)\cup (I^{n-1}\times 1),X]\, \text{rel}\, (\partial I^n)=[S^n,X]=\pi _n(X,x)$ 
(i.e. formal homotopy groups coincide with the usual ones).
\item{$\bullet $} In the case when functors $\tilde \pi ^I_n$ are {\bf representable} (by certain cogroup objects $\tilde S^n_I$) we call the 
representing objects $\tilde S^n_I$  ({\bf formal}) {\bf spheres}. It makes sense to define (as usual) $\tilde \pi ^I_n(a):=[\tilde S^n_I,a]$, but these two 
definitions will not  always be equivalent. The first one is more internal, and the only external parameter is $I$. 
\item{$\bullet $} Any $\infty $-functor $F:L\to L'$, preserving $\sim $\, , induces (weak) {\bf group homomorphisms} 
$F_{I,a}:\tilde \pi ^I_n(a)\to \tilde \pi ^{FI}_n(Fa)$. So, for example, every $\infty $-equivalence between full 
subcategories of $\infty \text{-}\bold{TopMan}_b$, preserving the homotopy type of $\bold{1}$, will preserve the (usual) homotopy groups. \vskip 0.0cm \hfill $\square $
\vskip 0.2cm
{\bf Definition 2.1.2.} For a map \ $f:a\to b$ \ such that $f\circ x=y$, \ $x\in L^0(I,a)$, \ $y\in L^0(I,b)$ the 
{\bf induced map} \ $f_*\equiv \tilde \pi ^I_n(f):\tilde \pi ^I_{n}(a,x)\to \tilde \pi ^I_{n}(b,y)$ \ is determined by 
restriction of the functor $L(I,f):$\newline
$$\cases L^0(I,a)\to L^0(I,b):x'\mapsto f\circ _1x' & \ \ \text{if} \ n=0 \\
\bold{Aut}_{L(I,a)}(e^{n-1}x)\to \bold{Aut}_{L(I,b)}(e^{n-1}y):g\mapsto \mu _{I,a,b}(e^nf,g) & \ \ \text{if} \ n>0
\endcases $$
\vskip -0.6cm
\hfill $\square $

\vskip 0.2cm
{\bf Remark.} To be correctly defined,  induced maps $\tilde \pi ^I_n(f)$ for $n>1$ need commutativity of $*$ with $e$. 
The first two ``groups"
$\tilde \pi ^I_0(a,x), \tilde \pi ^I_1(a,x)$ always make sense and depend functorially on objects.    \hfill  $\square $

\vskip 0.2cm
\proclaim{\bf {Proposition 2.1.1 (homotopy invariance of homotopy groups)}} If $x:I\to a$, $f\sim f'\in L^0(a,b)$ such that
$f\circ _1x\sim f'\circ _1x$ is a trivial equivalence (all arrows for $\sim $ are identities) then 
$\tilde \pi ^I_n(f)/\sim \, =\tilde \pi ^I_n(f')/\sim \, :\tilde \pi^I_n(a,x)/\sim \, \to \tilde \pi ^I_n(b,f\circ x)/\sim \, $.
\endproclaim
\demo\nofrills{Proof\ \ } is immediate. \hfill $\square $
\enddemo

\proclaim{\bf Proposition 2.1.2} $\tilde \pi ^I_n/\hskip -0.1cm\sim \hskip 0.05cm:L^*_{1\sim }\to \bold{Grp}$ is an $(\infty ,1)$-invariant, 
where $L^*_{1\sim }:=\coprod\limits _{n\ge 0}L^{*n}_{1\sim }$, \, $L^{*n}_{1\sim }:=\cases  L^{*n}\text{ (pointed objects and maps) } & n=0,1  \\
\text{equivalences from }L^{n} \text{ with dom and codom in }L^{*(n-1)}_{1\sim }    & n>1     \endcases $. 
\endproclaim 
\demo{Proof} The partial functor $\tilde \pi ^I_{n}/\hskip -0.1cm\sim \hskip 0.05cm:L^{*0}\coprod L^{*1}\to \bold{Grp}$ is trivially extendable starting from 
equivalences on level 2 (because of proposition 2.1.1).    \hfill   $\square $
\enddemo 

\vskip 0.1cm
\centerline{{\bf Example} (Fundamental Group)}
\vskip 0.2cm
Let $2\text{-}\bold{Top}$ be the usual $\bold{Top}$ with homotopy classes of homotopies as $2$-cells. 
Define the {\bf fundamental groupoid} $2$-functor as the representable 
\, $\Pi (-):=Hom_{2\text{-}\bold{Top}}(1,-):2\text{-}\bold{Top}\to 2\text{-}\bold{Cat}:$\vskip 0.2cm \hskip -0.39cm
$\cases  X\to \Pi (X) & Ob\, (\Pi (X)) \text{ are its points, \, } Ar\, (\Pi (X)) \text{ are homotopy classes of paths} \\
(X@>f>>Y)\mapsto \Pi (f) & \text{transformation of fundamental groupoids, \ } \Pi (f):\cases x\mapsto f(x) & \\ 
[\gamma ]\mapsto [f\circ \gamma ] & \endcases \\
(f@>[H]>>f')\mapsto \Pi ([H]) &  \text{nat. trans. } \Pi ([H])=\underset {{x\in X}}\to {\{[H]*i_x\}}:Hom_{2\text{-}\bold{Top}}(1,f)@>\sim>>Hom_{2\text{-}\bold{Top}}(1,f')  \hskip -0.2cm  \endcases $ 
\vskip 0.0cm
(where \, $\underset {{x\in X}}\to {\{[H]*i_x\}}=\underset {x\in X}\to {\{[H(x,-)]\}}$ \, are homotopy classes of paths between $f(x)$ and $f'(x)$ natural in $x\in X$).
\vskip 0.15cm\hskip -0.39cm 
$\pi _1(X,x_0):=\bold{Aut}_{\Pi (X)}(x_0)\hookrightarrow \Pi (X)$ is the {\bf fundamental group} of the space $X$ at point $x_0\in X$, 
$\pi _1((X,x_0)@>f>>(Y,y_0)):=\bold{Aut}_{\Pi (X)}(x_0)@>\Pi (f)>>\bold{Aut}_{\Pi (Y)}(y_0)$.

\proclaim{\bf Proposition 2.1.3} \vskip 0.05cm
\item{$\bullet $} If \, $[H]:f@>\sim >>f':X\to Y$ \, is a $2$-cell in \, $2\text{-}\bold{Top}$ \, then \, $\pi _1(f')([\gamma ])=[H(x_0,-)]\circ \pi _1(f)([\gamma ])\circ [H(x_0,-)]^{-1}$ for all \, $[\gamma ]\in \pi _1(X,x_0)$.
\item{$\bullet $} In the case \, $[H]:f@>\sim >>f':(X,x_0)\to (Y,y_0)$ \, is a pointed $2$-cell \, (\, $[H(x_0,-)]=1_{f(x_0)}:f(x_0)\to f(x_0)=f'(x_0)$) \, then \, $\pi _1(f)=\pi _1(f')$. 
\endproclaim
\demo\nofrills{Proof \ \ } follows from the naturality square\hskip 0.7cm
$\vcenter{\xymatrix{{f(x_0) \ } \ar[r]_-{\sim }^-{[H(x_0,-)]} \ar[d]_-{\Pi (f)([\gamma ])} & { \ f'(x_0)} \ar[d]^-{\Pi (f')([\gamma ])} \\
{f(x_0) \ } \ar[r]^-{\sim }_{[H(x_0,-)]} & { \ f'(x_0)}}}$\vskip -0.55cm
\hfill $\square $
\enddemo

\subhead 2.2. Duality and Invariant Theory\endsubhead

\vskip 0.1cm
{\bf Proposition 2.2.1.} Let $\bold{K}$ be $\bold{Set}$, $\bold{Top}$ or $\bold{Diff}^+$ (spectra of smooth completion (see {\bf 2.3}) of commutative algebras with Zariski topology), $G$ a group. Then there exists
a concrete natural dual adjunction $\xymatrix{\bold{ComAlg}^{op} \ar@/^/[r]^-F & G\text{-}\bold{K} \ar@/^/[l]^-H_-{\top }}$ with $k$ 
($\Bbb R$ or $\Bbb C$), its schizophrenic object, such that $k\in Ob\, G\text{-}\bold{K}$ has trivial action of $G$, and 
$F\circ H:G\text{-}\bold{K}\to G\text{-}\bold{K}$ is a functor ``taking the quotient space generated by the equivalence relation $x\sim y$ iff 
$x,y \in \text{Closure}(\text{the same orbit})$" (it is essentially the orbit space).  \hfill $\square $ 

\vskip 0.2cm
{\bf Definition 2.2.1.} \item{$\bullet $} The adjoint object $\Cal A_X=HX$ for an object $X$ in $G\text{-}\bold{K}$ is called the 
{\bf algebra of invariants}. 
\item{$\bullet $} If $U:G\text{-}\bold{K}\to G\text{-}\bold{K}$ is an endofunctor then $\Cal A_{U(X)}$ is called the {\bf algebra of 
$U$-invariants} of the object $X$. \hfill $\square $ 

\vskip 0.2cm
{\bf Remarks.}
\item{$\bullet $} For $U=(-)^n$, the $n$-fold Cartesian product, $\Cal A_{U(X)}$ is the
{\bf algebra of $n$ point   invariants}.
\item{$\bullet $} For $\bold{K}=\bold{Diff}$, $U=\bold{Jet}^n$, $\bold{Jet}^n(X):=\{j^n_0f \, |\, f\in \bold{Diff}(k,X)\}$, the set of all 
$n$-jets of all maps from $k$ to $X$ at point $0$ (with a certain manifold structure obtained from local trivializations), we get {\bf differential invariants}. 
\item{$\bullet $} The functor $U=\bold{Jet}^{\infty }:\bold{Diff}\to \bold{Diff}^+$ does not fit into the above scheme, but everything is still correct if 
$U:G\text{-}\bold{K}\to G\text{-}\bold{K}_1$ is an extension to $G\text{-}\bold{K}_1$, a category concretely adjoint to $\bold{ComAlg}$.
\item{$\bullet $} $G$ can be, of course, $\bold{Aut}(X)$. 

\vskip 0.2cm
According  to Klein's Erlangen Program, every group acting on a space determines a geometry and, conversely, every geometry hides 
a group of transformations. Properties of geometric objects which are invariant under all transformations are called {\it geometric }
(or invariant  or absolute) for the given $G$-space and a class of geometric objects.

\vskip 0.2cm
{\bf The equivalence problem} \cite{Car1, Car2, Vas, Olv, Gar} consists of a $G$-space $X$ and two ``geometric objects" $S_1, S_2$ of the same type on the space $X$. It is required
to determine if these two objects can be mapped to one another by an element of $G$. An approach is to find a (complete) system of 
invariants of each object. 
 
\vskip 0.2cm
\subhead {} 2.2.1. Classification of covariant geometric objects\endsubhead

By {\it covariant geometric objects} we mean objects like submanifolds, foliations or systems of differential equations, i.e., objects which
behave contravariantly (!)  from the categorical viewpoint and which can be described by a {\bf differential ideal} $I$ ($dI\subset I$) in 
$\Lambda (X)$, the exterior differential algebra of $X$.

\proclaim{\bf Proposition 2.2.1.1} Let $G$ be a Lie-like group (i.e., there exists an algebra of invariant forms on $G$)  
\cite{A-V-L, Car1, Car2}. Then any $G$-equivariant map $\sigma :G\to X$ ($G$ is given with left shift action and $X$ is a left $G$-space) produces a system of invariants of the differential 
ideal $I\subset \Lambda (X)$ (with generators of degree $0$ and $1$) in the following way:
\item{$\bullet $} Take the image $\bar {\Lambda }_{inv}:=Im\, (\Lambda _{inv}(G)\hookrightarrow \Lambda (G)\twoheadrightarrow \Lambda (G)/{\sigma ^*}(I))$,
where $\Lambda _{inv}(G)$ is a subalgebra of left-invariant forms on $G$, $\sigma ^*:\Lambda (X)\to \Lambda (G)$ is the induced map of 
exterior differential algebras, $\sigma ^*(I)$ is the smallest differential ideal in $\Lambda (G)$ containing the image of $I$ under $\sigma ^*$.
\item{$\bullet $} Take the module $\Lambda ^0(G)\cdot \bar {\Lambda }^1_{inv}$ generated by $1$-forms in $\bar {\Lambda }_{inv}$ over $\Lambda ^0(G)$.
There is an open set $\Cal O\subset G$ and a basis $\{\omega ^{\alpha }_{inv}\}_{\alpha \in A}\subset \bar {\Lambda }^1_{inv}$ for the module 
$\Lambda ^0(G)\cdot \bar {\Lambda }^1_{inv}$ restricted to $\Cal O$, i.e., $\forall \, \omega ^i_{inv}\in \bar {\Lambda }^1_{inv}$ 
$\exists \, ! \, \text{ functions } f^i_{\alpha }\in C^{\infty }(\Cal O)$ such that $\omega ^i_{inv}=\sum\limits _{\alpha }f^i_{\alpha }\omega ^{\alpha }_{inv}$.
Form set $J_0:=\{f^i_{\alpha }\}$.
\item{$\bullet $} Take the expansion of differentials $df^i_{\alpha }=\sum\limits _{\beta }f^i_{\alpha \beta }\omega ^{\beta }_{inv}$ (over $\Cal O$). Form the set 
$J_1:=\{f^i_{\alpha \beta }\}$. 
\item{$\bullet $} Continue this process to get $J_2:=\{f^i_{\alpha \beta \gamma }\}, \dots , J_n:=\{f^i_{\alpha _1\dots \alpha _{n+1}}\}\dots $
Form the set $J:=\bigcup\limits _{n}J_n$. Its elements are relative invariants of the differential ideal $I\subset \Lambda (X)$.
\item{$\bullet $} Take the algebra $\Cal A_J\subset C^{\infty }(\Cal O)$, generated by $J$, and take its smooth completion $\overline {\Cal A_J}$ (see {\bf 3.2}). 
Then the ideal $\xymatrix{\bold{Rel}(\Cal A_J) \ar@{>->}[r] & \overline {\bold{Alg}(J)} \ar@{->>}[r] & \overline {\Cal A_J}}$, of all relations of $\Cal A_J$, gives absolute invariants of 
the differential ideal $I\subset \Lambda (X)$, 
where $\overline {\bold{Alg}(J)}$ is the smooth completion of the free algebra generated by $J$.
\endproclaim   
\demo\nofrills{Proof} \ \ follows from the diagrams

\hbox{\hskip 2cm
\hbox{$\xymatrix{G \ar[r]^-{l_g} \ar[d]_-{\sigma } & G \ar[d]^-{\sigma } \\
X \ar[r]_-{l_g} & X}$}
\hskip 2cm
\hbox{$\xymatrix{\Lambda _{inv}(G) & \Lambda _{inv}(G) \ar[l]_-{id} \\
\Lambda (X) \ar[u]^-{\sigma ^*} & \Lambda (X) \ar[l]^-{l_g^*} \ar[u]_-{\sigma ^*}}$}
}

and equations $\omega ^i_{inv}=\sum\limits _{\alpha }f^i_{\alpha }\omega ^{\alpha }_{inv}\, mod(\sigma ^*(I))$.
\hfill $\square $
\enddemo

{\bf Remark.} $G\text{-}\bold{Diff}(G,X)$ is in $1$-$1$-correspondence with all sections of the orbit space $X_G$. So, if $X$ is homogeneous then 
it is exactly the set of all points of $X$ and $\sigma :G\to X=G@>\sim >>G\times \{x_0\}@>1\times i_{x_0}>>G\times X@>\rho >>X$ is a $G$-equivariant 
map corresponding to the point $x_0\in X$, where $\rho $ is the given $G$-action on $X$. 

The following result can be found in \cite{Lap}. Although not well-known, it is a fundamental classification
of analytic geometric objects.

\vskip 0.1cm
\proclaim{\bf Proposition 2.2.1.2 (Exterior differential algebra associated to a group of analytic automorphisms)} Let $X$ be an analytic $n$-dimensional manifold, $\bold{An}(X)$, its group of automorphisms, $H^{\infty }(X):=\{j^{\infty }_0f\, |\, f\in \bold{Diff}(k^n,X), \ X \text{ is analytic}, \text{ Jacobian}(f)\ne 0 \}$, 
the $\infty $-frame bundle over $X$ (with a usual topology and manifold structure). Then there is an exterior differential $k$-algebra $\Lambda _{inv}(H^{\infty }(X))$
of invariant forms on $H^{\infty }(X)$ freely generated by elements of degree 1 obtained by the following process:
\item{$\bullet $} $\omega ^i:=x^i_jdx^j$ are any  ``shift"  forms on $X$
\item{$\bullet $} $\omega ^i_{j}$ are the most general solutions of Maurer-Cartan equations $d\omega ^i=\omega ^i_j\wedge \omega ^j$
\item{$\bullet $} $\omega ^i_{jk}$ are the most general solutions of Maurer-Cartan equations $d\omega ^i_j=\omega ^i_k\wedge \omega ^k_j+\omega ^i_{jk}\wedge \omega ^k$
\item{$\bullet $} $\omega ^i_{jkl}$, $\cdots $, $\omega ^i_{i_1\dots i_n}$, $\cdots $

All forms are symmetric in the lower indices. They characterize the underlying space of $\bold{An}(X)$ uniquely up to analytic isomorphisms. \hfill $\square $
\endproclaim   

\vskip 0.0cm
{\bf Remark.} At each point $x_0\in X$, \, $\omega ^i=0$, \, and forms \, $\bar \omega ^i_{i_1\dots i_n}:=\bigl .\omega ^i_{i_1\dots i_n}\bigr |_{\omega ^i=0}\, , \ \, n\ge 1$, \, are free generators of the exterior differential 
algebra of the {\bf differential group} acting simply transitively on each fiber of $H^{\infty }(X)$.

\vskip 0.2cm
\subhead 2.2.2. Classification of smooth embeddings into a Lie group\endsubhead

\vskip 0.1cm
This is often the last step of smooth classification of geometric objects \cite{Car2, Fin, Kob}. The process of finding   
differential invariants is similar to that in Proposition 2.2.1.1.  The following is essentially in \cite{Vas0, Vas, Lap}.

\proclaim{\bf Proposition 2.2.2.1} For a smooth embedding $f:X\to G$ of a smooth manifold $X$ into a Lie group $G$,  a complete system of 
differential invariants of $f$ can be obtained in the following way:
\vskip 0.1cm
\item{$\bullet $} $Im(f^*:\Lambda ^1_{inv}(G)\to \Lambda (X))$ is locally free, so, has as its basis $\omega _{inv}^i, i=1,\dots ,n$, $n=dim(X)$, near each point.
\item{$\bullet $} Coefficients of linear combinations $\omega ^I_{inv}=\sum \limits _{i=1}^na^I_i\omega ^i_{inv},\, I=n+1,\dots ,dim(G)$, are differential invariants 
of first order (of the map $f$).
\item{$\bullet $} Coefficients of differentials of invariants of first order $da^I_i=\sum \limits _{j=1}^na^I_{ij}\omega _{inv}^j$ are 
differential invariants of second order (of the map $f$).
\item{$\bullet $} $\dots $ Coefficients of differentials of invariants of $(k-1)$ order $da^I_{i_1\dots i_{k-1}}=\sum \limits _{i_k=1}^na^I_{i_1\dots i_k}\omega _{inv}^{i_k}$ 
are differential invariants of order $k$ $\dots $
\vskip 0.2cm
Such calculated invariants characterize an orbit $G\cdot f$ uniquely up to ``changing the parameter space" $X@>\sim >>X'$.  \hfill  $\square $
\endproclaim

\subhead {\bf 2.3. Tangent functor for smooth algebras}\endsubhead

\vskip 0.1cm
This is an example of the dual construction for the main functor of Differential Geometry (which suggests how it can be extended over spectra of commutative algebras).

Let $T:\bold{Diff}\to \bold{Diff}$ be the {\bf tangent} functor on the category of real $\infty $-smooth manifolds. 
In local 
coordinates it is of the form $\cases X\to TX:(x^i)\to (x^i,\xi ^j) & X\in Ob\, \bold{Diff}\\
f\to Tf:(f^i(x))\to (f^i(x),\frac {\partial f^j}{\partial x^k}\, \xi ^k) & f\in Ar\, \bold{Diff}
\endcases $
\vskip 0.2cm
$\bold{Diff}\hookrightarrow \Bbb R\text{-}\bold{Alg}^{op}$ is a subcategory of the opposite of 
the category of real commutative algebras. 
Working in $\bold{Diff}$, it is hard (if possible at all) to give a coordinate-free characterization of $T$. The question is 
how to characterize the image in  $\Bbb R\text{-}\bold{Alg}$?
\vskip 0.2cm
{\bf Definition 2.3.1.} Let $\Cal A\in Ob\, \Bbb R\text{-}\bold{Alg}$.  
\item{$\bullet $} $\rho :\Cal A\to \bold{Top}(\bold{Spec}_{\Bbb R}(\Cal A),\Bbb R)$
is called {\bf functional representation} homomorphism of $\Cal A$, where $\bold{Spec}_{\Bbb R}(\Cal A)=\Bbb R\text{-}\bold{Alg}(\Cal A,\Bbb R)$ has the initial 
topology with respect to all functions $\rho (a), \ a \in \Cal A$, \ $\rho (a)(f):=ev(f,a):=|f|(a)$. 
\item{$\bullet $} $\Cal A$ is called {\bf smooth} if $\forall a_1,a_2,...,a_n\in \Cal A$ and 
$\forall f:\Bbb R^n\to \Bbb R\in C^{\infty }(\Bbb R^n)$ the composite \newline
$f\, \circ <\rho (a_1),\rho (a_2),...,\rho (a_n)>\ \in \ \text{Im}\, (\rho )$. \hfill $\square $
\vskip 0.2cm
Denote by $\Bbb R\text{-{\tt Sm}-}\bold{Alg}\hookrightarrow \Bbb R\text{-}\bold{Alg}$ full subcategory of smooth algebras.
\vskip 0.2cm
\proclaim{\bf Lemma 2.3.1} $\Bbb R\text{-{\tt Sm}-}\bold{Alg}\hookrightarrow \Bbb R\text{-}\bold{Alg}$ is a reflective subcategory, 
i.e. the inclusion has a left adjoint $\text{\tt Sm}:\Bbb R\text{-}\bold{Alg}\to \Bbb R\text{-{\tt Sm}-}\bold{Alg}$, {\bf smooth completion} of 
$\Bbb R$-algebras.
\endproclaim
\demo{Proof} Just take for each $\Bbb R$-algebra $\Cal A$ $\Bbb R$-algebra $\text{\tt Sm}(\Cal A)$ of all terms $\{f(a_1,...,a_n)\, |\, f:\Bbb R^n\to \Bbb R,\ a_1,...,a_n\in \Cal A\}$
(all smooth operations are admitted). Each morphism $f$ from an $\Bbb R$-algebra $\Cal A$ to a smooth algebra $\Cal B$ is uniquely 
extendable to $\tilde  f:\text{\tt Sm}(\Cal A)\to \Cal B$. \hfill $\square $
\enddemo

Let $\bold{Sym}\text{-}\bold{Alg}$ be the category of symmetric partial differential algebras. $Ob\, (\bold{Sym}\text{-}\bold{Alg})$ are 
graded commutative algebras over commutative $\Bbb R$-algebras with a differential $d:\Cal A^0\to \Cal A^1$ of degree $1$  
determined only on elements of degree $0$ ($d$ is $\Bbb R$-linear and satisfies the Leibniz rule). $Ar\, (\bold{Sym}\text{-}\bold{Alg})$ are 
graded degree $0$ algebra homomorphisms which respect $d$. 

\proclaim{\bf Lemma 2.3.2} There is an adjunction $\vcenter{\xymatrix{\Bbb R\text{-}\bold{Alg} \ar@/^/[r]^-{\bold{Sym}}_-{\bot } & \bold{Sym}\text{-}\bold{Alg} \ar@/^/[l]^-{p_0}}}$, \newline 
where: $p_0$ is the projection onto the $0$-component \ $\cases p_0(\Cal A):=\Cal A^0 & \\
p_0(\Cal A@>f>>\Cal B):=(\Cal A^0@>f^0>>\Cal B^0) & \endcases $, \vskip 0.2cm
$\bold{Sym}$ is the functor forming the graded symmetric algebra over the module of differentials of the given algebra \vskip 0.2cm
$\cases \bold{Sym}(\Cal C):=\bold{Sym}(\Lambda ^1(\Cal C)) & \\
\bold{Sym}(\Cal C@>h>>\Cal D):=(\bold{Sym}(\Cal C)@>\tilde h>>\bold{Sym}(\Cal D)) & \\
\tilde h\bigl (\sum c_{i^1\dots i^k}(dc_1)^{i_1}\cdots (dc_k)^{i_k}\bigr ):=\sum h(c_{i^1\dots i^k})(dh(c_1))^{i_1}\cdots (dh(c_k))^{i_k} & \endcases $ \vskip -0.35cm \hfill $\square $
\endproclaim

\vskip 0.2cm
\proclaim{\bf Lemma 2.3.3} \ $\vcenter{\xymatrix{\Bbb R\text{-}\bold{Alg} \ar[r]^-{\text{\tt Sm}} \ar[dr]_-{\bold{Spec}_{\Bbb R}} & \Bbb R\text{-{\tt Sm}-}\bold{Alg} \ar[d]^-{\bold{Spec}_{\Bbb R}} \\
 & \bold{Top}}}$ \ (smooth completion does not change spectrum).
\endproclaim
\demo{Proof} $\forall \alpha :\Cal A\to \Bbb R$ $\exists ! \text{ an extension } \tilde \alpha :\text{\tt Sm}(\Cal A)\to \Bbb R:f(a_1,\dots ,a_n)\mapsto f(\alpha (a_1),\dots ,\alpha (a_n))$. And conversely,
each such $\tilde \alpha $ is restricted uniquely to $\alpha $. Initial topology on $\Bbb R\text{-}\bold{Alg}((\text{\tt Sm})(\Cal A),\Bbb R)$ 
does not change because new functions are functionally (continuously) dependent on old ones.    \hfill $\square $ 
\enddemo

{\bf Remark.} With the Zariski topology on spectra, the  smooth completion yields the same set with a weaker topology. 
For \, $C^{\infty }(X),\, X\in Ob\, \bold{Diff}$ \, the Zariski and initial topologies coincide.  

\vskip 0.2cm
\proclaim{\bf Proposition 2.3.1} \,  
$\bullet $ The tangent functor $T:\Bbb R\text{-{\tt Sm}-}\bold{Alg}\to \Bbb R\text{-{\tt Sm}-}\bold{Alg}$ is equal to the composite
$\Bbb R\text{-{\tt Sm}-}\bold{Alg}\hookrightarrow \Bbb R\text{-}\bold{Alg}@>\bold{Sym}>> \bold{Sym}\text{-}\bold{Alg}@>U>> \Bbb R\text{-}\bold{Alg}@>\text{\tt Sm}>> \Bbb R\text{-{\tt Sm}-}\bold{Alg}$, 
where $U$ forgets the differential $d$ and grading.
\item{$\bullet $} To the canonical projection $\vcenter{\xymatrix{TX \ar[d]_-{p_X} \\ X}}$ there corresponds a canonical embedding 
$\vcenter{\xymatrix{T(C^{\infty }(X))  \\ C^{\infty }(X) \ar[u]^-{i_{C^{\infty }(X)}}}}$.
\endproclaim
\demo{Proof} 
\item{$\bullet $} If $X\in Ob\, \bold{Diff} $ \ $TX\sim \bold{Spec}_{\Bbb R}(U\circ \bold{Sym}(C^{\infty }(X)))\sim \bold{Spec}_{\Bbb R}(\text{\tt Sm}\hskip 0.05cm\circ U\circ \bold{Sym}(C^{\infty }(X)))$.                                   
\item{$\bullet $} This is immediate.                 \hfill $\square $
\enddemo
\vskip 0.2cm
{\bf Remark.} It is reasonable to define $T$ on $\Bbb R\text{-}\bold{Alg}$ as $T:=U\hskip -0.05cm\circ \bold{Sym}$ and transfer 
it to spectra via duality $\xymatrix{\Bbb R\text{-}\bold{Alg}^{op} \ar@/^/[r]_-{\bot }^-{F} & \bold{Spec}_{\Bbb R} \ar@/^/[l]^-{G}}$ (as $F\circ T^{op}\circ G$).

\head {\bf 3. Representable $\infty $-functors} \endhead

{\bf Definition 3.1.} $\infty $-categories $L$ and $L'$ are {\bf equivalent} if $L\sim L'$ in {\bf $\infty $-CAT}. \hfill $\square $

If equivalence $L\sim L'$ is given by functors $\xymatrix{L \ar@/^2ex/[r]^-{F}_-{\sim } & L' \ar@/^2ex/[l]^-{G}_-{\sim }}$ then  
$\forall a\in L^0 \  a\sim G\circ F(a)$,  $\forall b\in L^{'0}$\linebreak
$b\sim F\circ G(b)$  naturally in $a$ and $b$.

\vskip 0.2cm
{\bf Definition 3.2.} $\infty $-functor $F:L\to L'$ is (weakly)
\item{$\bullet $} {\bf faithful} if $\forall a,a'\in L^0$ $\forall f^n,g^n\in L^n(a,a')$ \, $F(f^n)\sim F(g^n)\Rightarrow f^n\sim g^n$, 
\item{$\bullet $} {\bf full} if $\forall a,a'\in L^0$ $\forall h^n\in L^{'n}(F(a),F(a'))$ \, $\exists f^n\in L^n(a,a')$ such that $F(f^n)\sim h^n$,
\item{$\bullet $} {\bf surjective on objects} if $\forall b\in L^{'0}$ \ $\exists a\in L^0$ such that $F(a)\sim b$.
\hfill $\square $
\vskip 0.2cm
Unlike first order equivalence, there is no simple criterion of higher order equivalence. 
\proclaim{\bf {Proposition 3.1}} If the functor $\xymatrix{L \ar[r]^-{F}_-{\sim } & L'}$ is an equivalence then $F$ is (weakly) faithful, full and 
surjective on objects.
\endproclaim
\demo{Proof} "$\Rightarrow $" \ Regard the diagram (where $G$ is a quasiinverse of $F$)
$$\xymatrix{a \ar@/^2ex/[rr]^-{e^n\rho _a}_-{\sim } \ar[d]_-{f^n} && G\circ F(a) \ar@/^2ex/[ll]^-{e^n\theta _a}_-{\sim } \ar[d]^-{G(F(f^n))} \\
a' \ar@/^2ex/[rr]^-{e^n\rho _{a'}}_-{\sim } && G\circ F(a') \ar@/^2ex/[ll]^-{e^n\theta _{a'}}_-{\sim }
}$$

where: $f^n\in L^n(a,a')$, $e^n\rho _{a}\in L^n(a,G(F(a)))$, $e^n\theta _a\in L^n(G(F(a)),a)$, $n\ge 0$.
\vskip 0.2cm
Take $f^n,g^n:a\to a'\in L^n(a,a')$ such that $F(f^n)\sim F(g^n)$. Then $f^n\sim e^n\theta _{a'}\circ _{n+1}G(F(f^n))\circ _{n+1}e^n\rho _a\sim e^n\theta _{a'}\circ _{n+1}G(F(g^n))\circ _{n+1}e^n\rho _a\sim g^n$, \, i.e., $F$ is faithful ($G$ is faithful by symmetry).
\vskip 0.1cm
Take $\alpha ^n:F(a)\to F(a')\in L^{'n}(F(a),F(a'))$. Then $\beta ^n:=e^n\theta _{a'}\circ _{n+1}G(\alpha ^n)\circ _{n+1}e^n\rho _a:a\to a'\in L^n(a,a')$ is such that $G(F(\beta ^n))\sim G(\alpha ^n)$. So, $F(\beta ^n)\sim \alpha ^n$ because $G$ is faithful. Therefore, $F$ is full ($G$ is full by symmetry).
\vskip 0.1cm
$F$ and $G$ are obviously surjective on objects. \hfill $\square $
\enddemo

\vskip 0.1cm
{\bf Remark.} The inverse direction "$\Leftarrow $" for the above proposition works only partially. Namely, 
for each $b\in L^{'0}$ choose $G(b)\in L^0$ and equivalence 
$\xymatrix{b \ar@/^2ex/[rr]^-{\rho _b}_-{\sim } && F(G(b)) \ar@/^2ex/[ll]_-{\sim }^-{\theta _b}}$ 
(which is possible since $F$ is surjective on objects), moreover, if 
$b=F(a)$ choose $G(b)=a, \ \rho _b=eb, \ \theta _b=e(F(G(b)))=eb$. For each $f^n:b\to b'\in L^{'n}(b,b')$ choose an element 
$G(f^n)\in  L^n(G(b),G(b'))$ such that $e^n\rho _{b'}\circ _{n+1}f^n\circ _{n+1}e^n\theta _b\sim F(G(f^n))$ 
(which is possible since $F$ is fully faithful). Then $G:L'\to L$ is obviously a (weak) functor. $a=G(F(a))$ is natural in $a$ 
by construction, but $b\sim F(G(b))$ is natural in $b$ for only first order arrows $\rho _b, \, \theta _b$ presenting $\sim $\, . So, $F$ 
should be somehow 'naturally surjective on objects' which does not make sense yet when the functor $G$ is not defined.  \hfill $\square $

\vskip 0.2cm
{\bf Definition 3.3.} An $\infty $-functor $F:L\to L'$ is called
\item{$\bullet $} an {\bf isomorphism} if it is a bijection (on sets $L,L'$) and the inverse map is a functor,
\item{$\bullet $} a {\bf quasiisomorphism} if there exists a functor $G:L'\to L$ such that 
$\forall a^n\in L^n \ \, G(F(a^n))\sim a^n$ \, and \, $\forall b^n\in L^{'n} \ \, F(G(b^n))\sim b^n$, \, $n\ge 0$.  \hfill $\square $

\proclaim{\bf Proposition 3.2} The notions of (functor) isomorphism and quasiisomorphism coincide.
\endproclaim
\demo{Proof} Each isomorphism is a quasiisomorphism. Conversely, if $\xymatrix{L \ar@/^/[r]^-{F} & L' \ar@/^/[l]^-{G} }$ is 
a quasiisomorphism then $\forall a^n\in L^n$, \, $n\ge 0$, \, $G(F(ea^n))\sim ea^n$. So, $d(G(F(ea^n)))=dea^n$, i.e. $G(F(dea^n))=dea^n$ and $G(F(a^n))=a^n$ (instead of $d$, $c$ could be used). 
The same, $\forall b^n\in L^{'n}$, \, $n\ge 0$, \, $F(G(b^n))=b^n$.          \hfill $\square $
\enddemo

\vskip 0.1cm
Denote (quasi)isomorphism (equivalence) relation by \, $\simeq $\, . 

\vskip 0.35cm
\centerline{{\bf Examples} (isomorphic $\infty $-categories)}
\vskip 0.1cm
\item{1.} Assume, $\xymatrix{f^n \ar[r]^-{\simeq }_-{\alpha } & g^n }$ are 
isomorphic elements of degree $n$ (in a strict category $L$) then $L(f^n,f^n)\simeq L(g^n,g^n)$ are isomorphic $\infty $-categories.
Indeed, there is an isomorphism $F:L(f^n,f^n)\to L(g^n,g^n):x\mapsto \alpha *(x*\alpha ^{-1})$,  
where $*$ means a horizontal composite. $F$ is a functor. Its inverse is 
$G:L(g^n,g^n)\to L(f^n,f^n):y\mapsto \alpha  ^{-1}*(y*\alpha )$. 
[For $\alpha $ just an equivalence, it is not true] 
\item{2.} $\infty\text{\bf -CAT}(L(-,a),F)\simeq F(a)$ (see below, the Yoneda Lemma).    \hfill $\square $

\vskip 0.2cm
{\bf Definition 3.4.} Two $n$-modifications $\alpha ^n, \beta ^n:L\to \infty \text{-}\bold{CAT}$, $n\ge 0$, are called 
{\bf quasiequivalent} of depth $k$, \, $0\le k\le n+1$, (denote it by $\alpha ^n\approx _k\beta ^n$) if their corresponding components 
are quasiequivalent of depth $k-1$, \, i.e. $\forall a\in L^0$ \, $\alpha ^n(a)\approx _{k-1}\beta ^n(a)$. 
$\approx _0$ means $\sim $ by definition. [In other words, $\alpha ^n\approx _k\beta ^n$ if all their components 
of components on depth $k$ are equivalent, i.e. $\alpha ^n\approx _0\beta ^n$ if they are equivalent $\alpha ^n\sim \beta ^n$;
$\alpha ^n\approx _1\beta ^n$ if their components are equivalent $\forall \, a\in L^0 \ \alpha ^n(a)\sim \beta ^n(a)$; 
$\alpha ^n\approx _2\beta ^n$ if components of all components are equivalent; etc.]. 
If $\alpha ^n, \beta ^n:L\to L'$ are proper $n$-modifications (living in $\infty \text{-}\bold{CAT}$) for them only
$\approx _0$ and $\approx _1$ make sense.        \hfill   $\square $

\vskip 0.2cm
\proclaim{\bf Lemma 3.1} 
\item{$\bullet $} $\approx _k$ is an equivalence relation.
\item{$\bullet $} $\approx _{k_1}\ \Rightarrow \ \approx_{k_2}$ if $k_1\le k_2$.
\item{$\bullet $} If $\alpha ^n\approx _k\beta ^n$ then $d\alpha ^n=d\beta ^n$, \, $c\alpha ^n=c\beta ^n$.
\item{$\bullet $} If $(L_1,\approx _{k_1})$, $(L_2,\approx_{k_2})$ are two $\infty $-categories 
(not necessarily proper, i.e. living in $\infty \text{-}\bold{CAT}$) for which given equivalence relations make sense for 
all elements, and $F:L_1\to L_2$, $G:L_2\to L_1$ are maps (not necessarily functors) such that $\forall \, l_1\in L_1$ 
$G(F(l_1))\approx _{k_1}l_1$ and $\forall \, l_2\in L_2$ $F(G(l_2))\approx _{k_2}l_2$, \, and $F,G$ both preserve $d$ (or $c$) then 
$F,G$ are bijections inverse to each other.
\item{$\bullet $} For $L,L'\in Ob\, (\infty \text{-}\bold{CAT})$ and $a\in L^0$ the map 
$ev_{a}:\infty \text{-}\bold{CAT}(L,L')\to L':f^n\mapsto f^n(a)$ is a strict functor. [Similar statement holds when $L,L'$ are 
not proper, e.g. $\infty \text{-}\bold{CAT}$, but we need to formulate it for a bigger universe containing $\infty \text{-}\bold{CAT}$] 
\endproclaim
\demo{Proof} The first two statements are obvious. The third one follows from the fact $x\sim y \, \Rightarrow \, dx=dy, \ cx=cy$ \, and that 
$d,c$ are taken componentwise. The fourth statement follows by the same argument as in the proof of proposition 1.3.2. The last 
statement holds because, again, all operations in $\infty \text{-}\bold{CAT}(L,L')$ are taken componentwise.  \hfill  $\square $
\enddemo 

\vskip 0.0cm
{\bf Remark.} For the proof of the Yoneda lemma, a double evaluation functor is needed. For two functors $F,G:L\to \infty \text{-}\bold{CAT}$ 
take the restriction of the evaluation functor $ev_a$ on the hom-set between $F$ and $G$, i.e. 
$ev_{a\, F,G}:\infty \text{-}\Bbb {CAT}(L,\infty \text{-}\bold{CAT})(F,G)\to \infty \text{-}\bold{CAT}(F(a),G(a)):f^n\mapsto f^n(a)$, 
where $\infty \text{-}\Bbb {CAT}$ is a bigger (and weaker) universe containing $\infty \text{-}\bold{CAT}$ as an object. Now, take a second 
evaluation functor $ev_x:\infty \text{-}\bold{CAT}(F(a),G(a))\to G(a):g^n\mapsto g^n(x)$, \, $x\in (F(a))^0$. Then the double 
evaluation functor is the composite $ev_x\circ _1ev_{a\, F,G}:\infty \text{-}\Bbb {CAT}(L,\infty \text{-}\bold{CAT})(F,G)\to G(a):f^n\mapsto f^n(a)(x)$.
It is a strict functor.      \hfill    $\square $

\vskip 0.2cm
$\infty \text{-}\bold{CAT}$-valued functors, natural transformations and modifications live now in a bigger universe 
$\infty \text{-}\Bbb {CAT}$, and we do not yet have for them appropriate definitions. 

\vskip 0.2cm
{\bf Definition 3.5.} $\infty \text{-}\bold{CAT}$-valued functors, natural transformations and modifications are introduced in a 
similar way as the usual ones by changing all occurrences of $\sim $ with (one degree weaker relation) $\approx _1$, i.e.
\item{$\bullet $} a map $F:L\to \infty \text{-}\bold{CAT}$ of degree $0$ is a {\bf functor} if $F$ strictly preserves $d$ and $c$,
$Fdx=dFx, \ Fcx=cFx$, and weakly up to $\approx _1$ preserves $e$ and composites, $Fex\approx _1eFx$, $F(x\circ _ky)\approx _1F(x)\circ _kF(y)$, 
\item{$\bullet $} For a given sequence of two functors $F,G:L\to \infty \text{-}\bold{CAT}$, $\dots $, two $(n-1)$-modifications 
$\alpha ^{n-1}_1, \alpha ^{n-1}_2:\alpha ^{n-2}_1\to \alpha ^{n-2}_2$ strict (or weak) {\bf $n$-modification} 
$\alpha ^{n}:\alpha ^{n-1}_1\to \alpha ^{n-1}_2$ is a map $\alpha ^{n}:L^0\to \infty \text{-}\bold{CAT}^{n+1}$ such that
$\forall a,b\in L^0$ $\alpha ^n(b)*F(-)\approx _1G(-)*\alpha ^n(a):L^{\ge n}(a,b)\to L^{'\ge n}(F(a),G(b))$ (components of values 
of functors are equivalent).      \hfill   $\square $

\vskip 0.2cm
{\bf Definition 3.6.} A covariant (contravariant) functor $F:L\to \infty \text{\bf -CAT}$ is 
\item{$\bullet $} {\bf weakly representable} if $\exists a\in L^0$
such that $L(a,-)\sim F$ \, ($L(-,a)\sim F$). It means there is an equivalence of two $\infty $-categories $L(a,b)\sim F(b)$ \, ($L(b,a)\sim F(b)$) natural in $b$,
\item{$\bullet $} {\bf strictly representable} if there exists $ a\in L^0$ such that $L(a,-)\simeq F$ \, ($L(-,a)\simeq F$), 
i.e. $\forall b\in L^0$ \, $\exists $ an isomorphism \, $L(a,b)\simeq F(b) \ \, (L(b,a)\simeq F(b))$ natural in $b$.   \hfill  $\square $

\vskip 0.25cm
\proclaim{\bf Lemma 3.2} For given representable $L(-,a):L^{op}\to \infty \text{-}\bold{CAT}$ and functor 
$F:L^{op}\to \infty \text{-}\bold{CAT}$
\item{$\bullet $} all natural transformations  $\tau ^0:L(-,a)\to F$  are of the form  $\forall \, b\in Ob\, L$ the $b$-component is a functor $\tau ^0_b:L(b,a)\to F(b)$,  $\tau ^0_b(f^m)\sim F(f^m)(\tau ^0_a(ea))$,  $f^m\in L^m(b,a)$,
\item{$\bullet $} all $n$-modifications  $\tau ^n:L(-,a)\to F$, $n\ge 1$, are of the form  $\forall \, b\in Ob\, L$ the $b$-component is a
$(n-1)$-modification $\tau ^n_b:L(b,a)\to F(b)$, $\tau ^n_b(f^0)\sim F(f^0)(\tau ^n_a(ea))$, $f^0\in L^0(b,a)$.
\endproclaim
\demo\nofrills{Proof \ }:  follows from the naturality square \hskip 0.35cm$\vcenter{\xymatrix{a & L(a,a) \ar[r]^-{\tau ^n_a} \ar[d]_-{L(f^m,a)} & F(a) \ar[d]^-{F(f^m)} \\
b \ar[u]^-{f^m} & L(b,a) \ar[r]_-{\tau ^n_b} & F(b)}}$\hskip 0.35cm$n\ge 0$   \hfill  $\square $
\enddemo

\proclaim{\bf Lemma 3.3} For a given $n$-cell $\beta ^n\in (F(a))^n$, $n\ge 0$, $n$-modification $\tau ^n:L(-,a)\to F$ such that
$\tau ^n_a(ea)=\beta ^n$ exists and unique up to \, $\approx _2$\, .
\endproclaim
\demo{Proof:} Uniqueness follows from lemma 1.3.2, existence from the definition of $n$-modification 
$\tau ^n_b(f^m):=F(f^m)(\beta ^n)$ (for $n>0$, $m=0$ only) and the naturality square showing correctness of the definition \hskip 1cm
$\vcenter{\xymatrix{b & L(b,a) \ar[d]_-{L(g^k,a)} \ar[r]^-{\tau ^n_b} & F(b) \ar[d]^-{F(g^k)}\\
c \ar[u]^-{g^k} & L(c,a) \ar[r]_-{\tau ^n_c} & F(c)}}$  \vskip 0.0cm
\hskip 0.0cm$\bigl (\mu _{c,b,a}(f^m,g^k):=\mu _{c,b,a}(e^{max(m,k)-m}f^m,e^{max(m,k)-k}g^k)\bigr )$                   \hfill $\square $
\enddemo

\vskip 0.2cm
{\bf Corollary 1.} All $n$-modifications $\tau ^n:L(-,a)\to F$, $n\ge 0$, have strict form 
$\tau ^n_b(f^0)=F(f^0)(\tau ^n_a(ea))$, \, $f^0\in L^0(b,a)$.    \hfill   $\square $

\vskip 0.2cm
{\bf Corollary 2 (criterion of representability).} A $\infty \text{-}\bold{CAT}$-valued presheaf $F\hskip -0.035cm:\hskip -0.035cmL^{op}\hskip -0.035cm\to \hskip -0.035cm\infty \text{-}\bold{CAT}$ is
\vskip 0.05cm 
\item{$\bullet $} {\bf strictly representable} (with representing object $a\in L^0$) \, iff \, 
there exists an object \, $\beta ^0\in (F(a))^0$ \, such that \, $\forall \, \gamma ^n\in (F(b))^n, \ n\ge 0$, \, $\exists !$ 
$n$-arrow \, $(f^n:b\to a)\in L^n(b,a)$ \, with \, $\gamma ^n=F(f^n)(\beta ^0)$,
\vskip 0.05cm 
\item{$\bullet $} {\bf weakly representable} (with representing object $a\in L^0$) \, iff \,
there exists an object \, $\beta ^0\in (F(a))^0$ \, such that \, $\forall b\in Ob\, L$ \, the functor \, 
$L(b,a)\to F(b):f^n\mapsto F(f^n)(\beta ^0)$ is an equivalence of categories.
\item{}(Similar statements hold for a covariant presheaf \, $F:L\to \infty \text{-}\bold{CAT}$)
\hfill $\square $ 

\vskip 0.2cm

\proclaim{\bf {Proposition 3.3 (Yoneda Lemma)}} 
For the functor $F:L^{op}\to \infty\text{\bf -CAT}$ and the object $a\in L^0$, there is a strict isomorphism 
$\infty \text{-}\hskip 0.0125cm\Bbb {CAT}(L(-,a),F)\simeq F(a)$ natural in $a$ and $F$.
\endproclaim
\demo{Proof} Strict functoriality of the correspondence $\tau ^n\mapsto \tau ^n_a(ea)$ is straightforward (because it is a 
double evaluation functor). The map $\beta ^n\mapsto F(-)(\beta ^n)$ is quasiinverse to the first map (with respect to $\approx _2$ and $=$ 
equivalence relations in $\infty\text{-}\Bbb {CAT}(L(-,a),F)$ and $F(a)$ respectively), and it strictly preserves $d$ and $c$. So, these both 
maps are strict isomorphisms. 

Naturality is given by \hskip 0.5cm$\vcenter{\xymatrix{a & F \ar[d]_-{\alpha ^k} & & \infty \text{-}\Bbb {CAT}(L(-,a),F) \ar[r]^-{\simeq } \ar[d]_-{\infty \text{-}\Bbb {CAT}(L(-,f^m),\alpha ^k)} & F(a) \ar[d]^-{\alpha ^k(f^m)} \\
b \ar[u]^-{f^m} & G & & \infty \text{-}\Bbb {CAT}(L(-,b),G) \ar[r]_-{\simeq } & G(b)}}$  \vskip 0.0cm
$\bigl (\text{where \, } \alpha ^k(f^m):=\mu _{F(a),F(b),G(b)}(e^{max(k,m)-k}\alpha ^k_b,e^{max(k,m)-m+1}F(f^m)), \ k,m\ge 0
\bigr )$
\hfill $\square $
\enddemo

\vskip 0.2cm
{\bf Remark.} The Yoneda lemma for $\infty $-categories is similar to the one for first order categories with the difference that elements 
$\beta ^n\in (F(a))^n$ of degree $n$ now determine  higher degree arrows ($n$-modifications) 
$\beta ^n:L(-,a)\to F$ in a $\infty \text{-}\bold{CAT}$-valued presheaf category.         \hfill    $\square $

\vskip 0.2cm
\proclaim{\bf Proposition 3.4 (Yoneda embedding)} There is a Yoneda embedding $\bold Y\hskip -0.05cm:\hskip -0.05cmL\to \infty\text{\bf -CAT}^{L^{op}}\hskip -0.15cm:\alpha \mapsto L(-,\alpha )$, \ $\alpha \in L$,
which is an extension of the isomorphisms from the Yoneda lemma determined on hom-sets $L(a,b)$, $a,b\in L^0$. The Yoneda embedding preserves and reflects
equivalences $\sim $\, .      
\endproclaim
\demo{Proof} By the Yoneda isomorphism for a given $f^n\in L^n(a,b)$,  the corresponding $n$-modification is 
$L(-,b)(f^n):L(-,a)\to L(-,b)$ which is the same as $L(-,f^n):L(-,a)\to L(-,b)$; i.e. the functor 
$\bold Y\hskip -0.05cm:\hskip -0.05cmL\to \infty\text{\bf -CAT}^{L^{op}}\hskip -0.15cm:\alpha \mapsto L(-,\alpha )$, \ $\alpha \in L$, 
locally coincides with isomorphisms from the Yoneda lemma. By lemma 1.1.3 this functor preserves and reflects equivalences \, $\sim $\, .  \hfill    $\square $
\enddemo 

\vskip 0.2cm
{\bf Remark.} Under the assumption that the  category $\infty \text{-}\bold{CAT}$ of {\bf weak} categories, functors and $n$-modifications
exists, all the above reasons remain essentially the same, i.e. the Yoneda lemma and embedding seem to hold in a weak situation.   \hfill  $\square $

\head 4. {\bf (Co)limits} \endhead

{\bf Definition 4.1.} An {\bf $\infty $-graph} is a graded set $G=\coprod\limits _{n\ge 0}G^n$ with two unary operations 
$d,c:\coprod\limits _{n\ge 1}G^n\to \coprod\limits _{n\ge 0}G^n$ of degree $-1$ such that $d^2=dc, \ c^2=cd$. \hfill $\square $  
\vskip 0.1cm
{\bf Definition 4.2.} An {\bf $\infty $-diagram} $D:G\to L$ from $\infty $-graph $G$ to $\infty $-category $L$ is a function of degree $0$
which preserves operations $d,c$. \hfill $\square $

\proclaim{\bf {Proposition 4.1}} All diagrams from $G$ to $L$, natural transformations, modifications form an $\infty $-category 
\, $\bold{Dgrm}_{G,L}$ \, in the same way as the functor category $\infty \text{-}\bold {CAT}(L',L)$. \hfill $\square $
\endproclaim
For a given object $a\in L^0$ \, the \, {\bf constant diagram} to \, $a$ \, is \, $\Delta (a):G\to L: g\mapsto e^na$ \, if \, $g\in G^n$. 
$\Delta :L\to \bold{Dgrm}_{G,L}$ is an $\infty $-functor.
\vskip 0.1cm

Denote $\{e\}\alpha :=\{\alpha , e\alpha , e^2\alpha ,...,e^n\alpha ,...\}$, \, $\alpha \in L$.

\vskip 0.1cm
{\bf Definition 4.3.} Diagram $D:G\to L$ has
\item{$\bullet $} a {\bf limit} if the functor $\bold{Dgrm}_{G,L}(\Delta (-),D):L^{op}\to \infty \text{\bf -CAT}$ is representable. \newline
If $\xymatrix{\nu :L(-,a) \ar[r]^-{\sim } & \bold{Dgrm}_{G,L}(\Delta (-),D)}$ is an equivalence then \newline
$\nu _a(\{e\}ea)\subset \bold{Dgrm}_{G,L}(\Delta (a),D)$ is called a {\bf limit cone} over $D$, $a$ is its {\bf vertex} (or diagram {\bf limit $lim\, D$}), $\nu _a(ea)$ are its {\bf edges}, $\nu _a(e^ka),\, k>1$, are identities
\item{$\bullet $} {\bf colimit} if functor $\bold{Dgrm}_{G,L}(D,\Delta (-)):L\to \infty \text{\bf -CAT}$ is representable. \newline
If $\xymatrix{\nu :L(a,-) \ar[r]^-{\sim } & \bold{Dgrm}_{G,L}(D,\Delta (-))}$ is the equivalence then \newline
$\nu _a(\{e\}ea)\subset \bold{Dgrm}_{G,L}(D,\Delta (a))$ is called {\bf colimit cocone} over $D$, $a$ is its {\bf vertex} (or diagram {\bf colimit $colim\, D$}), $\nu _a(ea)$ are its {\bf edges}, $\nu _a(e^ka),\, k>1$, are identities \hfill $\square $

\vskip 0.2cm
{\bf Remark.} The conditions on equivalence $\nu $ in the above definition can be strengthened. If it is a (natural) isomorphism then 
(co)limits are called {\bf strict} and as a rule they are different from {\bf weak} ones \cite{Bor1}.

\vskip 0.2cm
\proclaim{\bf Proposition 4.2} For strict (co)limits the following is true  \vskip 0.25cm
\item{$\bullet $} $
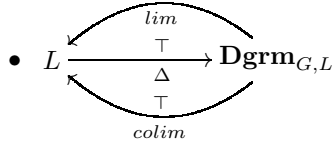
\xymatrix{L \ar[rr]_(0.49){\Delta }^(0.49){\top } & & \bold{Dgrm}_{G,L} \ar@/^5ex/[ll]^-{colim}_-{\top } \ar@/_5ex/[ll]^-{lim}}$ 
\vskip 0.1cm
\item{$\bullet $} Strict right adjoints preserve limits (strict left adjoints preserve colimits).
\endproclaim
\demo{Proof}
\item{$\bullet $} It is immediate from definition 4.3 and proposition 5.1. 
\item{$\bullet $} The argument is the same as for first order categories (see example 13, point {\bf 1.a}) [the essential thing is that
a strict adjunction is determined by (triangle) identities which are preserved under $\infty $-functors].    \hfill $\square $
\enddemo

\vskip 0.2cm
\centerline{\bf Examples}
\vskip 0.2cm
\item{$1.$} (strict binary products in $2\text{-}\bold{Top}$ and $2\text{-}\bold{CAT}$) They coincide with '$1$-dimensional' products. The mediating $2$-cell arrow is given componentwise 
\hskip 1.5cm$\vcenter
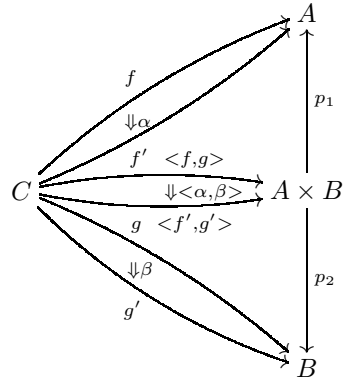
{\xymatrix{ & & & A \\
 & & & \\
C \ar@/^0.6pc/[rrruu]^(0.45){f}_(0.39){\Downarrow \alpha } \ar@/_/[rrruu]_(0.35){f'} \ar@/^/[rrrdd]^(0.35){g}_(0.43){\Downarrow \beta } \ar@/_0.7pc/[rrrdd]_(0.47){g'} \ar@/^/[rrr]^(0.6){<f,g>}_(0.65){\Downarrow <\alpha ,\beta >} \ar@/_/[rrr]_(0.6){<f',g'>} & & & A\times B \ar[uu]_-{p_1} \ar[dd]^-{p_2}\\
 & & & \\
 & & & B}}$
\vskip 0.1cm
\item{$2.$} (``equalizer" of a $2$-cell in $2\text{-}\bold{CAT}$) \cite{Bor1} For a given $2$-cell $\xymatrix{\bold{A} \ar@/^/[r]^-{F}_-{\ \Downarrow \alpha } \ar@/_/[r]_{G} & \bold{B}}$ in $2\text{-}\bold{CAT}$ 
its strict limit is a subcategory $\bold{E}\hookrightarrow \bold{A}$ such that $F(A)=G(A)$ and $\alpha _A=1_{F(A)}:F(A)\to G(A)$ (on objects), and $F(f)=G(f)$ (on arrows).
\vskip 0.1cm
\item{$3.$} (strict and weak pullbacks in $2\text{-}\bold{CAT}$) \cite{Bor1} Let $\Cal P$ be a ``$2$-dimensional" graph 
$1@>x>>0@<y<<2$ with trivial $2$-cells, $F:\Cal P\to 2\text{-}\bold{CAT}$ be a $2$-functor. Then its limit is a pullback
diagram in $2\text{-}\bold{CAT}$ \hskip 0.3cm$\vcenter{\xymatrix{F(1)\times _{F(0)}F(2) \ar[r]^-{p_2} \ar[d]_-{p_1} \ar[dr]^-{p_3} & F(2) \ar[d]^-{F(y)} \\
F(1) \ar[r]_-{F(x)} & F(0)}}$. \hskip 0.3cm When the limit is taken {\bf strictly} $F(1)\times _{F(0)}F(2)$ coincides with 
the ``$1$-dimensional" pullback, 
i.e. $F(1)\times _{F(0)}F(2)\hookrightarrow F(1)\times F(2)$ is a subcategory consisting of objects $(A,B), \ A\in Ob\, F(1), \ B\in Ob\, F(2), \ F(x)(A)=F(y)(B)$ 
and arrows $(f,g), \ f\in Ar\, F(1), \ g\in Ar\, F(2), \ F(x)(f)=F(y)(g)$. When the limit is taken {\bf weakly} 
$F(1)\times _{F(0)}F(2)$ is not a subcategory of product $F(1)\times F(2)$. It consists of $5$-tuples $(A,B,C,f,g)$, $A\in Ob\, F(1)$, 
$B\in Ob\, F(2)$, $C\in Ob\, F(0)$, \ $f:F(x)(A)@>\sim >>C$, \ $g:F(y)(B)@>\sim >>C$ are isomorphisms, with arrows $(a,b,c)$, 
$a:A\to A'$, $b:B\to B'$, $c:C\to C'$ such that $c\circ f=f'\circ F(x)(a)$, $c\circ g=g'\circ F(y)(b)$. 
Projections $p_1, p_2, p_3$ are obvious. The pullback square commutes up to isomorphisms \ $f:F(x)\circ p_1\Rightarrow p_3$, \ $g:F(y)\circ p_2\Rightarrow p_3$.   \hfill  $\square $

\head {\bf 5. Adjunction}\endhead

{\bf Definition 5.1.} The situation $\xymatrix{L\ar@/^/[r]^-{F} & L' \ar@/^/[l]^-{G}_-{\perp }}$ 
(where $L,L'$ are $\infty $-categories, $F,G$ are $\infty $-functors) is called
\item{$\bullet $} {\bf weak $\infty $-adjunction} if there is an equivalence
$L(-,G(+))\sim L'(F(-),+):L^{op}\times L'\to \infty \text{\bf -CAT}$
(i.e. $L(a,G(b))\sim L'(F(a),b)$ natural in $a\in L^0, \, b\in L^{'0}$),
\item{$\bullet $} {\bf strict $\infty $-adjunction} if there is an isomorphism                                            
$L(-,G(+))\simeq L'(F(-),+):L^{op}\times L'\to \infty \text{\bf -CAT}$ 
(i.e. $L(a,G(b))\simeq L'(F(a),b)$ natural in $a\in L^0, \, b\in L^{'0}$).           \hfill $\square $

\proclaim{\bf {Proposition 5.1}} The following are equivalent
\item{1.} $\xymatrix{L\ar@/^/[r]^-{F} & L' \ar@/^/[l]^-{G}_-{\perp }}$ is a strict $\infty $-adjunction
\item{2.} $\forall b\in L^{'0}$ \, $L'(F(-),b)$ is strictly representable \vskip 0.1cm
\item{3.} $\forall a\in L^0$ \, $L(a,G(-))$ is strictly representable
\endproclaim
\demo{Proof} 
\item{$\bullet $}$1. \Longrightarrow  2., 3.$ \, is immediate
\item{$\bullet $}$2. \Longrightarrow  1.$ \, From the criterion of strict representability (see point 1.3) it follows that $\forall \, b\in L^{'0}$ 
there exists a ``universal element" $(\beta ^0_b:F(G(b))\to b)\in L^{'0}(F(G(b)),b)$ such that $\forall \, (f^n:F(c)\to b)\in L^{'n}(F(c),b)$
$\exists \, !$ $n$-arrow $(g^n:c\to G(b))\in L^n(c,G(b))$ with $f^n=\mu _{F(c),F(G(b)),b}(e^n\beta ^0_b,F(g^n))$
\hskip 0.8cm$\vcenter{\xymatrix{G(b) & F(G(b)) \ar[r]^-{e^n\beta ^0_b} & b \\
c \ar@{-->}[u]^-{\exists ! g^n} & F(c) \ar[u]^-{F(g^n)} \ar[ur]_-{\forall f^n} & }}$ \vskip 0.0cm
Consequently, $\forall \, (f^n:b'\to b)\in L^{'n}(b',b)$ the diagram holds 
\hskip 0.0cm$\vcenter{\xymatrix{G(b) & F(G(b)) \ar[r]^-{e^n\beta ^0_b} & b \\
G(b') \ar@{-->}[u]^-{G(f^n)} & F(G(b')) \ar[r]_-{e^n\beta ^0_{b'}} \ar[u]^-{F(G(f^n))} & b' \ar[u]_-{f^n} }\hskip -0.1cm}$
\item{} It shows that assignment $Ob\, L'\ni b\mapsto G(b)\in Ob\, L$ is extendable to a functor $G:L'\to L$ (in an essentially unique way) 
and that $\beta ^0:FG\to 1_{L'}$ is a natural transformation ({\bf counit} \, $\varepsilon $ of the adjunction $F\dashv G$). 
\item{} Isomorphism $\varphi _{c,b}:L'(F(c),b)\to L(c,G(b))$ such that $\vcenter{\xymatrix{F(G(b)) \ar[r]^-{e^n\beta ^0_b} & b \\
F(c) \ar[u]^-{F(\varphi _{c,b}(f^n))} \ar[ur]_-{f^n} & }}$
is natural in $c\in Ob\, L, \, b\in Ob\, L'$ because of the naturality square 
\item{} $\vcenter{\xymatrix{c & b \ar[d]_-{f^n} & & L'(F(c),b) \ar[r]^-{\varphi _{c,b}} \ar[d]_-{L'(F(g^n),f^n)} & L(c,G(b)) \ar[d]^-{L(g^n,G(f^n))} \\
c' \ar[u]^-{g^n} & b' & & L'(F(c'),b') \ar[r]_-{\varphi _{c',b'}} & L(c',G(b'))}}$
\item{} (indeed, $\forall h^n\in L^{'n}(F(c),b)$ \, $G(f^n)*\varphi _{c,b}(h^n)*g^n\sim \varphi _{c',b'}(f^n*h^n*F(g^n))$, where $*$ is the horizontal composite,
since $e^n\beta ^0_{b'}*F(G(f^n)*\varphi _{c,b}(h^n)*g^n)\sim f^n*e^n\beta ^0_b*F(\varphi _{c,b}(h^n))*F(g^n)\sim f^n*h^n*F(g^n)$)
\item{$\bullet $}$3. \Longrightarrow  1.$ \, is similar to \, 2. $\Longrightarrow $ 1. \hfill $\square $
\enddemo

\vskip 0.1cm
{\bf Remark.} The analogous statement for a weak $\infty $-adjunction is not true. In the above proof  ``universal elements" were used in an essential way.       \hfill $\square $

\vskip 0.1cm
{\bf Definition 5.2.} For a given {\bf strict} $\infty $-adjunction $\xymatrix{L\ar@/^/[r]^-{F} & L' \ar@/^/[l]^-{G}_-{\perp }}$
\item{$\bullet $} universal elements $\varepsilon _b:F(G(b))\to b$ representing functors $L'(F(-),b)$ ($b\in Ob\, L'$ is a parameter)
form a natural transformation $\varepsilon :FG\to 1_{L'}$ which is called the {\bf counit} of the adjunction,
\item{$\bullet $} Universal elements $\eta _a:a\to G(F(a))$ representing functors $L(a,G(-))$ ($a\in Ob\, L$ is a parameter) 
form a natural transformation $\eta :1_{L}\to GF$ which is called the {\bf unit} of the adjunction.   \vskip 0.0cm \hfill $\square $

\vskip 0.1cm
{\bf Remark.} For a {\bf weak} $\infty $-adjunction no useful unit and counit exist.   \hfill  $\square $

\proclaim{\bf Proposition 5.2}
\item{$\bullet $} For both weak and strict adjunctions: the
composition of left adjoints is a left adjoint (the composition of right adjoints is a right adjoint).  
\item{$\bullet $} For a weak (strict) adjunction, a right or left adjoint is determined uniquely up to equivalence $\sim $ 
(up to isomorphism $\simeq $).                    
\endproclaim
\demo{Proof} 
\item{$\bullet $} If $\xymatrix{L \ar@/^/[r]^-{F}_-{\perp} & L' \ar@/^/[r]^-{F'}_-{\perp } \ar@/^/[l]^-{G} & L'' \ar@/^/[l]^-{G'}}$ 
then $L''(F'Fl,l'')\sim L'(Fl,G'l'')\sim L(l,GG'l'')$ (composite of natural equivalences). \ [For a strict adjunction the same reason works]
\item{$\bullet $} Assume, $L'(a,G'b)\hskip -0.05cm\sim \hskip -0.05cmL(Fa,b)\hskip -0.05cm\sim \hskip -0.05cmL'(a,Gb)$ are natural equivalences then $L'(-,G'b)\hskip -0.05cm\sim \hskip -0.05cmL'(-,Gb)$ is a natural
transformation (equivalence) natural in $b$. Then, by the Yoneda embedding, $G'b\sim Gb$ naturally in $b$, i.e. $G'\sim G$.
[Again, changing $\sim $ with $\simeq $ still works].              \hfill $\square $                               
\enddemo 

\proclaim{\bf Proposition 5.3} For a strict adjunction $\xymatrix{L \ar@/^/[r]^-{F} & L' \ar@/^/[l]^-{G}_-{\perp } }$ the Kan definition and the 
definition via ``unit-counit" coincide, i.e. the following are equivalent
\item{$\bullet $} $\varphi _{a,b}:L(a,G(b))\simeq L'(F(a),b):\varphi ^*_{a,b}$ natural in $a\in L^0, \, b\in L^{'0}$,
\item{$\bullet $} $\exists $ natural transformations $\eta :1_L\to GF$ and $\varepsilon :FG\to 1_{L'}$ satisfying the triangle 
identities $\varepsilon F\circ _1F\eta =1_F$ and $G\varepsilon \circ _1\eta G=1_G$. 
\endproclaim
\demo{Proof} For a strict adjunction, the same proof as for first order categories works.\vskip 0.1cm
\item{$\bullet $} Universal elements $\eta _a, \varepsilon _b$ for functors $L(a,G(-)),L'(F(-),b)$ mean that they are images of
$1_{F(a)}, 1_{G(b)}$ under functors $\varphi ^*_{a,F(a)},\varphi _{G(b),b}$\, , \, i.e. 
\hskip 0.5cm$\xymatrix{FGFa \ar[r]^-{\varepsilon _{Fa}} & Fa \\
Fa \ar@{-->}[u]^-{F\eta _a} \ar[ur]_-{1_{Fa}} & }$\hskip 0.5cm
$\xymatrix{Gb \ar[r]^-{\eta _{Gb}} \ar[dr]_-{1_{Gb}} & GFGb \ar@{-->}[d]^-{G\varepsilon _b} \\
 & Gb }$\vskip -0.1cm
\item{}\hskip 9.7cm (strict equalities)
\vskip 0.15cm\item{$\bullet $} Define maps $\cases \varphi _{a,b}(f^n):=e^n(\varepsilon _b)\circ _{n+1}F(f^n),  &  f^n\in L^n(a,G(b)) \\
\varphi ^*_{a,b}(g^n):=G(g^n)\circ _{n+1}e^n(\eta _a),   &  g^n\in L^{'n}(F(a),b)   \endcases $
\vskip 0.25cm\item{}They are functors $\cases \varphi _{a,b}:=\varepsilon _b*F(-):L(a,G(b))\to L'(F(a),b) &  \\
\varphi ^*_{a,b}:=G(-)*\eta _a:L'(F(a),b)\to L(a,G(b)) &   \endcases $
\hskip -0.3cmand inverses to each other:\vskip 0.15cm 
\item{}$\varphi ^*_{a,b}(\varphi _{a,b}(f^n))=\varphi ^*_{a,b}(e^n\varepsilon _b\circ _{n+1}F(f^n))=
G(e^n\varepsilon _b\circ _{n+1}F(f^n))\circ _{n+1}e^n\eta _a=e^nG(\varepsilon _b)\circ _{n+1}(GF(f^n)\circ _{n+1}e^n\eta _a)=
e^nG(\varepsilon _b)\circ _{n+1}(e^n\eta _{G(b)}\circ _{n+1}f^n)=e^{n+1}G(b)\circ _{n+1}f^n=f^n$, and similar, 
$\varphi _{a,b}(\varphi ^*_{a,b}(g^n))=g^n$.        
\item{}Naturality (e.g., of $\varphi _{a,b}$) follows from the square 
$\vcenter{\xymatrix{L(a,G(b))  \ar[r]^-{\varphi _{a,b}} \ar[d]_-{L(x^m,G(y^m))} &     L'(F(a),b)  \ar[d]^-{L'(F(x^m),y^m)}  \\
L(a',G(b'))   \ar[r]_-{\varphi _{a',b'}}       &         L'(F(a'),b')    }}$       
\vskip 0.1cm\item{}$\bigl (\varphi _{a',b'}(L(x^m,G(y^m))(f^n))=\varphi _{a',b'}(G(y^m)*f^n*x^m)=
\varepsilon _{b'}*FG(y^m)*F(f^n)*F(x^m)=y^m*\varepsilon _{b}*F(f^n)*F(x^m)=L'(F(x^m),y^m)(\varphi _{a,b}(f^n))$, where $n=0$ or 
$m=0\bigr )$.                    \hfill     $\square $
\enddemo

\vskip 0.35cm
\centerline{\bf Examples of higher order adjunctions}

\vskip 0.2cm
\item{1.} Every usual 1-adjunction $\xymatrix{A \ar@/^/[r] & B \ar@/^/[l]_-{\bot }}$ is an $\infty $-1-adjunction for trivial $\infty $-extensions of $A$ and $B$.
\item{2.} Gelfand-Naimark dual 1-adjunction $\xymatrix{\bold{C^*Alg}^{op} \ar@/^/[r] & \bold{CHTop} \ar@/^/[l]_-{\bot }}$ is extendable to $\infty $-2-adjunction (see {\bf 9}).
\proclaim{\hskip -0.0cm{\rm 3.}\ {\sl Quillen theorem} \cite{Mac}} Let $\bold{\Delta }$ be a category of finite linearly ordered sets, $\bold{Set}^{\bold{\Delta }^{op}}$ the category of simplicial sets, $Ho(\bold{Top}):=(2\text{\bf -Top})^{(1)}$, $Ho(\bold{Set}^{\Delta ^{op}}):=(2\text{\bf -Set}^{\Delta ^{op}})^{(1)}$. Then
$$\xymatrix{\Delta \ar@{{>}->}[dr] \ar@{{>}->}[d]|-{\text{\rm Yoneda}} &  \\
\bold{Set}^{\Delta ^{op}} \ar@{-->}[r] \ar[d] & \bold{Top} \ar@/^2.5ex/[l]_-{\bot} \ar[d]\\
Ho(\bold{Set}^{\Delta ^{op}}) \ar[r] & Ho(\bold{Top}) \ar@/^2.5ex/[l]_-{\bot}
}$$
\endproclaim \hfill $\square $
\vskip 0.0cm 
So, the top adjunction is actually a 2-adjunction (or $\infty $-2-adjunction).
\vskip 0.0cm 
All the above adjunctions are strict.

\head {\bf 6. Concrete duality for $\infty$-categories}\endhead

\vskip 0.0cm
Duality preserves all categorical properties.  It is significant that concrete duality for $\infty $-categories behaves the same as for $1$-categories.

\vskip 0.25cm
{\bf Definition 6.1.} 
\item{$\bullet $} {\bf Duality} is an equivalence $L^{op}\sim L'$.
\item{$\bullet $} A {\bf Concrete duality} over $\Bbb B\hookrightarrow \infty \text{-}\bold{CAT}$ is a duality 
$\xymatrix{L^{op} \ar@/^0.6pc/[r]^-{G}_-{\sim } & L' \ar@/^0.6pc/[l]^-{F}_-{\sim } }$
such that there exist  (faithful) forgetful functors $U:L\to \Bbb B$, \, $V:L'\to \Bbb B$ and objects $\tilde A\in L^0$, \, 
$\tilde B\in L^{'0}$ such that 
\vskip 0.15cm\itemitem{$\bullet $} $U(\tilde A)\sim V(\tilde B)$, 
\itemitem{$\bullet $} \hbox{$V\circ _1G\sim L(-,\tilde A)$, \ $U\circ _1F^{op}\sim L'(-,\tilde B)$} \hskip 1cm
$\vcenter{\hbox{\xymatrix{L^{op} \ar[r]^-{G} \ar[rd]_-{L(-,\tilde A)} & L' \ar[d]^-{V}\\
 & {\Bbb B}}}}$ \hskip 1cm
$\vcenter{\hbox{\xymatrix{L^{'op} \ar[r]^-{F^{op}} \ar[rd]_-{L'(-,\tilde B)} & L \ar[d]^-{U}\\
 & {\Bbb B}}}}$ 
\vskip 0.15cm
\item{} Representing objects $\tilde A\in L^0,\ \tilde B\in L^{'0}$ are called   {\bf dualizing} or {\bf schizophrenic objects} for the given concrete duality\cite{P-Th}.
\vskip 0.15cm\item{} \hskip 0cm[for a {\bf concrete dual adjunction} the definition is similar]           \hfill    $\square $ 

\vskip 0.25cm
\proclaim{\bf {Proposition 6.1} (representable forgetfuls $\Rightarrow $ concrete dual adjunction)} Let $(L,U)$, $(L',V)$ be (weakly) dually adjoint 
$\infty $-categories $\xymatrix{L^{op} \ar@/^/[r]^-{G}_-{\top } & L' \ar@/^/[l]^-{F} }$ 
with representable forgetful functors \ 
$U\sim L(A_0,-):L\to \Bbb B$, \, $V\sim L'(B_0,-):L'\to \Bbb B$ (where $\Bbb B\hookrightarrow \infty \text{-}\bold{CAT}$ is a subcategory). 
\, Then this adjunction is concrete over
$\Bbb B$ with dualizing object $(\tilde A, \tilde B)$, where $\tilde A:=F(B_0)$, \, $\tilde B:=G(A_0)$, i.e. 
\vskip 0.3cm
\item{$\bullet $} $U(\tilde A)\sim V(\tilde B)$
\item{$\bullet $} \hbox{$V\circ _1G\sim L(-,\tilde A)$, \ $U\circ _1F^{op}\sim L'(-,\tilde B)$} \hskip 1cm
$\vcenter{\hbox{\xymatrix{L^{op} \ar[r]^-{G} \ar[rd]_-{L(-,\tilde A)} & L' \ar[d]^-{V}\\
 & {\Bbb B}}}}$ \hskip 1cm
$\vcenter{\hbox{\xymatrix{L^{'op} \ar[r]^-{F^{op}} \ar[rd]_-{L'(-,\tilde B)} & L \ar[d]^-{U}\\
 & {\Bbb B}}}}$
\endproclaim

\demo{Proof}
\vskip 0.1cm
\item{$\bullet $} $U(\tilde A)=UF(B_0)\sim L(A_0,FB_0)\sim L'(B_0,GA_0)\sim VGA_0=V\tilde B$
\item{$\bullet $} $VG(-)\sim L'(B_0,G(-))\sim L(-,FB_0)=L(-,\tilde A)$ \, (and similar, $UF(-)\sim L'(-,\tilde B)$)        \hfill $\square $
\enddemo

\vskip 0.2cm
{\bf Remarks.}
\item{$\bullet $} Concrete duality as above should be called {\bf weak}. {\bf Strict} variants of definition 6.1 and 
proposition 6.1 also exist (by changing $\sim $ to isomorphism $\simeq $ and weak dual adjunction to the strict one).
\item{$\bullet $} (Weak or strict) concrete duality (dual adjunction) is given essentially by hom-functors 
which admit lifting along forgetful functors (to obtain proper values). Representing objects of these functors have equivalent 
(or isomorphic) underlying objects.      
\item{$\bullet $} For the usual 1-dimensional categories $\Bbb B=\bold{Set}\hookrightarrow \infty \text{-}\bold{CAT}$ 
($\infty $-1-subcategory). For dimension $n$, as a rule, $\Bbb B=n\text{-}\bold{Cat}\hookrightarrow \infty \text{-}\bold{CAT}$ 
($\infty \text{-}n$-subcategory of small $(n-1)$-categories).    \hfill    $\square $

\subhead 6.1. Natural and non natural duality \endsubhead

\vskip 0.2cm
{\bf Definition 6.1.1.} \item{$\bullet $} For hom-set $L(A,\tilde A)$ and element $(x:A_0\to A)\in L^0(A_0,A)$ the {\bf evaluation functor} 
at the point $x$ is \, $ev_{A,x}:=L(x,\tilde A):L(A,\tilde A)\to L(A_0, \tilde A)$ \, ($ev_{A,x}\in \Bbb B^1\hookrightarrow \infty \text{-}\bold{CAT}^1$).
\item{}Similarly, the {\bf evaluation $(n-1)$-modification} $ev_{A,x^n}$, $n=1,2,\dots $, for $x^n\in L^n(A_0,A)$ is $L(x^n,\tilde A)\in \Bbb B^n(L(A,\tilde A),L(A_0,\tilde A))$.
\item{$\bullet $} For a forgetful functor $V:L'\to \Bbb B$ an arrow $f^n:V(Y)\to V(Y')\in \Bbb B^n(V(Y),V(Y'))$ is called 
an $L'${\bf -arrow} if $\exists \Phi ^n:Y\to Y'\in L^{'n}(Y,Y')$ such that $V(\Phi ^n)=f^n$.
\item{$\bullet $} A lifting of hom-functor $V\circ G\sim L(-,\tilde A)$ is called 
{\bf initial} \cite{A-H-S} if $\forall A\in L^0\ \forall Y\in L^{'0}\ \forall f^n:V(Y)\to L(A,\tilde A)\in \Bbb B^{n}(V(Y),L(A,\tilde A))$ $f^n$ 
is an $L'$-arrow iff $\forall (x^n:A_0\to A)\in L^n(A_0,A)$ 
$ev_{A,x^n}\circ _{n+1}f^n:V(Y)\to L(A_0,\tilde A)\in \Bbb B^n(V(Y),L(A_0,\tilde A))$ is an $L'$-arrow.
\item{$\bullet $} If liftings of hom-functors $V\circ G\sim L(-,\tilde A),\ U\circ F\sim L'(-,\tilde B)$ are both initial, 
then the concrete dual adjunction $\xymatrix{L^{op} \ar@/^/[r]^-{G}_-{\top } & L' \ar@/^/[l]^-{F} }$, if it exists, 
is called {\bf natural} \cite{Hof, P-Th}, and otherwise, non-natural. \hfill $\square $

\vskip 0.2cm
Even if $U\tilde A\sim V\tilde B$ and $\forall A\in L^0, B\in L^{'0}$ $\Bbb B$-objects $L(A,\tilde A), \, L'(B,\tilde B)$ can be 
lifted to $L', L$, the hom-functors $L(-,\tilde A),\, L'(-,\tilde B)$ need not (which happens only if lifting of the
assignments $A\mapsto L(A,\tilde A),\, B\mapsto L'(B,\tilde B)$ can be extended functorially over all cells). 

\vskip 0.15cm
We introduce the following concept. {\bf The initial lifting condition for the evaluation cones} \newline 
$\{ev_{A,x^n}\in \Bbb B^n(L(A,\tilde A),L(A_0,\tilde A))\}_{x^n\in L^n(A_0,A)}^{\, n\, \in \, \Bbb N}, \, 
\{ev_{B,y^n}\in \Bbb B^n(L'(B,\tilde B),L'(B_0,\tilde B))\}_{y^n\in L^{'n}(B_0,B)}^{\, n\, \in \, \Bbb N}$ consists of the following requirements:
\item{$\bullet $} hom-categories of the form $L(A,\tilde A),\, L'(B,\tilde B)\in Ob\, (\Bbb B)$ can be lifted to $L',L$
\item{$\bullet $} evaluation cones \newline 
$\{ev_{A,x^n}\hskip -0.05cm\in \hskip -0.05cm\Bbb B^n(L(A,\tilde A),L(A_0,\tilde A))\}_{x^n\in L^n(A_0,A)}^{\, n\, \in \, \Bbb N}, \, 
\{ev_{B,y^n}\hskip -0.05cm\in \hskip -0.05cm\Bbb B^n(L'(B,\tilde B),L'(B_0,\tilde B))\}_{y^n\in L^{'n}(B_0,B)}^{\, n\, \in \, \Bbb N}$
can be lifted to $\{ev_{A,x^n}\in L^{'n}(G(A),\tilde B)\}_{x^n\in L^n(A_0,A)}^{\, n\, \in \, \Bbb N}\, , \ \{ev_{B,y^n}\in L^n(F(B),\tilde A)\}_{y^n\in L^{'n}(B_0,B)}^{\, n\, \in \, \Bbb N}$ \, in $L',L$   
\item{$\bullet $} $\forall f^n\in \Bbb B^n(VX,L(A,\tilde A))$ $f^n$ is $L'$-arrow iff $\forall x^n\hskip -0.05cm\in \hskip -0.05cmL^n(A_0,A)$ \, $\mu (ev_{A,x^n},f^n)\hskip -0.05cm\in \hskip -0.05cm\Bbb B^n(VX,L(A_0,\tilde A))$ is $L'$-arrow
(and, symmetrically, $\forall g^n\in \Bbb B^n(UY,L'(B,\tilde B))$ $g^n$ is an $L$-arrow iff $\forall y^n\hskip -0.05cm\in \hskip -0.05cmL^{'n}(B_0,B)$ \, $\mu (ev_{B,y^n},g^n)\hskip -0.05cm\in \hskip -0.05cm\Bbb B^n(UY,L'(B_0,\tilde B))$ is an $L$-arrow) \hfill $\square $

\vskip 0.1cm
In the following proof, we denote lifted evaluation maps by $ev_{A,x}$ 
(or something similar) and underlying evaluation maps in $\Bbb B$ by $|ev_{A,x}|$.

\vskip 0.2cm
\proclaim{\bf Proposition 6.1.1} If two strict $\infty $-categories $L,L'$ concrete over 
$\Bbb B\hookrightarrow \infty \text{-}\bold{CAT}$ with representable (strictly faithful) forgetful functors
$U=L(A_0,-), \ V=L'(B_0,-)$ have objects $\tilde A\in L^0, \, \tilde B\in L^{'0}$ such that
\item{$\bullet $} $U\tilde A\sim V\tilde B$
\item{$\bullet $} the hom-functors $L(-,\tilde A):L^{op}\to \Bbb B, \  L'(-,\tilde B):L^{'op}\to \Bbb B$ satisfy 
the {\bf initial lifting condition for evaluation cones} \newline 
then there exists a  natural {\bf strict} concrete dual adjunction \hskip 0.3cm 
$\xymatrix{L^{op} \ar@/^/[r]^-{G}_-{\top } & L' \ar@/^/[l]^-{F} }$ \hskip 0.3cm 
$L(A,FB)\underset {\text{nat. iso}}\to {\simeq }L'(B,GA)$\hskip 0.5cm
$\vcenter{\hbox{\xymatrix{L^{op} \ar[r]^-{G} \ar[rd]_-{L(-,\tilde A)} & L' \ar[d]^-{V}\\
 & {\Bbb B}}}}$ \hskip 0.5cm
$\vcenter{\hbox{\xymatrix{L^{'op} \ar[r]^-{F^{op}} \ar[rd]_-{L'(-,\tilde B)} & L \ar[d]^-{U}\\
 & {\Bbb B}}}}$ \hskip 0.5cm
with $(\tilde A,\tilde B)$  its schizophrenic object.
\endproclaim
\demo{Proof} 
\item{$\bullet $} $L(A,\tilde A), \, L'(B,\tilde B)$ are lifted to $L',L$ by the assumed condition.
\item{$\bullet $} Let $f^n\in L^n(A,A')$, then $L(f^n,\tilde A):L(A',\tilde A)\to L(A,\tilde A)$ is 
an $L'$-arrow since
$ev_{A,a^n}\circ _{n+1}L(f^n,\tilde A):=L(a^n,\tilde A)\circ _{n+1}L(f^n,\tilde A)=L(f^n\circ _{n+1}a^n,\tilde A)=:ev_{A',f^n\circ _{n+1}a^n}$, 
which is liftable $\forall a^n\in L^n(A_0,\tilde A)$. Therefore, $L(f^n,\tilde A)$ is an $L'$-arrow, and similarly, $L'(g^n,\tilde B)$ is an $L$-arrow, i.e., there exist  
maps $\xymatrix{L^{op} \ar@/^/[r]^-{G} & L' \ar@/^/[l]^-{F}}$, which are obviously functorial. 

Why do they give an adjunction?
\item{$\bullet $} (unit and counit) $1$-arrow (unit) $\eta _B:B\to GFB$ is given by 
$|\eta _B|=:V\eta _B:|B|\to |GFB|:b\mapsto [ev_{B,b}:FB\to \tilde A]$, 
$b\in |B|=L'(B_0,B)$, $|GFB|=L(FB,\tilde A)$, $|ev_{B,b}|:|FB|\to |\tilde A|$, $|FB|=L'(B,\tilde B)$, 
$|\tilde A|=L(A_0,\tilde A)\sim L'(B_0,\tilde B)$.
Why can $|\eta _B|$ be lifted to $L'$? Take the composite with evaluation maps 
$|ev_{FB,c}|\circ _1|\eta _B|(b)=|ev_{FB,c}|(ev_{B,b})=|ev_{B,b}|(c)=|c|(b)$, 
where $c\in |FB|^0=L^{'0}(B,\tilde B)=L^0(A_0,FB)$, $b\in |B|^n$. So, $|ev_{FB,c}|\circ _1|\eta _B|=|c|$ is an $L'$-arrow.
Therefore, $|\eta _B|$ is an $L'$-arrow. The counit is given symmetrically $\varepsilon _A\to FGA$, 
$|\varepsilon _A|:|A|\to |FGA|:a\mapsto [ev_{A,a}:GA\to \tilde B]$, $|A|=L(A_0,A)$, $|FGA|=L'(GA,\tilde B)$, 
$|ev_{A,a}|:|GA|\to |\tilde B|$, $|GA|=L(A,\tilde A)$, $|\tilde B|=L'(B_0,\tilde B)\sim L(A_0,\tilde A)$.
By the same argument $|\varepsilon _A|$ is an $L$-arrow.
\item{$\bullet $} (triangle identities) $G\varepsilon _A\circ _1\eta _{GA}=1_{GA}$, $F\eta _B\circ _1\varepsilon _{FB}=1_{FB}$. 
It is sufficient to prove them for underlying maps. Since forgetful functors are faithful this follows. \newline 
$|G\varepsilon _A|\circ _1|\eta _{GA}|\overset ?\to =|1_{GA}|$, where $|\eta _{GA}|:|GA|\to |GFGA|$, $|GA|=L(A,\tilde A)$, 
$|GFGA|=L(FGA,\tilde A)$, $\varepsilon _A:A\to FGA$, $|G\varepsilon _A|:|GFGA|\to |GA|$.\newline 
Take $(f^n:A\to \tilde A)\in |GA|=L^n(A,\tilde A)$, $a^m\in |A|=L^m(A_0,A)$. Two cases are possible 
$\cases \hskip -0.1cm(a)\ (f^n,n>0)\ \& \ (a^0):& \hskip -0.25cm||G\varepsilon _A|\circ _1|\eta _{GA}|(f^n)|(a^0)=|L(\varepsilon _A,\tilde A)(ev_{GA,f^n})|(a^0)=|ev_{GA,f^n}\circ _{n+1}\\
\hskip -0.1cm(b)\ (f^0)\ \& \ (a^n,n\ge 0):& \hskip -0.25cm||G\varepsilon _A|\circ _1|\eta _{GA}|(f^0)|(a^n)=|L(\varepsilon _A,\tilde A)(ev_{GA,f^0})|(a^n)\hskip 0.1cm=\hskip 0.1cm|ev_{GA,f^0}\circ _{1}\endcases $ \newline
$\cases \hskip -0.1cm(a) & \hskip -0.25cme^n\varepsilon _A|(a^0)=|ev_{GA,f^n}|\circ _{n+1}e^n|\varepsilon _A|(a^0)=|ev_{GA,f^n}|(ev_{A,e^na^0})=|ev_{A,e^na^0}|(f^n)=|f^n|(a^0)\hskip -0.1cm\\
\hskip -0.1cm(b) & \hskip -0.0cm\varepsilon _A|(a^n)\hskip 0.05cm=\hskip 0.05cm|ev_{GA,f^0}|\circ _{1}|\varepsilon _A|(a^n)\hskip 0.15cm=\hskip 0.15cm|ev_{GA,f^0}|(ev_{A,a^n})\hskip 0.15cm=\hskip 0.15cm|ev_{A,a^n}|(f^0)\hskip 0.15cm=\hskip 0.15cm|f^0|(a^n)\hskip -0.1cm\endcases $ \newline 
$\cases \hskip -0.1cm(a) & \hskip -0.25cm=:\mu ^L_{A_0,A,\tilde A}(f^n,e^na^0)=||1_{GA}|(f^n)|(a^0) \\
\hskip -0.1cm(b) & \hskip -0.25cm=:\mu ^L_{A_0,A,\tilde A}(e^nf^0,a^n)=||1_{GA}|(f^0)|(a^n) \endcases $

\vskip 0.1cm
The second triangle identity holds similarly.
\vskip 0.05cm
\item{$\bullet $} (naturality of $\eta _B, \varepsilon _A$) \ Again, it is sufficient to prove naturality for underlying maps \newline 
$\vcenter{\xymatrix{|B| \ar[r]^-{|\eta _B|} \ar[d]_-{|f|} & |GFB| \ar[d]|-{|GFf|=L(Ff,\tilde A)}\\
|B'| \ar[r]_-{|\eta _{B'}|} & |GFB'|}}$
\hskip 1cmTwo cases are $\cases (a)\ (b^n\in |B|^n, n\ge 0)\ \& \ (f^0\in L^{'0}(B,B')) & \\
(b)\ (b^0\in |B|^0)\ \& \ (f^n\in L^{'n}(B,B')) & \endcases $ \newline 
$(a)$ \hskip 1.5cm $\vcenter{\xymatrix{|B| \ar[r]^-{|\eta _B|} \ar[d]_-{|f^0|} & |GFB| \ar[d]|-{|GFf^0|=L(Ff^0,\tilde A)}\\
|B'| \ar[r]_-{|\eta _{B'}|} & |GFB'|}}$
\hskip 1.5cm $\vcenter{\xymatrix{b^n \ar@{|->}[r] \ar@{|->}[dd] & ev_{B,b^n} \ar@{|->}[d]\\
  & ev_{B,b^n}\circ _{n+1}e^n(Ff^0) \ar@{==>}[d]^-{=}\\
|f^0|(b^n) \ar@{|->}[r] & ev_{B',|f^0|(b^n)}}}$ \newline 
$(b)$ \hskip 1.5cm $\vcenter{\xymatrix{|B| \ar[r]^-{e^n|\eta _B|} \ar[d]_-{|f^n|} & |GFB| \ar[d]|-{|GFf^n|=L(Ff^n,\tilde A)}\\
|B'| \ar[r]_-{e^n|\eta _{B'}|} & |GFB'|}}$
\hskip 1.5cm $\vcenter{\xymatrix{b^0 \ar@{|->}[r] \ar@{|->}[dd] & ev_{B,e^nb^0} \ar@{|->}[d]\\
  & ev_{B,e^nb^0}\circ _{n+1}(Ff^n) \ar@{==>}[d]^-{=}\\
|f^n|(b^0) \ar@{|->}[r] & ev_{B',|f^n|(b^0)}}}$ \newline  
\vskip 0.1cm
\centerline{(recall \ $|f^n|(b^0)\equiv \mu (f^n,e^nb^0)$, \ $|f^0|(b^n)\equiv \mu (e^nf^0,b^n)$)} 
\vskip 0.1cm
Why \ \ $\cases \hskip -0.05cm(a)\ \ ev_{B,b^n}\circ _{n+1}e^n(Ff^0)=ev_{B',|f^0|(b^n)} & \\
\hskip -0.05cm(b)\ \ ev_{B,e^nb^0}\circ _{n+1}(Ff^n)=ev_{B',|f^n|(b^0)} & \endcases $? \newline 
\vskip 0.1cm
Take underlying maps:
\vskip 0.1cm
$\cases \hskip -0.1cm(a)& \hskip -0.2cm|ev_{B,b^n}|\circ _{n+1}e^n|Ff^0|(h^n)=|ev_{B,b^n}|(h^n\circ _{n+1}e^nf^0)=|h^n\circ _{n+1}e^nf^0|(b^n)=\\
\hskip -0.1cm(b)& \hskip -0.2cm|ev_{B,e^nb^0}|\circ _{n+1}|Ff^n|(h^0)=|ev_{B,e^nb^0}|(e^nh^0\circ _{n+1}f^n)=|e^nh^0\circ _{n+1}f^n|(e^nb^0)= \endcases $ 
\vskip 0.1cm
$\cases \hskip -0.1cm(a)& \hskip -0.2cm=|h^n|\circ _{n+1}|e^nf^0|(b^n)=|ev_{B',|f^0|(b^n)}|(h^n), \ \ h^n\in L^{'n}(B',\tilde B)\\
\hskip -0.1cm(b)& \hskip -0.2cm=e^n|h^0|\circ _{n+1}|f^n|(e^nb^0)=|ev_{B',|f^n|(b^0)}|(h^0), \ \ h^0\in L^{'0}(B',\tilde B)\endcases $ 
\vskip 0.1cm
(the types of the above arrows are $Ff:FB'\to FB$, $ev_{B,b}:FB\to \tilde A\ (\text{$L$-map})$, $ev_{B',|f|(b)}:FB'\to \tilde A\ (\text{$L$-map})$, 
$|ev_{B,b}|:L'(B,\tilde B)\to |\tilde B|=L'(B_0,\tilde B)$, $|ev_{B',|f|(b)}|:L'(B',\tilde B)\to |\tilde B|=L'(B_0,\tilde B)$, 
$|Ff|:L'(B',\tilde B)\to L'(B,\tilde B)$, $|Ff|=L'(f,\tilde B)$). 
\vskip 0.1cm
Therefore, $\eta _B$ is natural. Similarly, $\varepsilon _A$ is natural.
\item{$\bullet $} (isomorphisms-functors $\xymatrix{{L(A,FB)\ \ } \ar@/^0.9pc/[r]^-{\theta _{A,B}} & {\ \ L'(B,GA)} \ar@/^0.9pc/[l]^-{\theta ^*_{A,B}}}$)\newline 
Define \ $\cases \theta _{A,B}(f^n):=G(f^n)\circ _{n+1}e^n(\eta _B), & f^n\in L^n(A,FB)\\
\theta ^*_{A,B}(g^n):=F(g^n)\circ _{n+1}e^n(\varepsilon _A), & g^n\in L^{'n}(B,GA) \endcases $\newline 
Let $g^n\in L^{'n}(B,GA)$. Then $\theta _{A,B}(\theta ^*_{A,B}(g^n))\hskip -0.05cm:=G(Fg^n\circ _{n+1}e^n(\varepsilon _A))\circ _{n+1}e^n(\eta _B)=e^n(G\varepsilon _A)\circ _{n+1}GFg^n\circ _{n+1}e^n(\eta _B)\underset {\text{nat. of $\eta _B$}} \to {=} e^n(G\varepsilon _A)\circ _{n+1}e^n(\eta _{GA})\circ _{n+1}g^n\underset {\text{triangle id.}} \to {=}e^n(1_{GA})\circ _{n+1}g^n=e^{n+1}(GA)\circ _{n+1}g^n=g^n$.
Similarly, \ $\theta ^*_{A,B}(\theta _{A,B}(f^n))=f^n, \ f^n\in L^n(A,FB)$. 
$\theta _{A,B},\, \theta ^*_{A,B}$ are obviously functors. So, they are isomorphisms.
\item{$\bullet $} (naturality of $\theta _{A,B},\, \theta ^*_{A,B}$) \ We need to prove the diagram \newline 
${\vcenter{\xymatrix{ A & B & L(A,FB) \ar[r]^-{e^n\theta _{A,B}} \ar[d]_-{L(x^n,Fy^n)} & L'(B,GA) \ar[d]^-{L'(y^n,Gx^n)}\\
A' \ar[u]^-{x^n} & B' \ar[u]^-{y^n} & L(A',FB') \ar[r]^-{e^n\theta _{A',B'}} & L'(B',GA') }}}$ \hskip 0.2cm commutes.
\vskip 0.15cm
$L'(y^n,Gx^n)\circ _{n+1}e^n\theta _{A,B}\overset {?}\to {=}e^n\theta _{A',B'}\circ _{n+1}L(x^n,Fy^n)$
\vskip 0.25cm
Two cases are: $\cases \hskip -0.1cm(a) & \hskip -0.2cm(f^0\in L(A,FB))\ \& \ (x^n,y^n,n>0)\\ 
\hskip -0.1cm(b) & \hskip -0.2cm(f^n\in L(A,FB), n\ge 0)\ \& \ (x^0,y^0)\endcases $
\vskip 0.1cm
$(a)$\hskip 0.7cm $\vcenter{\xymatrix{f^0 \ar@{|->}[r] \ar@{|->}[dd] & e^nG(f^0)\circ _{n+1}e^n(\eta _B) \ar@{|->}[d]\\
 & Gx^n\circ _{n+1}(e^nG(f^0)\circ _{n+1}e^n(\eta _B))\circ _{n+1}y^n \ar@{==>}[d]^-{\ ?}_-{=\ } \ar@/^1.5pc/@{==>}@<3cm>[dd]_-{=}^-{({\text{$\eta _B$ is nat.}})} \\
Fy^n\circ _{n+1}e^nf^0\circ _{n+1}x^n \ar@{|->}[r] & G(Fy^n\circ _{n+1}e^nf^0\circ _{n+1}x^n)\circ _{n+1}e^n(\eta _{B'}) \ar@{==>}[d]_-{=\ }\\
 & Gx^n\circ _{n+1}e^nGf^0\circ _{n+1}GFy^n\circ _{n+1}e^n(\eta _{B'})}}$
\vskip 0.1cm
$(b)$\hskip 0.25cm $\vcenter{\xymatrix{f^n \ar@{|->}[r] \ar@{|->}[dd] & G(f^n)\circ _{n+1}e^n(\eta _B) \ar@{|->}[d]\\
 & e^nGx^0\circ _{n+1}(G(f^n)\circ _{n+1}e^n(\eta _B))\circ _{n+1}e^ny^0 \ar@{==>}[d]^-{\ ?}_-{=\ } \ar@/^1.5pc/@{==>}@<3cm>[dd]_-{=}^-{({\text{$\eta _B$ is nat.}})} \\
e^nFy^0\circ _{n+1}f^n\circ _{n+1}e^nx^0 \ar@{|->}[r] & G(e^nFy^0\circ _{n+1}f^n\circ _{n+1}e^nx^0)\circ _{n+1}e^n(\eta _{B'}) \ar@{==>}[d]_-{=\ }\\
 & e^nGx^0\circ _{n+1}Gf^n\circ _{n+1}e^nGFy^0\circ _{n+1}e^n(\eta _{B'})}}$
\vskip 0.15cm
Therefore, $L$ and $L'$ are concretely dually adjoint. This correspondence is natural (by condition) and strict ($\theta _{A,B}$ 
and $\theta ^*_{A,B}$ are isomorphisms).
\hfill $\square $
\enddemo

\vskip 0.2cm
{\bf Corollary.} Concrete natural duality is a {\bf strict} adjunction.                \hfill $\square $

\vskip 0.3cm
\centerline{{\bf Well-known dualities} \cite{P-Th, Bel, A-H-S}}
\vskip 0.25cm

All dualities below are of first order, natural \cite{P-Th}, and obtained by restriction of appropriate dual adjunctions.

\item{1.} $\bold{Vec}_k$ is dually equivalent to itself $\xymatrix{\bold{Vec}_k^{op} \ar@/^/[rr]^-{\bold{Vec}_k(-,k)} && \bold{Vec}_k
\ar@/^/[ll]^-{\bold{Vec}_k(-,k)}_-{\perp }}$, where $\bold{Vec}_k$ is a category of vector spaces over field $k$
\item{2.} $\bold{Set}^{op}\sim \text{\bf Complete Atomic Boolean Algebras}$
\item{3.} $\bold{Bool}^{op}\sim \text{\bf Boolean Spaces}$ (Stone duality), where $\bold{Bool}$ is a category of Boolean rings (every element is idempotent). It is obtained from the dual adjunction $\xymatrix{\bold{CRing} \ar@/^/[rr]^-{\bold{CRing}(-,\bold{2})} && \bold{Top} \ar@/^/[ll]^-{\bold{Top}(-,\bold{2})}_-{\perp }}$, where $\bold 2$ is two-element ring and discrete topological space. $\xymatrix{\bold{CRing}(A,\bold 2) \ar@{^{(}->}[r] & \bold 2^A}$ (subspace in Tychonoff topology)
\item{4.} $\text{\rm hom}(-,\Bbb R/\Bbb Z):\bold{CompAb}^{op}\sim \bold{Ab}$ (Pontryagin duality), where $\bold{CompAb}$, $\bold{Ab}$ are categories of compact abelian groups and abelian groups respectively
\item{5.} $\text{hom}(-,\Bbb C):\text{\bf C*Alg}^{op}\sim \bold{CHTop}$ (Gelfand-Naimark duality), where $\text{\bf C*Alg}$, $\bold{CHTop}$ are categories of commutative $\Bbb C^*$-algebras and compact Hausdorff spaces. $\xymatrix{\text{\bf C*Alg}(A,\Bbb C) \ar@{^{(}->}[r] & {\Bbb C}^A}$ (subspace in Tychonoff topology)

\head {\bf 7. Vinogradov duality} \endhead

Let $K$ be a commutative ring, $A$ a commutative algebra over $K$, $A\text{-}\bold{Mod}\hookrightarrow K\text{-}\bold{Mod}$ be the categories
of modules over $A$ and $K$ respectively. 

\vskip 0.1cm
{\bf Definition 7.1.} \cite{V-K-L} For $P,Q \in Ob\, (A\text{-}\bold{Mod})$
\item{$\bullet $} $K$-linear maps \newline
$l(a):=a\, \cdot \, -, r(a):=-\, \cdot \, a, \delta (a):=l(a)-r(a):K\text{-}\bold{Mod}(P,Q)\to K\text{-}\bold{Mod}(P,Q)$ 
are called {\bf left, right multiplications} and {\bf difference operator} (by element $a\in A$),
\item{$\bullet $} A $K$-linear map $\Delta :P\to Q$ is a {\bf differential operator of order $\le r$} \ if \ $\forall a_0, a_1,\dots , a_r \in A$ 
$\delta _{a_0,a_1,\dots ,a_r}(\Delta )=0$, where $\delta _{a_0,a_1,\dots ,a_r}:=\delta _{a_0}\circ \delta _{a_1}\circ \cdots \circ \delta _{a_r}$. \hfill $\square $ 

\vskip 0.1cm
\proclaim{\bf Lemma 7.1} 
\item{$\bullet $} If $\Delta _1\in K\text{-}\bold{Mod}(P,Q), \Delta _2\in K\text{-}\bold{Mod}(Q,R)$ are differential operators of order $\le r$ and $\le s$ 
respectively, then $\Delta _2\circ \Delta _1:K\text{-}\bold{Mod}(P,R)$ is a differential operator of order $\le r+s$,
\item{$\bullet $} $\forall a\in A, \ P\in Ob\, (A\text{-}\bold{Mod})$ module multiplication (by $a$) \, $l_a:P\to P:p\mapsto ap$ is a differential operator of order $0$.  \hfill $\square $ 
\endproclaim 

The differential operators between $A$-modules form the arrows of a category $A\text{-}\bold{Diff}$, such that 
$A\text{-}\bold{Mod}\hookrightarrow A\text{-}\bold{Diff}\hookrightarrow K\text{-}\bold{Mod}$,  
 and the first two categories have the same objects. $A\text{-}\bold{Diff}$ is enriched in $(K\text{-}\bold{Mod}, \otimes _K)$  and 
enriched in two different ways in $(A\text{-}\bold{Mod}, \otimes _K)$, except that composition is not  an $A$-module map. Module 
multiplication for the first enrichment $A\text{-}\bold{Diff}$ in $(A\text{-}\bold{Mod}, \otimes _K)$ is given by 
$A\times A\text{-}\bold{Diff}(P,Q)\to A\text{-}\bold{Diff}(P,Q):(a,\Delta )\mapsto l_a\circ \Delta $, for the second enrichment by
$A\times A\text{-}\bold{Diff}(P,Q)\to A\text{-}\bold{Diff}(P,Q):(a,\Delta )\mapsto \Delta \circ l_a$. Denote $A\text{-}\bold{Diff}$ with 
left module multiplication in hom-sets $l_a\circ -$ by the same name $A\text{-}\bold{Diff}$ and with right multiplication in hom-sets 
$-\circ l_a$ by $A\text{-}\bold{Diff^+}$.

\vskip 0.1cm
\proclaim{\bf Proposition 7.1} 
\item{$\bullet $} $\forall P,Q\in Ob\, (A\text{-}\bold{Mod})$ \, $A\text{-}\bold{Diff}(P,Q)=\bigcup \limits _{s=0}^{\infty }\bold{Diff}_s(P,Q)$,
\, $A\text{-}\bold{Diff^+}(P,Q)=\bigcup \limits _{s=0}^{\infty }\bold{Diff}_s^+(P,Q)$ are filtered $A$-modules by submodules of differential 
operators of order $\le s, \, s=0,1,...$,
\item{$\bullet $} $\forall P\in Ob\, (A\text{-}\bold{Mod})$ \, $A\text{-}\bold{Diff}(P,P)$ is an associative $K$-algebra.  \hfill $\square $
\endproclaim

\proclaim{\bf Proposition 7.2} 
\item{$\bullet $} $\bold{Diff}_s(P,-),\, \bold{Diff}_s^+(-,P):A\text{-}\bold{Mod}\to A\text{-}\bold{Mod}$ are $A$-linear functors,
\item{$\bullet $} $\forall P\in Ob\, (A\text{-}\bold{Mod})$ functor $\bold{Diff}_s^+(-,P)$ is representable by object 
$\bold{Diff}_s^+(P):=\bold{Diff}_s^+(A,P)$, i.e. $\forall Q\in Ob\, (A\text{-}\bold{Mod})$ $A\text{-}\bold{Mod}(Q,\bold{Diff}_s^+(P))@>\sim >>\bold{Diff}_s^+(Q,P)$, 
\item{$\bullet $} $\forall P\in Ob\, (A\text{-}\bold{Mod})$ functor $\bold{Diff}_s(P,-)$ is representable by object 
$\bold{Jet}^s(P):=A\otimes _KP\, \text{\rm mod} \, \mu ^{s+1}$, where $\mu ^{s+1}$ is a submodule of $A\otimes _KP$ generated by elements 
$\delta ^{a_0}\circ \cdots \circ \delta ^{a_{s+1}}(a\otimes p)$ {\rm [$\delta ^{b}(a\otimes p):=ab\otimes p-a\otimes bp$]}, i.e. $\forall Q\in Ob\, (A\text{-}\bold{Mod})$ $A\text{-}\bold{Mod}(\bold{Jet}^s(P),Q)@>\sim >>\bold{Diff}_s(P,Q)$,
\item{$\bullet $} inclusion $A\text{-}\bold{Mod}\hookrightarrow A\text{-}\bold{Diff}^+$ is an (enriched) left adjoint with counit
$ev:\bold{Diff}^+(P)\to P:\Delta \mapsto \Delta (1)$, i.e. $\forall \Delta \in \bold{Diff}^+(Q,P)$ 
$\exists ! f_{\Delta }\in A\text{-}\bold{Mod}(Q,\bold{Diff}^+(P))$ such that 
$$\xymatrix{\bold{Diff}^+(P) \ar[r]^-{ev} & P\\
Q \ar@{-->}[u]^-{f_{\Delta }} \ar[ru]_-{\Delta } & }$$
and this correspondence is $A$-linear, \ $f_{\Delta }:q\mapsto (a\mapsto \Delta (aq))$,
\item{$\bullet $} inclusion $A\text{-}\bold{Mod}\hookrightarrow A\text{-}\bold{Diff}$ is an (enriched) right adjoint with unit 
$j^{\infty }:P\to \bold{Jet}^{\infty }(P):p\mapsto 1\otimes p\, \text{\rm mod}\, \mu ^{\infty }$ {\rm [$\mu ^{\infty }:=\bigcap \limits _{s=0}^{\infty }\mu ^s$]}, i.e. 
$\forall \Delta \in \bold{Diff}(P,Q)$ 
$\exists ! f^{\Delta }\in A\text{-}\bold{Mod}(\bold{Jet}^{\infty }(P),Q)$ such that 
$$\xymatrix{P \ar[r]^-{j^{\infty }} \ar[dr]_-{\Delta } & \bold{Jet}^{\infty }(P) \ar@{-->}[d]^-{f^{\Delta }}\\
 & Q  }$$
and this correspondence is $A$-linear, \ $f^{\Delta }:(a\otimes p)\, \text{\rm mod}\, \mu^{\infty } \mapsto a\Delta (p)$,
\item{$\bullet $} subcategory $A\text{-}\bold{Mod}$ is reflective and coreflective in $A\text{-}\bold{Diff}$ (enriched in $K\text{-}\bold{Mod}$). \hfill $\square $
\endproclaim

\vskip 0.15cm
$\forall s\in \Bbb N$ introduce two full subcategories of $A\text{-}\bold{Mod}$: 
\item{$\bullet $} $A\text{-}\bold{Mod}\text{-}\bold{Diff}_s$, consisting of all $A$-modules of type $\bold{Diff}_s(P,A), \ P\in Ob\, (A\text{-}\bold{Mod})$,  
\item{$\bullet $} $A\text{-}\bold{Mod}\text{-}\bold{Jet}^s$, consisting of all $A$-modules of type $\bold{Jet}^s(P), \ P\in Ob\, (A\text{-}\bold{Mod})$.

\vskip 0.3cm
\proclaim{\bf Proposition 7.3 (Vinogradov Duality)} 
For a commutative 
algebra $A$ there is a concrete natural dual adjunction  
$\xymatrix{A\text{-}\bold{Mod}\text{-}\bold{Diff}_s^{op} \ar@/^1.5ex/[rr]_-{\perp } & & A\text{-}\bold{Mod}\text{-}\bold{Jet}^s \ar@/^1.5ex/[ll] }$, $s\in \Bbb N$, 
obtained by restriction of $\xymatrix{A\text{-}\bold{Mod}^{op} \ar@/^1.5ex/[rr]^-{A\text{-}\bold{Mod}(-,A)}_-{\perp } & & A\text{-}\bold{Mod} \ar@/^1.5ex/[ll]^-{A\text{-}\bold{Mod}(-,A)} }$. $A$ is a schizophrenic object. \hfill $\square $ 
\endproclaim 

\vskip 0.1cm
{\bf Remarks.}
\vskip 0.15cm
\item{$\bullet $} The above duality theorem is not stated explicitly in \cite{V-K-L} but the result is implicitly there.
\item{$\bullet $} The above proposition states a formal analogue of duality between differential operators and jets over a fixed manifold $X$. 
Geometric modules of sections of vector bundles over $X$ correspond to modules $P$ over $C^{\infty }(X)$ with the property 
$\bigcap \limits _{x\in X}\mu _xP=0$, where $\mu _x$ is a maximal ideal at point $x\in X$. Functors $\bold{Diff}_s(-,A)$ and $\bold{Jet}^s(-)$ preserve the module property to be geometric \cite{V-K-L}.
\item{$\bullet $} This duality is an alternative (algebraic) way to introduce jet-bundles in Geometry (instead of 
the classical approach
due to Grothendieck and Ehresmann as equivalence classes of maps which tangent of order $s$ at a point). When $A=C^{\infty }(X)$ and $P$ is a geometric 
module realizable as a vector bundle $V(P)$ over $X$, then $\bold{Jet}^s(P)$ is realizable as $\bold{Jet}^s(V(P))$ over $X$ in the 
classical sense \cite{V-K-L, Vin1, Vin2}.   \hfill $\square $

\head {\bf 8. Duality for differential equations}\endhead
\vskip 0.0cm
\proclaim{\bf {Proposition 8.1}} Let \, {\bf UAlg} be a category of universal algebras with a representable forgetful functor.
Then every topological algebra ${\goth A}$ is a schizophrenic object (see \cite{P-Th}), and so  yields
a natural dual adjunction between {\bf UAlg} and {\bf Top}.
\endproclaim
\demo{Proof}
\item{$\bullet $} The initial topology on {\bf UAlg}(B,$\goth A$) gives the initial lifting with respect to evaluation maps
$ev_{B,\, b}:\text{\bf UAlg}(B,\goth A)\to |\goth A|$, \ $b\in |B|$.
\item{$\bullet $} The algebra of continuous functions $\text{\bf Top}(X,\goth A)$ is initial with respect to 
the evaluation maps
$ev_{X,\, x}:\text{\bf Top}(X,\goth A)\to |\goth A|$, $x\in |X|$ (which are obviously homomorphisms) since operations in \linebreak
$\text{\bf Top}(X,\goth A)$ are pointwise and each arrow $f\in \text{\bf Top}(X,\goth A)$ is completely determined by all its values
$ev_{X,\, x}(f)=|f|(x)$, \ $x\in |X|$. Hence,  if $g:|B|\to \text{\bf Top}(X,\goth A)$ is a {\bf Set}-map such that $\forall x\in |X|$
$ev_{X,\, x}\circ g$ is a homomorphism,  ($\omega _n(ev_{X,\, x}\circ g)b_1,...,(ev_{X,\, x}\circ g)b_n=
ev_{X,\, x}\circ g\omega _nb_1,...,b_n=ev_{X,\, x}\omega _ngb_1,...,gb_n$, where $\omega _n$ is an $n$-ary operation.
The first equality holds because $ev_{X,\, x}\circ g$ is a homomorphism, the second equality because $ev_{X,\, x}$ is a homomorphism), then $g$
is a homomorphism since two maps whose values coincide at each point coincide themselves.
\hfill $\square $
\enddemo

\proclaim{\bf Corollary} Take {\bf UAlg}$=k\text{-}\Lambda \text{-}\bold{Alg}$, the category of exterior differential algebras over a field $k$ ($\Bbb R$ or $\Bbb C$).  These are thought of as presenting ``generalized differential
equations". Take $\goth A=\Lambda (C^{\infty }({\Bbb R}^n))$ or $\Lambda (C^{\omega }({\Bbb C}^n))$
(which acts as a parameter space) with a topology not weaker than $jet^{\infty }$. Then there exists a natural dual adjunction
$\xymatrix{k\text{-}\Lambda \text{-}\bold{Alg}^{op} \ar@/_/[r]^-{\perp } & {\text{\bf Top}} \ar@/_/[l]}$ (between differential equations and their solution spaces). \hfill $\square $
\endproclaim
{\bf Remark.} If we regard the  category $k\text{-}\Lambda \text{-}\bold{Alg}$ whose forgetful functor is representable, we will get a lot of extra  ``points"
which do not have geometric sense. Only graded maps of degree $0$ to $\goth A$
have geometric sense (they present integral manifolds of dimension not bigger than $n$). In this case,
the  representation of exterior
differential algebras, when it exists, will not be via their solution spaces but via much bigger spaces. If we restrict 
$k\text{-}\Lambda \text{-}\bold{Alg}$ to only
graded morphisms of degree $0$ then the forgetful functor is not representable. But the notion of 
``schizophrenic object" still makes sense
and the theorem for natural dual adjunction \cite{P-Th} still holds. So, there is a representation of exterior differential algebras via
their usual solution spaces. \hfill $\square $
\vskip 0.2cm
We denote concrete subcategories of $\bold{Top}$ dual to categories $k\text{-}\bold{Alg}$ (algebras over $k$) and 
$k\text{-}\Lambda \text{-}\bold{Alg}$ (exterior differential algebras over $k$ with graded degree $0$ morphisms) 
by ${\text{\tt alg-}}\bold{Sol}$ and ${\text{\tt diff-}}\bold{Sol}$ respectively, i.e., 
$k\text{-}\bold{Alg}^{op}\sim {\text{\tt alg-}}\bold{Sol}$, $k\text{-}\Lambda \text{-}\bold{Alg}^{op}\sim {\text{\tt diff-}}\bold{Sol}$. 
In particular, ${\text{\tt alg-}}\bold{Sol}$ contains all algebraic and all smooth $k$-manifolds 
($k=\Bbb R$ or $\Bbb C$), ${\text{\tt diff-}}\bold{Sol}$ contains all spaces of the form ${\text{\tt alg-}}\bold{Sol}(k^n,X)$ 
(with representing object $\goth A=\Lambda (C^{\infty }(k^n))$).

\vskip 0.3cm
\proclaim{\bf {Lemma 8.1} ({\rm rough structure of} $\text{\tt diff-}\bold{Sol}$)}
\item{$\bullet $} $Ob\, (\text{\tt diff-}\bold{Sol})$ are pairs $(X,\coprod\limits _{i=1}^n\Cal F_i)$ where $X\hskip -0.1cm:=k\text{-}\Lambda \text{-}\bold{Alg}(D,k)=k\text{-}\bold{Alg}(D,k)\in Ob\, (\text{\tt alg-}\bold{Sol})$, $\Cal F_i\subset \text{\tt alg-}\bold{Sol}(k^i,X),\ 1\le i\le n$ [$\Cal F_i$ are not arbitrary subspaces of $\text{\tt alg-}\bold{Sol}(k^i,X)$].
\item{$\bullet $} $Ar\, (\text{\tt diff-}\bold{Sol})$ are pairs $(f,\coprod\limits _{i=1}^n\text{\tt alg-}\bold{Sol}(k^i,f)):(X,\coprod\limits _{i=1}^n\Cal F_i)\to (X',\coprod\limits _{i=1}^n\Cal F'_i)$ where $f:X\to X'\in Ar\, (\text{\tt alg-}\bold{Sol})$, $\text{\tt alg-}\bold{Sol}(k^i,f):\Cal F_i\to \Cal F'_i, \ 1\le i\le n$. \hfill $\square $
\endproclaim

\vskip 0.0cm
\proclaim{\bf {Proposition 8.2}} There are the following adjunctions
\item{$\bullet $} $\xymatrix{k\text{-}\bold{Alg} \ar@/^/[r]^-{\ \Lambda _k}_-{\bot } & k\text{-}\Lambda \text{-}\bold{Alg} \ar@/^/[l]^-{p_0}}$ where
$\Lambda _k$ is the {\bf free exterior differential algebra functor},  $p_0$ is the projection onto the subalgebra of degree-$0$ elements,
\vskip 0.1cm
\item{$\bullet $} $\xymatrix{\text{\tt alg-}\bold{Sol} \ar@/^/[r]^-{\text{ \rm hom}(k^n,-)}_-{\top } & \text{\tt diff-}\bold{Sol} \ar@/^/[l]^-{b}}$
where \, $b$ \, is   the base space functor,

such that
$$\xymatrix{ k\text{-}\Lambda \text{-}\bold{Alg}^{op} \ar@/^2ex/[rr]^-{F}_-{\sim } \ar@/^2ex/[dd]^-{p_0^{op}} && \text{\tt diff-}\bold{Sol} \ar@/^2ex/[ll]^-{F'}_-{\sim } \ar@/^2ex/[dd]^-{b}\\
 && \\
k\text{-}\bold{Alg}^{op} \ar@/^2ex/[rr]^-{G}_-{\sim } \ar@/^2ex/[uu]^-{\Lambda _k^{op}}_-{\ \vdash } && \text{\tt alg-}\bold{Sol} \ar@/^2ex/[ll]^-{G'}_-{\sim } \ar@/^2ex/[uu]^-{\text{hom}(k^n,-)}_-{\ \vdash }
}$$
\endproclaim
\demo{Proof}
\item{$\bullet $} $\xymatrix{k\text{-}\Lambda \text{-}\bold{Alg}(\Lambda _k(A),D) \ar[r]^-{\sim }  & k\text{-}\bold{Alg}(A,p_0(D))\\
\rho  \ar@{{|}->}[r]^-{\sim } \ar[u]^-{\in } & \rho _0 \ar[u]_-{\in }
}$ \hskip 1cm(natural in $A$ and $D$)\newline
where  $\rho _0$ is the $0$-component of graded degree $0$ homomorphism $\rho =\bigoplus\limits_{i\ge 0}\rho _i$.
\item{$\bullet $} $\xymatrix{\text{\tt diff-}\bold{Sol}(S,\text{\rm hom}(k^n,X)) \ar[r]^-{\sim } & \text{\tt alg-}\bold{Sol}(b(S),X) \\
f \ar@{{|}->}[r]^-{\sim } \ar[u]^-{\in } & f \ar[u]_-{\in }
}$ \hskip 1cm(natural in $S$ and $X$)\newline
where:  $S$ is a pair $(b(S),\coprod\limits _{i=1}^n\Cal F_i), \ \Cal F_i\subset \text{hom}(k^i,b(S)), \ 1\le i\le n$, \
right $f:b(S)\to X$ is a usual map, and left $f:=(f,\coprod\limits _{i=1}^n\text{hom}(k^i,f)):(b(S),\coprod\limits _{i=1}^n\Cal F_i)\to (X,\coprod\limits _{i=1}^n\text{hom}(k^i,X))$.
\vskip 0.2cm
The above square of adjunctions is immediate.
\hfill $\square $
\enddemo

\vskip 0.2cm
\subhead 8.1. Cartan involution \endsubhead

\vskip 0.1cm
For systems in Cartan involution (as defined below) a (single) solution can be calculated recursively beginning from smallest $0$ dimension. By Cartan's theorem \cite{BC3G, Car1, Fin, Vas}
every system can be made into such a form by a sufficient number of differential prolongations \cite{BC3G, Car1, Fin, Vas}. There is a 
cohomological criterion for systems to be in involution. 

\vskip 0.1cm
{\bf Definition 8.1.1.} Let $\Cal A\in Ob\, (k\text{-}\Lambda \text{-}\bold{Alg})$, $\goth A_n$ be $\Lambda _{\Bbb R}(C^{\infty }(\Bbb R^n))$  or  $\Lambda _{\Bbb C}(C^{\omega }(\Bbb C^n))$, \, $n\ge 0$. 
\item{$\bullet $} Any (differential homomorphism of degree $0$) $\rho :\Cal A\to \goth A_n$ is called an {\bf integral manifold} of $\Cal A$ (of dimension not bigger than $n$).
\item{$\bullet $} $deg\, (\rho :\Cal A\to \goth A_n)=m, \ 0\le m\le n$, \ iff $\rho $ can be factored through a $\gamma :\Cal A\to \goth A_m$, i.e.,
$\vcenter{\xymatrix{\Cal A \ar@{..>}[r]^-{\gamma } \ar[dr]_-{\rho } & \goth A_m \ar@{..>}[d] \\
 & \goth A_n}}$ \, and $m$ is the smallest such number.
\item{$\bullet $} $deg\, (\Cal A)=n$ \, iff maximal degree of integral manifolds of \, $\Cal A$ \, is $n$.
\item{$\bullet $} $\Cal A$, $deg(\Cal A)=n$, is in {\bf Cartan involution} iff for each $m$-dimensional integral manifold $\rho :\Cal A\to \goth A_m$, $m<n$, 
there exists an $(m+1)$-dimensional integral manifold $\beta :\Cal A\to \goth A_{m+1}$ which {\bf contains} $\rho $, i.e., \hskip 0.1cm
$\vcenter
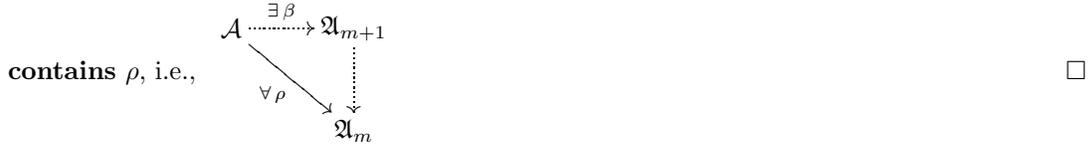
{\xymatrix{\Cal A \ar@{..>}[r]^-{\exists \, \beta } \ar[dr]_-{\forall \, \rho } & \goth A_{m+1} \ar@{..>}[d] \\
 & \goth A_{m}}}$  \hfill $\square $ 

\vskip 0.25cm
{\bf Remarks.}
\item{$\bullet $} $\goth A_0$ is just $k$ ($\Bbb R$ or $\Bbb C$) with the trivial differential. $\rho :\Cal A\to \goth A_0$ corresponds to a point $b\, (\rho ):b\, (\Cal A)\to k$.
Each point of $\Cal A$ is a $0$-dimensional integral manifold.
\item{$\bullet $}  The original Cartan definition was for classical algebras (quotient algebras of $\Lambda _{\Bbb R}(C^{\omega }(\Bbb R^N))$) 
and in terms of 'infinitesimal integral elements' (nondifferential homomorphisms of degree $0$ \, $f:\Cal A\to \Lambda _{k}(d\tau ^1,...,d\tau ^N)$) \cite{BC3G, Car1, Fin}. For that case, the two definitions coincide.
\item{$\bullet $} By a number of differential prolongations (adding new jet-variables with obvious relations), every classical 
system can be put into Cartan involution form (E. Cartan's theorem). 
\item{$\bullet $} The integration step (constructing an integral manifold  of $1$ higher dimension) is done by the method of ``Cauchy characteristics".  \hfill $\square $

\vskip 0.3cm
\proclaim{\bf Proposition 8.1.1} Let $\Cal A$ be a quotient algebra of $\Lambda _{\Bbb R}(C^{\omega }(\Bbb R^N))$, $deg(\Cal A)=n$, 
corresponding to a system of differential equations $\vcenter
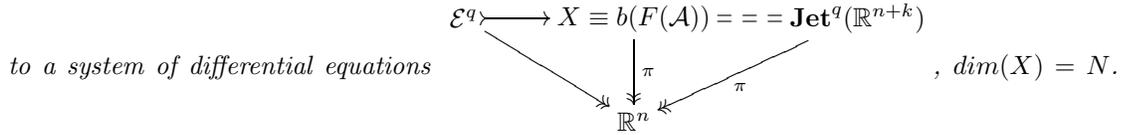
{\xymatrix{\Cal E^q \ar@{>->}[r] \ar@{->>}[dr] & X\equiv b(F(\Cal A)) \ar@{->>}[d]^-{\pi } \ar@{==}[r] & \bold{Jet}^q(\Bbb R^{n+k}) \ar@{->>}[dl]^-{\pi } \\
 & \Bbb R^n & }}$, $dim(X)=N$.  
Then $\Cal A$ is in Cartan involution iff 
the following {\bf Spencer $\delta $-complex} is acyclic:
\vskip 0.3cm
$\xymatrix{0 \ar[r] & g^{(r)} \ar[r]^-{\delta } & g^{(r-1)}\otimes \Lambda ^1(\Bbb R^n) \ar[r]^-{\delta } & g^{(r-2)}\otimes \Lambda ^2(\Bbb R^n) \ar[r]^-{\delta } & \cdots  \\
\cdots \ar[r]^-{\delta } & g^{(r-n)}\otimes \Lambda ^n(\Bbb R^n) \ar[r] & 0 & & }$
\vskip 0.2cm
where $g^{(r)}:=T(\bold{Jet}^r(\Cal E^q))\bigcap V\pi ^{q+r}_{q+r-1}\hookrightarrow S_{q+r}(T_*\Bbb R^n)\otimes V\pi $ is $r$-th prolongation of symbol $g$,  
$\pi ^{q+r}_{q+r-1}:\bold{Jet}^{q+r}(\Bbb R^{n+k})\to \bold{Jet}^{q+r-1}(\Bbb R^{n+k})$ is a natural projection of jet-bundles,
$V$ is the ``vertical" subbundle, $S_p$ is the $p$-th symmetric power,\newline
$\delta (\alpha _1\cdots \alpha _{q+r-l}\otimes v\otimes \beta _1\wedge \cdots \wedge \beta _l):=\sum \limits _{i=1}^{q+r-l}\alpha _1\cdots \hat {\alpha _i}\cdots \alpha _{q+r-l}\otimes v\otimes \alpha _i\wedge \beta _1\wedge \cdots \wedge \beta _l$,
\ $v$ is a section of $V\pi $.
\endproclaim 
\demo{Proof} See \cite{A-V-L, Sei, Vin2, V-K-L, Ver}   \hfill $\square $
\enddemo 

\vskip 0.1cm
The original Cartan  involutivity test was in terms of certain dimensions of  ``infinitesimal integral elements". The above theorem is 
due to J.P. Serre \cite{A-V-L, La-Se}.

\head {\bf 9. Gelfand-Naimark 2-duality}\endhead

Let $\xymatrix{\bold {C^{*}Alg^{op}} \ar@/^/[r]^{F}& \bold {CHTop}\ar@/^/[l]^{G}_{\perp}}$ be the usual Gelfand-Naimark duality between commutative $C^*$-algebras and compact Hausdorff spaces. Both categories are strict 2-categories with homotopy classes of homotopies as 2-cells (homotopy of $C^*$-algebras is a homotopy in $\bold {Top}$ each instance of which is a $C^*$-algebra homomorphism). The reasonable question is: can it be extended to a 2-duality? The answer is yes.

By definition 
$$\xymatrix{
{\bold {C^*Alg}}(A,B)\times |A| \ar[r]^-{ev} & |B| & {\bold {C^*Alg}}(B,\Bbb C)\times {\bold {C^*Alg}}(A,B)\ar[r]^-{c_{A,B,\Bbb C}} & {\bold {C^*Alg}}(A,\Bbb C)\\
|I|\times |A|\ar[u]^{f\times 1} \ar[ur]_{\bar f} && {\bold {C^*Alg}}(B,\Bbb C)\times |I| \ar[u]^{1\times f} \ar[ur]_{F(\bar f)}\\
}$$
So that, if $f:|I|\times |A|\to |B|$ is a homotopy in $\bold {C^*Alg}$, then its image in $\bold {CHTop}$ is $F(\bar f):|F(B)|\times |I|\to |F(A)|$ (where $|\ |$ denotes underlying set or map).

We need to prove that such extended $F$ preserves 2-categorical structure (for $G$ proof is symmetric).

\centerline{\bf Preserving homotopies}

\proclaim{\bf Lemma 9.1}If $B$ is locally compact then $\xymatrix{{\bold {Top}}(B,C)\times {\bold {Top}}(A,B)\ar[r]^-{c_{A,B,C}} & {\bold {Top}}(A,C)}$ is continuous {\rm (with} compact-open {\rm topology in all hom-sets)}.
\endproclaim
\demo\nofrills{Proof \ }is standard. Let $f=g\circ h=c_{A,B,C}(g,h)$. Take $U^K$ be a (subbase) nbhd of $f$. Sufficient to show that $\exists$ (subbase) nbhds $U_1^{K^1}\ni g,\  U_2^{K^2}\ni h$, s.t. $U_1^{K^1}\circ U_2^{K^2}=c_{A,B,C}(U_1^{K^1},U_2^{K^2})\subset  U^{K}$. Take $U_1=U,\ K_2=K,\ K_1$ be a compact nbhd of $h(K)$, s.t. $K_1\subset g^{-1}(U)$ ($K_1$ exists by local compactness of $B$), $U_2=\text{int}(K_1)$. \hfill $\square $
\enddemo
\proclaim{\bf Corollary} If $A$ is locally compact then $ev_{A,B}:\bold {Top}(A,B)\times |A|\to |B|$ is continuous.
\endproclaim 
\demo{Proof} Each space $A$ is homeomorphic to $\bold {Top}(1,A)$ (with compact-open topology), and $ev_{A,B}$ corresponds to $c_{1,A,B}$. \hfill $\square $
\enddemo
\proclaim{\bf Lemma 9.2}  $\bullet $ Initial topology on $|F(A)|=\bold {C^*Alg}(A,\Bbb C)$ w.r.t. evaluation maps $\forall a\in A$  $\xymatrix{{\bold {C^*Alg}}(A,\Bbb C)\times 1
\ar[r]_{1\times a} & {\bold {C^*Alg}}(A,\Bbb C)\times |A|\ar[r]_-{ev} & |\Bbb C|}$ is point-open.

\hskip 1.9cm$\bullet$ Initial topology on $|G(X)|=\bold {CHTop}(X,\Bbb C)$ w.r.t. evaluation maps $\forall x\in X$  $\xymatrix{{\bold {CHTop}}(X,\Bbb C)\times 1
\ar[r]_{1\times x} & {\bold {CHTop}}(X,\Bbb C)\times |X|\ar[r]_-{ev} & |\Bbb C|}$ is compact-open.
\endproclaim
\demo{Proof} See \cite{P-Th}, \cite{Joh}, \cite{Eng}. \hfill $\square$
\enddemo

\proclaim{\bf Lemma 9.3} If ${\Cal A, \Cal B\subset \bold {LCTop}}$ are naturally dual subcategories of locally compact spaces (let $D$ be a dualizing object) then if $\Cal A(X,D)$ has compact-open topology (as initilal topology w.r.t. evaluation maps) then initial topology of $|X|\cong \Cal B(\Cal A(X,D),D)$ is compact-open as well.
\endproclaim
\demo{Proof} Evaluation map $ev:\Cal A(X,D)\times |X|\to |D|$ is continuous (since $X$ is locally compact and $\Cal A(X,D)$ has compact-open topology). It implies that initial (point-open) topology on $|X|\cong \Cal B(\Cal A(X,D),D)$ is actually compact-open [by assumption, topology of $|X|$ is initial w.r.t. all maps ${'f'}:|X|@>\sim >>1\times |X|@>f\times 1>>\Cal A(X,D)\times |X|@>ev>>|D|$. It means that topology on $|X|\cong \Cal B(\Cal A(X,D),D)$ is point-open since subbase open sets in point-open and initial topologies are the same $U^{'f'}:=\{ x\in |X|\bigm | {'f'}(x)\in \underset {open}\to U\subset D\} ={'f'}^{-1}(U)$].

We need to show that $\{x\in |X|\bigm | \forall f\in \underset {compact}\to K\subset \Cal A(X,D).\hphantom{ol} 'f'(x)\in \underset {open}\to U\subset D\}=\underset {f\in K}\to \bigcap {{'f'}^{-1}(U)}$ is open in point-open topology on $|X|$.

Take $x\in \underset {f\in K}\to \bigcap {'f'}^{-1}(U)$, then $ev\, (K,x)\subset U$. By continuity of $ev$, $\forall y\in K.\hphantom{o} \exists \underset {open}\to {V_y}\ni ~ y.\mathbreak
\hphantom{o} \exists \underset {open}\to {W_y}\ni x$, s.t. $ev\, (V_y,W_y)\subset U$. $\underset {open \atop covering}\to {\underset {y \in K}\to \bigcup V_y}\supset K$, so, by compactness, $\underset {j=1,\dots,n}\to \bigcup V_{y_j}\supset K$. Therefore, \,$ev\, \bigl (V_{y_j},\underset {j=1,\dots,n}\to \bigcap W_{y_j}\bigr )\subset U$, \,$ev\, \bigl (\underset {j=1,\dots,n}\to \bigcup V_{y_j},\underset {j=1,\dots,n}\to \bigcap W_{y_j}\bigr )\subset U$, \, $ev\, \bigl (K,\underset {j=1,\dots,n}\to \bigcap W_{y_j}\bigr )\subset U$, i.e., $x$ is internal. \hfill $\square$
\enddemo

\proclaim{\bf Corollary} Gelfand-Naimark duality preserves homotopies.
\endproclaim
\demo{Proof} $|A|=\bold {CHTop}(X,\Bbb C)$ has compact-open topology. $|X|=\bold {C^*Alg}(A,\Bbb C)$ has point-open topology, so, by {\bf {Lemma 9.3}} compact-open topology.

Multiplication $c_{A,B,\Bbb C}$ is continuous (since all hom-sets have compact-open topology). Therefore, $F(\bar f)$ is continuous.

[In inverse direction $G:{\bold {CHTop}}\to {\bold {C^*Alg}}$ there is no problem because $\bold {CHTop}(X,\Bbb C)$ has compact-open
topology. See also \cite{Loo}]. \hfill $\square$
\enddemo

\vskip 0.2cm
\centerline{\bf Preserving homotopy relation between homotopies}
\vskip 0.3cm

Let $\Bar {\Bar f}:|I|\times |I|\times |A|\to |B|$ be continuous, s.t. $\Bar {\Bar f}(0,t,a)={\Bar f}_0(t,a)$, \, $\Bar {\Bar f}(1,t,a)={\Bar f}_1(t,a)$.

$$\xymatrix{{\bold {C^*Alg}}(A,B)\times |A|\ar[r]^-{ev} & |B| & {\bold {C^*Alg}}(B,\Bbb C)\times {\bold {C^*Alg}}(A,B)\ar[r]^-{c_{A,B,\Bbb C}} & {\bold {C^*Alg}}(A,\Bbb C)\\
|I|\times |I|\times |A| \ar[u]^{\Bar {\Bar f}^T\times 1_{|A|}} \ar[ur]^{\Bar {\Bar f}} & & {\bold {C^*Alg}}(B,\Bbb C)\times |I|\times |I| \ar[u]^{1\times {\Bar {\Bar f}}^T} \ar[ur]^{F({\Bar {\Bar f}})}\\
1\times |I|\times |A| \ar@/^/[u]^{0\times 1_{|I|\times |A|}} \ar@/_/[u]_{1\times 1_{|I|\times |A|}} & \\
|I|\times |A| \ar[u]_{\sim}^{<!,1_{|I|}>\times 1_{|A|}} \ar@/_1pc/[uuur]^(0.7){{\bar f}_0} \ar@/_1.3pc/@<-1.2ex>[uuur]_(0.7){{\bar f}_1} & & {\bold {C^*Alg}}(B,\Bbb C)\times |I| \ar@/^/[uu]^{1\times ((0\times 1_{|I|})\circ <!,1_{|I|}>)} \ar@/_/[uu]_{1\times ((1\times 1_{|I|})\circ <!,1_{|I|}>)} \ar@/_3.5pc/[uuur]^(0.7){F({\bar f}_0)} \ar@/_3.7pc/@<-1.2ex>[uuur]_(0.7){F({\bar f}_1)}\\
}$$
So, $F({\Bar {\Bar f}})$ is a homotopy from $F({{\Bar f}_0})$ to $F({{\Bar f}_1})$. $F({\Bar {\Bar f}})$ is continuous since $c_{A,B,\Bbb C}$ is continuous in compact-open topology. ${\bold {C^*Alg}}(B,\Bbb C)$ has compact-open topology by {\bf Lemma 9.3}.
\vskip 5mm

\centerline{\bf Preserving unit 2-cells $i_f$}
\vskip -5mm
$$\xymatrix{{\bold {C^*Alg}}(A,B)\times |A|\ar[r]^-{ev} & |B| & {\bold {C^*Alg}}(B,\Bbb C)\times {\bold {C^*Alg}}(A,B)\ar[r]^-{c_{A,B,\Bbb C}} & {\bold {C^*Alg}}(A,\Bbb C)\\
1\times |A| \ar[u]^{f'\times 1} \ar[r]^-{\sim}_-{p_2} & |A| \ar[u]_{f} & {\bold {C^*Alg}}(B,\Bbb C)\times 1 \ar[u]^{1\times f'} \ar[r]^(0.5){\sim }_(0.5){p_1} \ar[ur]^{F(f\circ p_2)} & {\bold {C^*Alg}}(B,\Bbb C) \ar[u]_{-\circ f}\\
|I|\times |A| \ar[u]^{!\times 1} \ar[ur]_{p_2} && {\bold {C^*Alg}}(B,\Bbb C)\times |I| \ar[u]^{1\times !} \ar[ur]_{p_1} \ar@/_0.8pc/[uur]^(0.28){F(i_f)}\\
}$$

So, if $i_f=f\circ p_2\circ (!\times 1_{|A|})=f\circ p_2$, then $F(i_f)=F(f)\circ p_1=i_{F(f)}$.

\vskip 0.2cm
\centerline{\bf Preserving composites $i_g*{\bar f}:|I|\times |A|@>{\bar f}>>|B|@>g>>|C|$}
\centerline{\bf and ${\bar f*i_h:|I|\times |A'|@>{1\times h}>>|I|\times |A|}@>{\bar f}>>|B|$}

$$\xymatrix{{\bold {C^*Alg}}(A,C)\times |A|\ar[r]^-{ev} & |C| &{\bold {C^*Alg}}(C,\Bbb C)\times |I|\ar[dr]^{F(g\circ {\bar f})} \ar[d]|{1\times ({\bold {C^*Alg}}(A,g)\circ f)}   \ar@{->}@/_3.2pc/@<-8ex>[ddd]|{{\bold {C^*Alg}}(g,\Bbb C)\times 1} &\\
{\bold {C^*Alg}}(A,B)\times |A|\ar[u]^{(g\circ -)\times 1} \ar[r]^-{ev} & |B|\ar[u]_{g}& {\bold {C^*Alg}}(C,\Bbb C)\times {\bold {C^*Alg}}(A,C) \ar[r]_-{c_{A,C,\Bbb C}} & {\bold {C^*Alg}}(A,\Bbb C)\\
|I|\times |A| \ar[u]^{f\times 1} \ar[ur]_{\bar f}& &{\bold {C^*Alg}}(B,\Bbb C)\times {\bold {C^*Alg}}(A,B) \ar[r]^-{c{A,B,\Bbb C}} &{\bold {C^*Alg}}(A,\Bbb C) \ar[u]^{1}_{\sim }\\
&&{\bold {C^*Alg}}(B,\Bbb C)\times |I| \ar[u]|{1\times f} \ar[ur]_{F({\bar f})}&\\
}$$
$g\circ {\bar f}$ is a homotopy corresponding to ${\bold {C^*Alg}}(A,g)\circ f$. Outer perimeter of the right diagram commutes because of definition of $F({\bar f})$, $F(g\circ {\bar f})$ and associativity low [if $(s,t)\in {\bold {C^*Alg}}(C,\Bbb C)\times |I|$ then $s\circ (g\circ f(t))=(s\circ g)\circ f(t)$]. So, $F(g\circ {\bar f})=F({\bar f})\circ (F(g)\times 1_{|I|})$, i.e., $F(i_g*{\bar f})=F({\bar f})*i_{F(g)}$.

$$\xymatrix{ {\bold {C^*Alg}}(A,B)\times |A| \ar[r]^-{ev} & |B| & {\bold {C^*Alg}}(B,\Bbb C)\times {\bold {C^*Alg}}(A,B) \ar[r]^-{c_{A,B,\Bbb C}} & {\bold {C^*Alg}}(A,\Bbb C) \ar[dd]|{{\bold {C^*Alg}}(h,\Bbb C)}\\
|I|\times |A| \ar[u]^{f\times 1} \ar[ur]^{\bar f} & & {\bold {C^*Alg}}(B,\Bbb C)\times |I| \ar[u]^{1\times f} \ar[ur]^{F({\bar f})} \ar[d]_{1\times ({\bar f}\circ (1\times h))^T} \ar[dr]^{F({\bar f}\circ (1\times h))}&\\
|I|\times |A'| \ar[u]^{1\times h} \ar[uur]_{{\bar f}\circ (1\times h)} \ar[d]_{({\bar f}\circ (1\times h))^T\times 1} \ar[dr]|{{\bar f}\circ (1\times h)= \atop {ev\circ (f\times 1)\circ (1\times h)= \atop ev\circ (f\times h)}}& &{\bold {C^*Alg}}(B,\Bbb C)\times {\bold {C^*Alg}}(A',B) \ar[r]_-{c_{A',B,\Bbb C}} & {\bold {C^*Alg}}(A',\Bbb C)\\
{\bold {C^*Alg}}(A',B)\times |A'| \ar[r]_-{ev} & |B| \ar[uuu]_(0.4){\sim }^(0.4){1}\\
}$$

Right internal triangle of the right diagram commutes since if $(g,t)\in {\bold {C^*Alg}}(B,\Bbb C)\times |I|$ then ${\bold {C^*Alg}}(h,\Bbb C)\circ c_{A,B,\Bbb C}\circ (1\times f)(g,t)=(g\circ f(t))\circ h=g\circ (f(t)\circ h)=c_{A',B,\Bbb C}(g,f(t)\circ h)=c_{A',B,\Bbb C}(g,({\bar f}\circ (1\times h))^T(t))=c_{A',B,\Bbb C}\circ (1\times ({\bar f}\circ (1\times h))^T)(g,t)$. So, $F({\bar f}*i_h)=F({\bar f}\circ (1\times h))=F(h)\circ F({\bar f})=i_{F(h)}*F({\bar f})$.

\vskip 0.3cm
\centerline{\bf Preserving vertical composites}
\vskip 0.2cm

We need to show if ${\bar f}:{\bar f}\circ i_0\simeq {\bar f} \circ i_1$ and ${\bar g}:{\bar g}\circ i_0\simeq {\bar g} \circ i_1$ are homotopies in ${\bold {C^*Alg}}$ s.t. ${\bar f}\circ i_1={\bar g} \circ i_0$ then $F({\bar g}\odot {\bar f})=F({\bar g})\odot F(\bar f)$.

By definition, vertical composite ${\bar g}\odot {\bar f}$ is

$$\xymatrix@C=3pc{|A|\times |[0,{1\over 2}]| \ar@{^{(}->}@<-2ex>[dr]_{1\times i} & |A|\times |I| \ar[l]_-{1\times \alpha }^-{\sim } \ar[dr]^{\bar f}&\\
& |A|\times |I| \ar@{-->}[r]^{\exists !\, {\bar g}\odot {\bar f}} & |B|\\
|A|\times |[{1\over 2},1]| \ar@{^{(}->}@<1.5ex>[ur]^{1\times j} & |A|\times |I| \ar[l]_-{1\times \beta }^-{\sim } \ar[ur]_{\bar g}&\\
}$$

\hbox{
$$\xymatrix@C=0.5pc{
{\bold {C^*Alg}}(A,B)\times |A|\ar[rr]^-{ev} & & |B|\\
|I|\times |A| \ar[u]|{({\bar g}\odot {\bar f})^T\times 1} \ar[urr]_{{\bar g}\odot {\bar f}} & &\\
|[{1\over 2},1]|\times |A| \ar@{^{(}->}[u]^{j\times 1} & |[0,{1\over 2}]|\times |A| \ar@{^{(}->}@<-1ex>[ul]_{i\times 1}&\\
|I|\times |A| \ar@/^2pc/@<3ex>[uuu]|{f\times 1} \ar[ur]_{\alpha \times 1} \ar[uuurr]_(0.7){\bar f}&&\\
|I|\times |A| \ar@/^2pc/[uu]|(0.25){\beta \times 1} \ar@/^3pc/@<4ex>[uuuu]|{g\times 1} \ar@/_4pc/[uuuurr]_{\bar g}&&\\
}$$

$$\xymatrix@C=1.5pc{{\bold {C^*Alg}}(B,\Bbb C)\times {\bold {C^*Alg}}(A,B) \ar[r]^-{c_{A,B,\Bbb C}}& {\bold {C^*Alg}}(A,\Bbb C)\\
&\\
{\bold {C^*Alg}}(B,\Bbb C)\times |I| \ar[uu]^{1\times ({\bar g}\odot {\bar f})^T} \ar[uur]^{F({\bar g}\odot {\bar f})}& \\
&\\
{\bold {C^*Alg}}(B,\Bbb C)\times |I| \ar[uu]^{1\times (j\circ \beta)} \ar[uuuur]_{F({\bar g})} & {\bold {C^*Alg}}(B,\Bbb C)\times |I| \ar[uul]_(0.4){1\times (i\circ \alpha )} \ar[uuuu]_{F(\bar f)}\\
}$$
}

\vskip 0.35cm
By uniqueness $f\equiv {\bar f}^T=({\bar g}\odot {\bar f})^T\circ i\circ \alpha$, \  $g\equiv {\bar g}^T=({\bar g}\odot {\bar f})^T\circ j\circ \beta$.

\vskip 0.35cm
So, $\cases F({\bar g}\odot {\bar f})\circ (1\times (i\circ \alpha ))=F({\bar f})\\ F({\bar g}\odot {\bar f})\circ (1\times (j\circ \beta ))=F({\bar g})\endcases$. \ It means $F({\bar g}\odot {\bar f})=F({\bar g})\odot F({\bar f})$.

\centerline{\bf Preserving horisontal composites $\xymatrix{A \ar@/^/[r]^{f_0}_(0.4){{\bar f}\, \Downarrow } \ar@/_/[r]_{f_1}& B\ar@/^/[r]^{g_0}_(0.4){{\bar g}\, \Downarrow } \ar@/_/[r]_{g_1}&C}$}

${\bar g}*{\bar f}:=({\bar g}*i_{f_1})\odot (i_{g_0}*{\bar f})\simeq (i_{g_1}*{\bar f})\odot ({\bar g}*i_{f_0})$ (homotopic homotopies). \newline
$F({\bar g}*{\bar f})=F({\bar g}*i_{f_1})\odot F(i_{g_0}*{\bar f})=(i_{F(f_1)}*F({\bar g}))\odot (F({\bar f})*i_{F(g_0)})\simeq F({\bar f})*F({\bar g})$. \

\vskip 0.2cm
Proposition 9.1 completes the proof of Gelfand-Naimark 2-duality 
$\xymatrix{\bold {C^{*}Alg^{op}} \ar@/^/[r]^{F}& \bold {CHTop}\ar@/^/[l]^{G}_{\perp}}$.

\vskip 0.0cm
\proclaim{\bf Proposition 9.1} If $\xymatrix{\bold{C} \ar@/^/[r]^-{F} & \bold{D} \ar@/^/[l]^-{G} }$ are two strict $n$-categories
and two strict $n$-functors in the opposite directions such that the restriction
$\xymatrix{\bold{C}^{\le 1} \ar@/^/[r]^-{F^{\le 1}}_-{\perp } & \bold{D}^{\le 1} \ar@/^/[l]^-{G^{\le 1}} }$
is an adjunction with unit $\eta :1_{\bold{C}^{\le 1}}\to G^{\le 1}F^{\le 1}$ and counit 
$\varepsilon :F^{\le 1}G^{\le 1}\to 1_{\bold{D}^{\le 1}}$ which are still natural transformations for the extension 
(i.e. $\eta :1_{\bold{C}}\to GF$ and $\varepsilon :FG\to 1_{\bold{D}}$ are natural transformations) then the extended situation 
$\xymatrix{\bold{C} \ar@/^/[r]^-{F}_-{\perp } & \bold{D} \ar@/^/[l]^-{G} }$ is a strict adjunction.
\endproclaim 
\demo{Proof} A strict adjunction is completely determined by its 'unit-counit' (proposition 5.3). $\eta :1_{\bold{C}}\to GF$ and
$\varepsilon :FG\to 1_{\bold{D}}$ are natural transformations and satisfy triangle identities
$\varepsilon F\circ _1F\eta =1_F$ and $G\varepsilon \circ _1\eta G=1_G$ (because, e.g. $\varepsilon F=\varepsilon F^{\le 1}$, $1_{F}=1_{F^{\le 1}}$ (set-theoretically), etc.)  \vskip 0.0cm  \hfill   $\square $
\enddemo

\vskip 0.2cm
{\bf Corollary.} Any $1$-adjunction between a category of topological algebras and 
a subcategory of topological spaces is a $2$-adjunction if it can be extended functorially over 2-cells in the way that each instance
of the image of a homotopy is the image of this instance of the preimage-homotopy. 
\demo{Proof} Under given conditions unit and counit of 1-adjunction are automatically natural transformations for the extension. 
E.g., take unit $\eta $.  Naturality square $\vcenter{\xymatrix{A  \ar[r]^-{\eta _A}  \ar[d]_-{f^1} &   GFA  \ar[d]^-{GFf^1}  \\
B   \ar[r]_-{\eta _B}  &     GFB  }}$,
where $f^1:A\times I\to B$ is a homotopy, holds because each instance of it holds (since $\eta $ is a unit of 1-adjunction), i.e.
$\forall t\in I$ $\eta _B\circ f^1(-,t)=GF(f^1(-,t))\circ \eta _A$, it means $\eta _B\circ f^1=GF(f^1)\circ (\eta _A\times I)$, i.e.
$\eta _B*f^1=GF(f^1)*\eta _A$.            \hfill            $\square $
\enddemo

Gelfand-Naimark case is one of the above corollary.   End of proof of Gelfand-Naimark 2-duality.   \vskip 0.0cm  \hfill  $\square $

\vskip 0.2cm
{\bf Remark.} There are 'forgetful' functors $\bold {C^{*}Alg}\to 2\text{-}\bold{Set}$ and 
$\bold {CHTop}\to 2\text{-}\bold{Set}$ (where $2\text{-}\bold{Set}$ is the usual $\bold{Set}$ with just one iso-2-cell for each 
pair of maps with the same domain and codomain) but they are not faithful and forget too much in order $2\text{-}\bold{Set}$ 
could be an underlying category of Gelfand-Naimark 2-duality.     \hfill        $\square $

\vskip 0.3cm
\proclaim{\bf Proposition 9.2} $\bullet $ Gelfand-Naimark 2-duality is concrete over $2\text{-}\bold{Cat}$ ($2\text{-}\bold{Cat}$ is the 
usual 2-category of (small) categories, functors and natural transformations), i.e. $\exists $ (faithful) forgetful functors 
$U:\bold {C^{*}Alg}\to 2\text{-}\bold{Cat}$ and $V:\bold {CHTop}\to 2\text{-}\bold{Cat}$ such that 
$\vcenter{\xymatrix{\bold {C^{*}Alg^{op}}  \ar[r]^-{F} \ar[dr]_-{\bold {C^{*}Alg}(-,\Bbb C)\hskip 0.2cm} &   \bold {CHTop} \ar[d]^-{V} \\
    &          2\text{-}\bold{Cat}  }}$
and \hskip 0.15cm
$\vcenter{\xymatrix{\bold {CHTop^{op}} \ar[r]^-{G^{op}} \ar[dr]_-{\bold {CHTop}(-,\Bbb C)\hskip 0.2cm}  &    \bold {C^{*}Alg}  \ar[d]^-{U} \\
    &          2\text{-}\bold{Cat}  }}$ \hskip 0.1cm
where $U$ and $V$ are composites of inclusion and fundamental groupoid functors  
($U:\bold {C^{*}Alg}\hookrightarrow 2\text{-}\bold{Top}@>2\text{-}\bold{Top}(1,-)>>2\text{-}\bold{Cat}$ and  
$V:\bold {CHTop}\hookrightarrow 2\text{-}\bold{Top}@>2\text{-}\bold{Top}(1,-)>>2\text{-}\bold{Cat}$).  
\vskip 0.1cm\item{$\bullet $} This duality is natural, i.e. lifting of hom-functors $\bold {C^{*}Alg}(-,\Bbb C)$, $\bold {CHTop}(-,\Bbb C)$ along
$V$ and $U$ is initial.                           \hfill    $\square $
\endproclaim

\vskip 0.2cm
{\bf Remark.} 2-duality allows us to transfer (co)homology theories from one side to another.
Under a reasonable assumption that K-theory was determined in a universal way we could get {\bf M. Atiyah theorem} that 
{\it K-groups of commutative $C^*$-algebras and compact Hausdorff spaces coincide}. The problem, however, is that K-groups were determined technically
(not universally). But, there is a theorem by J. Cuntz \cite{Weg} that K-theory is universally determined on a large subcategory of 
$C^*$-algebras.            \hfill        $\square $

\head {\bf 10. Lukacs' extension of Pontryagin duality} \endhead 

The following is a new and recent example of a concrete duality, due to  G.Lukacs \cite{Luk}. The extension is natural with the same dualizing object $\Bbb R/\Bbb Z$, and establishes 
a concrete duality for abelian locally precompact groups.

\vskip 0.2cm
{\bf Definition 10.1.} 
\item{$\bullet $}The set $X$ in a topological group $G$ is called {\bf precompact} if $\forall \, U\ni e$ (neighbourhood of identity)
$\exists \text{ a finite subset } F\subset G$ such that $X\subset FU$.
\item{$\bullet $} The group $G$ is {\bf locally precompact} if it contains a precompact neighbourhood of the identity. \vskip 0.0cm  \hfill  $\square $

\vskip 0.1cm
(Locally) precompact groups are very closed to (locally) compact ones. Namely, their two-sided uniformity completions give (locally)
compact groups, and conversely, dense subgroups of (locally) compact groups are (locally) precompact.

\vskip 0.3cm
\proclaim{\bf Proposition 10.1 (Pontryagin-Lukacs)} There are the following natural dualities
\vskip 0.35cm
\item{} \hskip 3.5cm$\xymatrix{\bold{Ab}^{op} \hskip 0.25cm \ar@/^2ex/[r]_-{\sim } \hskip -0.35cm\ar@{^{(}->}[d] & \bold{CompAb} \ar@/^2ex/[l]_-{\sim } \hskip -0.35cm\ar@{^{(}->}[d] \\
\bold{locCompAb}^{op} \ar@/^2ex/[r]_-{\sim } \hskip -0.35cm\ar@{^{(}->}[d] & \bold{locCompAb} \ar@/^2ex/[l]_-{\sim } \hskip -0.35cm\ar@{^{(}->}[d] \\
\bold{locPreCompAb}^{op} \ar@/^2ex/[r]_-{\sim } & \bold{locCompAb}^{\Rightarrow } \ar@/^2ex/[l]_-{\sim } }$ 
\vskip 0.5cm
\item{} where $\bold{locCompAb}^{\Rightarrow }$ is a category of dense embeddings of locally compact abelian groups into compact abelian
groups (with commutative squares in \, $\bold{locCompAb}$ \, as arrows).    \hfill    $\square $
\endproclaim 

\vskip 0.2cm
{\bf Remarks.} 
\vskip 0.1cm
\item{$\bullet $} The main idea of this extension is that every locally precompact group $G$ can be represented as a 
dense injective $\bold{locCompAb}$-morphism \, $G_d\to compl(G)$, where $G_d$ is the same group with discrete topology, and 
$compl(G)$ is its completion with respect to two-sided uniformity on $G$. After that,  the usual Pontryagin duality is 
used \cite{Luk}.
\item{$\bullet $} The dualizing object in $\bold{locCompAb}^{\Rightarrow }$ is \hskip 0.0cm
$\vcenter{\xymatrix{\Bbb R/\Bbb Z \ar[d]^-{id} \\ \Bbb R/\Bbb Z}}$.   \hfill   $\square $

\head {\bf 11. Differential algebras as a dual to Lie calculus}\endhead 

\vskip 0.2cm
For Lie groups there is an equivalent alternative calculus via exterior differential algebras. For Lie groups of transformations,
it turns out to be more powerful than via Lie algebras. It was developed by E. Cartan and after him by the Russian School in 
Differential Geometry, mainly, by A.M. Vasiliev \cite{Vas0, Vas}.  

\vskip 0.2cm
{\bf Definition 11.1.} 
\item{$\bullet $} The exterior differential algebra $\Lambda \in Ob\, (k\text{-}\Lambda \text{-}\bold{Alg})$, $k=\Bbb C$ or $\Bbb R$, is 
called {\bf linear} if it finitely generated by elements of degree 1 (with possible linear (resolvable) relations between them over $k$).
\item{$\bullet $} The exterior differential algebra $\Lambda \in Ob\, (k\text{-}\Lambda \text{-}\bold{Alg})$, $k=\Bbb C$ or $\Bbb R$, is 
called {\bf quasilinear} if it is finitely generated by elements of degree 0 and 1 with relations between either elements of 
degree 0 or linear relations on elements of degree 1 with coefficients in $\Lambda ^0$.
\item{$\bullet $} A smooth map $f:X\to Y$ is called {\bf quasialgebraic} if there exist quasilinear subalgebras 
$\Lambda _1\hookrightarrow k\text{-}\Lambda (X)$ and $\Lambda _2\hookrightarrow k\text{-}\Lambda (Y)$ such that 
$f^*(\Lambda _2):=k\text{-}\Lambda (f)(\Lambda _2)\hookrightarrow \Lambda _1$.    \hfill    $\square $   

\vskip 0.2cm
Quasialgebraic maps admit an effective description.
All homomorphisms of Lie groups are quasialgebraic.

\vskip 0.2cm
\proclaim{\bf Proposition 11.1} There are equivalences 
$\bold{locLieGrp}\sim \bold{LieAlg}\sim k\text{-}\Lambda \text{-}\bold{ALg}_{lin}^{op}$ (local Lie groups $\sim $ Lie algebras 
$\sim $ (opposite of the category of) linear exterior differential algebras).      \hfill    $\square $
\endproclaim 

\vskip -0.35cm
\proclaim{\bf Lemma 11.1} 
\item{$\bullet $} The functor \, $\bold{Diff}^{op} @>C^{\infty }>> \bold{ComAlg} @>k\text{-}\Lambda >> k\text{-}\Lambda \text{-}\bold{Alg}@>compl>> k\text{-}\Lambda \text{-}\bold{Alg}_{compl}$ is monoidal
with respect to Cartesian product $\times $ in $\bold{Diff}$ and exterior product $\wedge $ in $k\text{-}\Lambda \text{-}\bold{Alg}_{compl}$, where $k\text{-}\Lambda $ 
is a free exterior differential algebra functor over $k$, $compl$ is a smooth (or analytic) completion of exterior differential algebras.  
\item{$\bullet $} Analogously, the functor \, $\bold{LieGrp}^{op} @>k\text{-}\Lambda _{inv}>> k\text{-}\Lambda \text{-}\bold{Alg}$ 
(assigning the algebra of (left)invariant forms) is monoidal
with respect to cartesian product $\times $ in $\bold{LieGrp}$ and exterior product $\wedge $ in $k\text{-}\Lambda \text{-}\bold{Alg}$.  \vskip 0.0cm  \hfill    $\square $
\endproclaim

\vskip -0.35cm
{\bf Remarks.}
\item{$\bullet $} The smooth (analytic) completion functor $compl:k\text{-}\Lambda \bold{Alg}\to k\text{-}\Lambda \text{-}\bold{Alg}_{compl}$ 
is a left adjoint to the inclusion (of the subcategory of smooth (analytic) exterior differential algebras) $k\text{-}\Lambda \text{-}\bold{Alg}_{compl}\hookrightarrow k\text{-}\Lambda \text{-}\bold{Alg}$ 
(it is given essentially by the smooth (analytic) completion of the algebra of coefficients of an exterior differential algebra).
\item{$\bullet $} The exterior product $\wedge $ in $k\text{-}\Lambda \text{-}\bold{Alg}_{compl}$ is bigger than in $k\text{-}\Lambda \text{-}\bold{Alg}$
and is equal to the smooth (analytic) completion of (the usual algebraic) exterior product in $k\text{-}\Lambda \text{-}\bold{Alg}$.   \hfill    $\square $

\vskip 0.2cm
{\bf Definition 11.2.} 
\item{$\bullet $} The exterior differential algebra $\Cal A$ is called {\bf smoothly realizable} if there exists a manifold $Y\in Ob\, \bold{Diff}$
and an embedding $\Cal A\hookrightarrow \Lambda (Y)$. It is {\bf fully} smoothly realizable if $\Cal A^1$ (locally) generates $T^*Y$. 
\item{$\bullet $} A {\bf Geometric triple} is a (locally trivial) fibre bundle $(G\times Y\to X)\in Ar\, \bold{Diff}$, 
equivariant with respect to a (left) action of (Lie group) $G$ on X and 
$\rho :G\times G\times Y\to G\times Y:(g,h,y)\mapsto (gh,y)$ an action on $G\times Y$. A {\bf morphism of geometric triples} is a morphism
of fibre bundles $\vcenter{\xymatrix{G_1\times Y_1 \ar[r]^-{\sigma \times F} \ar[d]_-{\rho _1} & G_2\times Y_2 \ar[d]^-{\rho _2}\\
X_1   \ar[r]_-{f} &   X_2  }}$ where $\sigma :G_1\to G_2$ is a Lie group homomorphism. A geometric triple $\rho :G\times Y\to X$ is 
{\bf local} if $G$ is a local Lie group and $X$ is a local $G$-space (admits a local group of transformations). 
\item{$\bullet $}An {\bf algebraic triple} is an exterior product of two differential algebras $\Cal A\wedge \Cal B$, where $\Cal A$ is linear,  
with a differential ideal $I\subset \Cal A\wedge \Cal B$ generated by elements of degree $1$. A {\bf morphism of algebraic triples}
$(\Cal A_1\wedge \Cal B_1,I_1)\to (\Cal A_2\wedge \Cal B_2,I_2)$ is a differential homomorphism 
$\alpha \wedge \beta :\Cal A_1\wedge \Cal B_1\to \Cal A_2\wedge \Cal B_2$ such that the differential ideal generated by the image 
$\alpha \wedge \beta \, (I_1)$ is $I_2$. An algebraic triple $(\Cal A\wedge \Cal B,I)$ is {\bf smoothly realizable} if 
$\Cal B$ is a smoothly realizable algebra.             \hfill    $\square $

\vskip 0.2cm
Lie groups of transformations are particular cases of geometric triples when $X=Y$ and the projection $\rho :G\times Y\to X$ coincides
with the action of $G$ on $X$.

\vskip 0.0cm
\proclaim{\bf Proposition 11.2} \cite{Vas} The smooth manifold $X$ admits a left action of the finite dimensional Lie group $G$ iff 
there exists a smooth manifold
$Y$, smoothly realizable algebra $\Cal B\hookrightarrow \Lambda (Y)$, and differential ideal 
$I\subset \Lambda _{inv}(G)\wedge \Cal B$ generated by $1$-forms such that the foliation in $G\times Y$ 
determined by $I$ is a (locally trivial) fibre bundle $G\times Y\to X$ with the base $X$.     \hfill  $\square $
\endproclaim

\proclaim{\bf Proposition 11.3} $\bold{locGeomTriple}\sim \bold{realAlgTriple}^{op}$ (local geometric triples $\sim $ 
(opposite to) smoothly realizable algebraic triples).    \hfill    $\square $
\endproclaim 

\vskip 0.2cm
{\bf Remark.} By proposition 11.2. Lie groups of transformations are in duality with a certain full subcategory of 
$\bold{realAlgTriple}$.

\head {\bf Bibliography}\endhead

\vskip 0.2cm
\refstyle{A}
\widestnumber\key{AAAAA}

\ref\key A-H-S
\by J. Adamek, H. Herrlich, G.E. Strecker
\book Abstract and Concrete Categories. The Joy of Cats
\yr online edition, 2004
\endref

\ref\key A-V-L
\by D.V. Alekseevskiy, A.M. Vinogradov, V.V. Lychagin
\book Main Ideas and Concepts of Differential Geometry
\yr 1988
\publ Moscow
\lang Russian
\endref

\ref\key Bel
\by J.L. Bell
\book Toposes and Local Set Theories: An Introduction
\yr 1988
\publ Clarendon Press, Oxford
\endref

\ref\key Bi-Cr
\by R.L. Bishop, R.J. Crittenden
\book Geometry of Manifolds
\publ Academic Press, New York and London
\yr 1964
\endref

\ref\key Bor1
\by F. Borceux
\book Handbook of Categorical Algebra 1. Basic Category Theory
\yr 1994
\publ Cambridge University Press
\endref

\ref\key Bor2
\by F. Borceux
\book Handbook of Categorical Algebra 2. Categories and Structures
\yr 1994
\publ Cambridge University Press
\endref

\ref\key Bor3
\by F. Borceux
\book Handbook of Categorical Algebra 3. Categories of Sheaves
\yr 1994
\publ Cambridge University Press
\endref

\ref\key Bru
\by U. Bruzzo
\book Introduction to Algebraic Topology and Algebraic Geometry
\yr 2002
\publ International School for Advanced Studies, Trieste
\endref

\ref\key BC3G
\by R.L. Bryant, S.S. Chern, R.B. Gardner, H.L. Goldschmidt, P.A. Griffiths
\book Exterior Differential Systems
\yr 1991
\publ Springer-Verlag
\endref

\ref\key Car1
\by E. Cartan
\book Exterior Differential Systems and their Geometric Applications
\yr 1962
\publ Moscow State University
\lang translated into Russian
\endref

\ref\key Car2
\by E. Cartan
\book Moving Frame Method and Theory of Groups of Transformations 
\yr 1963
\publ Moscow State University
\lang translated into Russian
\endref

\ref\key C-C-L
\by S.S. Chern, W.H.Chen, K.S. Lam
\book Lectures on Differential Geometry
\yr 2000
\publ World Scientific
\endref

\ref\key Cl-D
\by D.M. Clark, B.A. Davey
\book Natural Dualities for the Working Algebraist
\yr 1997
\publ Cambridge University Press
\endref

\ref\key C-L
\by E. Cheng, A. Lauda
\book Higher-Dimensional Categories: an illustrated guide book
\publ University of Cambridge
\yr 2004
\endref 

\ref\key D-N-F
\by B.A. Dubrovin, S.P. Novikov, A.T. Fomenko
\book Modern Geometry
\publ Moscow
\yr 1979
\lang Russian
\endref

\ref\key Eng
\by R. Engelking
\book General Topology
\yr 1977
\endref

\ref\key ELOS
\by L.E. Evtushik, U.G. Lumiste, N.M. Ostianu, A.P. Shirokov
\book Differential Geometric Structures on Manifolds
\yr 1979
\publ Moscow
\lang Russian
\endref

\ref\key Fin
\by S.P. Finikov
\book The Method of Cartan's Exterior Forms
\yr 1948
\publ Moscow
\lang Russian
\endref

\ref\key Gar
\by R.B. Gardner
\book The Method of Equivalence and Its Applications
\yr 1989
\publ Society for Industrial and Applied Mathematics
\endref

\ref\key Gol
\by R. Goldblatt
\book Topoi: The Categorial Analysis of Logic
\yr 1984
\publ North-Holland
\endref

\ref\key Hof
\by D. Hofmann
\paper Natural Dualities 
\endref

\ref\key Jac
\by B. Jacobs
\book Categorical Logic and Type Theory
\publ Elsevier, North-Holland
\yr 2001
\endref

\ref\key Joh
\by P.T. Johnstone
\book Stone spaces
\publ Cambridge University Press
\yr 1982
\endref

\ref\key Kob
\by S. Kobayashi
\book Groups of Transformations in Differential Geometry
\publ Moscow
\lang translated into Russian
\yr 1986
\endref

\ref\key Koc
\by J. Kock
\paper Weak Identity Arrows in Higher Categories
\inbook IMRP International Mathematics Research Papers, Volume 2006, Article ID 69163
\pages 1-54
\yr 2006
\endref 

\ref\key Kra
\by I.S. Krasil'shchik
\paper Calculus over Commutative Algebras: a Concise User Guide
\inbook The Diffiety Institute Preprint Series
\yr 1996
\endref 

\ref\key La-Se
\by L.A. Lambe and W.M. Seiler
\book Differential Equations, Spencer Cohomology, and Computing Resolutions
\publ online article
\endref

\ref\key Lap
\by G.F. Laptev
\paper The Main Infinitesimal Structures of Higher Order on a Smooth Manifold
\inbook Proceedings of Geometric Seminar, vol.\, 1
\publaddr Moscow
\yr 1966
\lang Russian
\pages 139-190
\endref

\ref\key Lei
\by T. Leinster
\book Higher Operads, Higher Categories
\publ Cambridge University Press
\yr 2003
\endref

\ref\key Lic
\by A. Lichnerowicz
\book Theorie Globale des Connexions et des Groupes d'Holonomie
\publ Roma, Edizioni Cremonese
\yr 1955
\lang translated into Russian
\endref

\ref\key Loo
\by L.H. Loomis
\book An Introduction to Abstract Harmonic Analysis
\publ D. Van Nostrand Company
\yr 1953
\endref

\ref\key Luk
\by G. Lukacs
\paper Duality Theory of Locally Precompact Groups
\inbook Category Theory Octoberfest
\publ Ottawa University
\yr 2006
\endref 

\ref\key Mac
\by S. MacLane
\book Categories for the Working Mathematician
\publ Springer-Verlag
\yr 1971
\endref

\ref\key M-M
\by S. MacLane, I. Moerdijk
\book Sheaves in Geometry and Logic
\publ Springer-Verlag
\yr 1992
\endref

\ref\key M-T
\by I. Madsen, J. Tornehave
\book From Calculus to Cohomology
\publ Cambridge University Press
\yr 1996
\endref

\ref\key Man
\by U.I. Manin
\book Lectures in Algebraic Geometry. Part 1. Affine Schemes
\publ Moscow State University
\yr 1970
\lang Russian
\endref

\ref\key May
\by J.P. May
\book A Concise Course in Algebraic Topology
\publ The University of Chicago Press
\yr 1999
\endref

\ref\key Nes
\by J. Nestruev
\book Smooth manifolds and observables
\publ Moscow
\yr 2003
\lang Russian
\endref

\ref\key Olv
\by P.J. Olver
\book Equivalence, Invariants and Symmetry
\publ Cambridge University Press
\yr 1995
\endref

\ref\key P-Th
\by H.-E. Porst, W. Tholen
\paper Concrete dualities
\inbook Category Theory at Work\eds H. Herrlich. H.-E. Porst
\publ Heldermann Verlag 
\publaddr Berlin \yr 1991
\pages 111-136
\endref

\ref\key Pos
\by M.M. Postnikov
\book Lectures in Algebraic Topology. Foundations of Homotopy Theory
\publ Moscow
\yr 1984
\lang Russian
\endref

\ref\key Sat
\by H. Sato
\book Algebraic Topology: An Intuitive Approach
\publ American Mathematical Society
\yr 1996
\endref

\ref\key Sei
\by W.M. Seiler
\book Differential Equations, Spencer Cohomology, and Pommaret Bases
\publ Heidelberg University, online lecture notes
\endref

\ref\key Sch
\by J.T. Schwartz
\book Differential Geometry and Topology
\publ New York
\yr 1968
\endref

\ref\key Ste
\by S. Sternberg
\book Lectures on Differential Geometry
\publ Prentice Hall
\yr 1964
\endref

\ref\key Str1
\by T. Streicher
\book Introduction to Category Theory and Categorical Logic
\publ online lecture notes
\yr 2003
\endref

\ref\key S-W
\by R. Sulanke, P. Wintgen
\book Differential Geometry and Fiber Bundles
\yr 1972
\publ Berlin
\lang translated into Russian
\endref

\ref\key Vas0
\by A.M. Vasiliev
\paper Differential Algebra as a Technics of Differential Geometry
\inbook Proceedings of Geometric Seminar, vol.\, 1
\publaddr Moscow
\yr 1966
\lang Russian
\endref

\ref\key Vas
\by A.M. Vasiliev
\book Theory of Differential-Geometric Structures
\yr 1987
\publ Moscow State University
\lang Russian
\endref

\ref\key Ver
\by A.M. Verbovetskii
\book Geometry of Finite Jets and Differential Equations: additional chapters
\publ Moscow Independent University
\yr 1999
\lang Russian
\endref 

\ref\key Vin1
\by A.M. Vinogradov
\paper The Logic Algebra for the Theory of Linear Differential Operators
\publ Dokl. Akad. Nauk SSSR (Soviet Math. Dokl.), vol. 13, No. 4
\yr 1972
\lang Russian
\endref

\ref\key Vin2
\by A.M. Vinogradov
\book Geometry of Nonlinear Differential Equations
\publ Moscow
\yr 1979
\lang Russian
\endref

\ref\key V-K-L
\by A.M. Vinogradov, I.C. Krasil'schik, V.V. Lychagin
\book Introduction to Geometry of Nonliniear Differential Equations
\publ Moscow
\yr 1986
\lang Russian
\endref

\ref\key Weg
\by N.E. Wegge-Olsen
\book K-Theory and $C^*$-Algebras, A Friendly Approach
\publ Oxford University Press
\yr 1993
\endref

\enddocument